\documentclass[11pt]{elsarticle}
\usepackage{color}
\usepackage{graphicx,multirow}
\usepackage{amssymb,mathtools}
\usepackage{epstopdf}
\usepackage{bm}
\usepackage{float}
\usepackage{subcaption}
\usepackage{mathtools}
\usepackage[margin=1in,papersize={8.5in,11in}]{geometry}
\usepackage{times}

\newtheorem{example}{\bf Example}[section]
\newtheorem{remark}{\bf Remark}[section]

\def\vs{\vspace{0.2cm}} 
\def\vse{\vspace{0.1cm}}
\def\L{\mathcal{L}} 
\def\K{\mathcal{K}} 

\journal{arXiv}

\begin{document}
\begin{frontmatter}
\title{Data-driven closures for stochastic dynamical systems}

\author[ucsc]{Catherine Brennan}
\author[ucsc]{Daniele Venturi\corref{correspondingAuthor}}
\address[ucsc]{Department of Applied Mathematics and Statistics \\University of California Santa Cruz\\ Santa Cruz, CA 95064}
\cortext[correspondingAuthor]{Corresponding author}
\ead{venturi@ucsc.edu}

\begin{abstract} 
In this paper we develop a new data-driven closure approximation method to compute the statistical properties of quantities of interest in high-dimensional stochastic dynamical systems. The new method relies on estimating conditional expectations from sample paths or experimental data, and it is independent of the dimension of the underlying phase space. We also address the important question of whether enough useful data is being injected into the reduced-order model governing the quantity of interest. To this end, we develop a new paradigm to measure the information content of data based on the numerical solution of hyperbolic systems of equations. The effectiveness of the proposed new methods is demonstrated in applications to nonlinear dynamical systems and models of systems biology evolving from random initial states.
\end{abstract}

\end{frontmatter}

\section{Introduction}
\label{sec:intro}
High-dimensional stochastic dynamical systems arise naturally in many 
areas of engineering, physical sciences and mathematics. Whether it is a 
physical system being studied in a lab or an equation being solved on a 
computer, the full state of the system is often intractable to handle in all its complexity. Instead, it is 
often desirable to reduce such complexity by moving from a full
model of the dynamics to a reduced-order model that involves only 
the observables of interest. Such observables may represent specific 
features of a stochastic system, e.g., the sensitivity of tumor 
populations to chemo-treatment in stochastic models tumoral cell 
growth \cite{Benzekry,Fiasconaro}, or the viscous dissipation in 
intertial range of fully developed turbulence \cite{McComb,VenturiPhysRep}. 
The dynamics of the observables may be simpler than that 
of the entire system, although the underlying law by which 
they evolve in space and time is often quite complex. 
Nevertheless, approximation of such law 
can in many cases allow us to avoid performing simulation of 
the full system and solve directly for the quantities of interest. 
In this paper, we aim at providing a new general framework 
to compute the probability density function (PDF) of 
such quantities of interest. To this end, we will employ formally 
exact PDF evolution equations and compute the unknown terms 
based on accurate data-driven closure approximations. 
To introduce the methodology, consider the following $N$-dimensional 
dynamical system evolving on a smooth manifold 
$\mathcal{M}\subseteq \mathbb{R}^N$
\begin{equation}
 \label{dyn}
 \begin{cases}
\displaystyle\frac{d{\bm x}}{dt} = \bm {G}(\bm x),\vse\\
{\bm x}(0) = \bm x_0(\omega). 
\end{cases}
\end{equation}
Here $\bm{x_0}(\omega)$ is a random initial state with known 
joint probability density function $p(\bm x_0)$. 
Non-autonomous systems driven by finite-dimensional 
random noise can be always written in the form \eqref{dyn}, 
by augmenting the number of phase variables (see, e.g., 
\cite{Venturi_PRS}). Suppose we are interested in the 
dynamics of a real-valued phase space function 
\begin{equation}
u(\bm x)=\mathcal{M}\rightarrow \mathbb{R}\qquad 
\textrm{(observable).}
\label{observable}
\end{equation}
In models of population biology, this phase space function 
may be represented by the population of a prey species, e.g., 
by the first component of the nonlinear predator-prey system \eqref{dyn}. 
In this case we set $u(\bm x(t))=x_1(t)$. 
The exact dynamics of any observable in 
the form \eqref{observable} can be expressed in 
terms of a semigroup of linear operators 
as \cite{Chorin,Dominy2017}
\begin{equation}
u(\bm x(t,\bm x_0))=\exp\left(t\K(\bm x_0)\right) u(\bm x_0),\qquad \textrm{where}\qquad \K(\bm x_0)=
\sum_{k=1}^N G_k(\bm x_0)\frac{\partial }{\partial x_{0k}}.
\label{koopman}
\end{equation}
Here, $\bm x(t,\bm x_0) $ represents the flow map \cite{Wiggins} 
generated by the system \eqref{dyn}. The linear operator $\exp(t\K)$ 
is known as the Koopman operator \cite{Koopman1931,Li2017}. 
Differentiation of \eqref{koopman} with respect to time yields 
the first-order linear partial differential 
equation (PDE)
\begin{equation}
\frac{\partial u(t,\bm x_0)}{\partial t}=\K(\bm x_0) u(t,\bm x_0).
\label{QOI_PDE}
\end{equation}

{\example $\,$} By setting $u(\bm x(t,\bm x_0))=x_{i}(t,\bm x_0)$ for $i=1,...,N$, 
we obtain the following $N$-dimensional system of linear PDEs  
\begin{equation}
\frac{\partial \bm x(t,\bm x_0)}{\partial t} = \bm G(\bm x_0)\cdot 
\nabla \bm x(t,\bm x_0),
\end{equation}
where the gradient is with respect to the variables $\bm x_0$. This 
system, together with the initial condition  $\bm x(0,\bm x_0)=\bm x_0$, 
allows us to compute the flow map generated by \eqref{dyn}.
\vs

\noindent
The dual of the Koopman operator $\exp(t\K)$ with 
respect to the inner product 
\begin{equation}
\left<f,g\right>=\int_{-\infty}^\infty\cdots \int_{-\infty}^\infty
f(\bm x_0)g(\bm x_0)p(\bm x_0)d\bm x_0
\end{equation}
is know as Frobenious-Perron operator. Such operator can be 
written in the form $\exp(t\L)$, where  
\begin{equation}
\L(\bm x) \phi =  -\nabla \cdot (\bm G(\bm x) \phi(\bm x)).
\end{equation}
The Frobenious-Perron operator pushes forward in time the joint
probability density function of the flow map $\bm x(t,\bm x_0)$, i.e., 
\begin{equation}
p(\bm x,t)=e^{t\L(\bm x)}p(\bm x,0).
\label{frobenious}
\end{equation}

\noindent
Differentiation of \eqref{frobenious} with respect to time yields 
the well-known Liouville transport equation \cite{Sobczyk,Venturi_PRS,Venturi_JCP2013}
\begin{equation}
\frac{\partial p(\bm x,t)}{\partial t}+\nabla\left(\bm G(\bm x)p(\bm x,t \right))=0. 
\label{Liouville}
\end{equation}
Computing the numerical solution of the Liouville equation can be 
quite challenging due to complications with 
high-dimensionality, multiple scales, lack of regularity, positivity 
and conservation properties. From a mathematical viewpoint, 
\eqref{Liouville} is a hyperbolic conservation law in as many variables 
as the dimension of the system \eqref{dyn}.

{\remark $\,$}By using the method of characteristics \cite{Rhee}, it 
is easy to obtain the following formal solution to \eqref{Liouville}
\begin{equation}
p(\bm x,t)=p_0\left(\bm x_0(\bm x,t)\right).
\exp\left(
-\int_0^t \nabla \cdot \bm G\left( \bm x(t,\bm x_0)\right)d\tau 
\right).
\end{equation}
 Here, $p_0(\bm x)=p(\bm x,0)$, while $\bm x_0(\bm x,t)$ denotes the 
 inverse flow map generated by \eqref{dyn}. This expression provides 
 a representation of the Frobenious-Perron semigroup \eqref{frobenious}.
 \vs
 
\noindent
This paper is organized as follows. In Section \ref{sec:reduced_order_PDF} 
we develop reduced order-PDF equations for arbitrary quantities of interest 
\eqref{observable} and discuss their mathematical properties. 
In particular, we discuss the closure problem arising from the 
dimension reduction procedure and relate it with the need of 
computing/estimating conditional expectations. In Section 
\ref{sec:estimating_conditional_expectations} we propose 
a robust procedure to compute conditional expectations 
based on sample trajectories of \eqref{dyn} or 
experimental data. 
Such procedure opens the possibility to compute 
data-driven solutions to reduced-order PDF equations and 
estimate the Mori-Zwanzig memory integrals \cite{Venturi_MZ} 
(Section \ref{sec:data-drivenMZ}). 
In Section \ref{sec:information_content} we develop
a new paradigm measure the information content of data. This 
allows us to infer, in particular, whether we have enough 
data to accurately close the reduced-order PDF equation 
for the quantity of interest. In Section \ref{sec:numerics} 
we demonstrate the effectiveness of the proposed data-driven 
closure approximation method in numerical applications to a 
high-dimensional nonlinear system and a drug resistant 
malaria propagation model.

\section{Reduced-order PDF equations}
\label{sec:reduced_order_PDF}
The Liouville equation \eqref{Liouville} describes the exact 
dynamics of the joint PDF of state variables $\bm x(t)$.
In most cases, however, we are only interested in a 
smaller subset of such  variables, or in the observable 
\eqref{observable} (phase space function).  
The probability density function of such observable can 
be represented as 
\begin{equation}
p(z,t)=\int_{-\infty}^\infty \cdots \int_{-\infty}^\infty 
\delta\left(z-u(\bm x)\right)p(\bm x,t)d\bm x,
\end{equation}
where $\delta(\cdot)$ is the Dirac's delta 
function (see \cite{Khuri,Venturi_MZ,Papoulis}). 
Multiplying the Liouville equation by  
$\delta\left(z-u(\bm x)\right)$ and integrating over 
all phase variables yields
\begin{equation}
\frac{\partial p(z,t)}{\partial t}+ 
\int_{-\infty}^\infty \cdots \int_{-\infty}^\infty 
e^{ia (z-u(\bm x))}
\nabla\cdot 
\left(\bm G(\bm x)p(\bm x,t)\right)\bm x da=0.
\label{reducedOrderPDF1}
\end{equation}
Note that here we employed the Fourier representation 
of the Dirac delta function. 
Equation \eqref{reducedOrderPDF1} governs the dynamics of 
the PDF any $u(\bm x(t))$. In general, it is {\em unclosed} in 
the sense that there are  terms at the right hand side that cannot 
be computed based on $p(z,t)$ alone. If we set 
$\bm u(\bm x(t))=x_k(t)$, i.e., we are interested in the $k$-th 
component of the dynamical system \eqref{dyn}) then 
\eqref{reducedOrderPDF1} reduces to\footnote{By using 
integration by parts and assuming that the joint PDF $p(\bm x,t)$ 
decays to zero sufficiently fast at infinity we obtain
\begin{align}
\int_{-\infty}^\infty \cdots \int_{-\infty}^\infty 
\nabla\cdot  
\left(\bm G (\bm x)p(\bm x,t)\right)
dx_1\dots dx_{k-1}dx_{k+1}\dots dx_N 
=\nonumber\\
\int_{-\infty}^\infty \cdots \int_{-\infty}^\infty 
\frac{\partial }{\partial x_k} 
\left(G_k(\bm x)p(\bm x,t)\right)dx_1\dots dx_{k-1}
dx_{k+1}\dots dx_N.\nonumber
\end{align}
}
\begin{equation}  \label{marginal_general}
\frac{\partial p(x_k,t)}{\partial t } + 
\int_{-\infty}^\infty \cdots \int_{-\infty}^\infty 
\frac{\partial }{\partial x_k} 
\left(G_k(\bm x)p(\bm x,t)\right)dx_1\dots dx_{k-1}dx_{k+1}\dots dx_N = 0.
\end{equation}
The specific form of this equation  depends on underlying dynamical system, i.e., on the nonlinear map $\bm{G}(\bm{x})$. Let us provide 
a simple example.

{\example $\,$} Consider the Kraichnan-Orszag three-mode problem \cite{Orszag,Xiao1} 
\begin{align}
\dot{x}_1 = x_1x_3 \qquad  
\dot{x}_2 = -x_2x_3 \qquad
\dot{x}_3 = -x_1^2 + x_2^2.
\label{KO-dyn}
\end{align}
The associated Liouville equation is 
\begin{align}
\frac{\partial p(\bm{x},t)}{\partial t} = -\frac{\partial }{\partial x_1} \left(x_1x_3p(\bm{x},t)\right)+ \frac{\partial }{\partial x_2} \left( x_2x_3 p(\bm{x},t)\right)  + \frac{\partial }{\partial x_3} \left((x_1^2 - x_2^2)p(\bm{x},t) \right). 
\label{KOLiouville}
\end{align}
Suppose we are interested in the PDF of the 
first component of the system, i.e.,  set $u(\bm x(t))=x_1(t)$ in 
equation \eqref{observable}. 
By integrating \eqref{KOLiouville} with respect to $x_2$ 
and $x_3$ and assuming that $p(\bm x,t)$ decays fast enough 
at infinity, we obtain
\begin{equation}
\frac{\partial p(x_1, t)}{\partial t} = - \frac{\partial }{\partial x_1}
 \int_{-\infty}^{\infty}  x_1x_3 p(x_1, x_3, t)dx_3.
\label{BBGKY_marginal}
\end{equation}
From this equation we see that the evolution 
of $p(x_1, t)$ depends on an integral involving 
$p(x_1, x_3, t)$. In other words,  
to solve an equation of this nature, we must find a 
way to approximate the term involving $p(x_1, x_3, t)$.  
To this end, it is convenient to first  transform the 
integral at the right hand side by using conditional 
probabilities. Specifically, we can write the joint PDF of 
$x_1(t)$ and $x_3(t)$ at time $t$ as
\begin{equation}
p(x_1,x_3,t) = p(x_1,t) p(x_3|x_1,t),
\label{conditional}
\end{equation}
where $p(x_3|x_1,t)$ is the conditional probability density 
of $x_3(t)$ given $x_1(t)$ \cite{Papoulis}. 
A substitution of \eqref{conditional} into \eqref{BBGKY_marginal} 
yields
\begin{equation}
\frac{\partial p(x_1, t)}{\partial t} = - \frac{\partial }{\partial x_1}
\left(x_1 p(x_1,t) \mathbb{E}[x_3(t)|x_1(t)]\right),
\label{BBGKY_marginal1}
\end{equation}
where
\begin{equation}
\mathbb{E}[x_3(t)|x_1(t)]= \int_{-\infty}^{\infty}  x_3 p(x_3 | x_1, t)dx_3
\label{condexpKO}
\end{equation}
is the conditional expectation of $x_3(t)$ given $x_1(t)$. 
As we will see in Section \ref{sec:estimating_conditional_expectations}, 
$\mathbb{E}[x_3(t)|x_1(t)]$ can be estimated from sample 
trajectories of \eqref{KO-dyn}. 
Note that the reduced-order PDF equation \eqref{BBGKY_marginal1} 
is a scalar conservation law where the (compressible) advection 
velocity field is $x_1 \mathbb{E}[x_3(t)|x_1(t)]$.  

\begin{figure}[t]
\centerline{\hspace{0.5cm}odd cumulants\hspace{4.5cm} even 
cumulants}
\centerline{
\includegraphics[height=5.8cm]{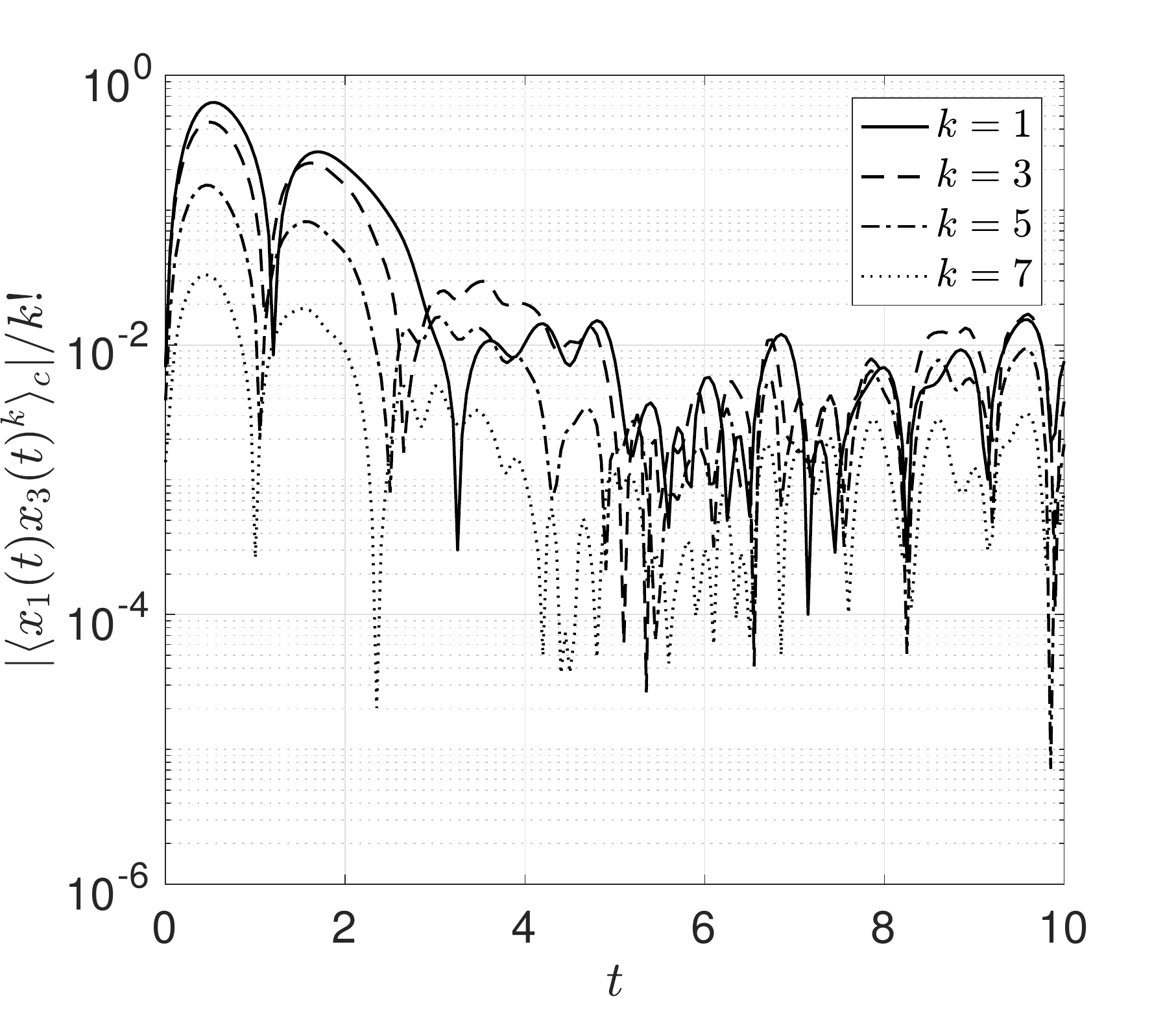}
\includegraphics[height=5.8cm]{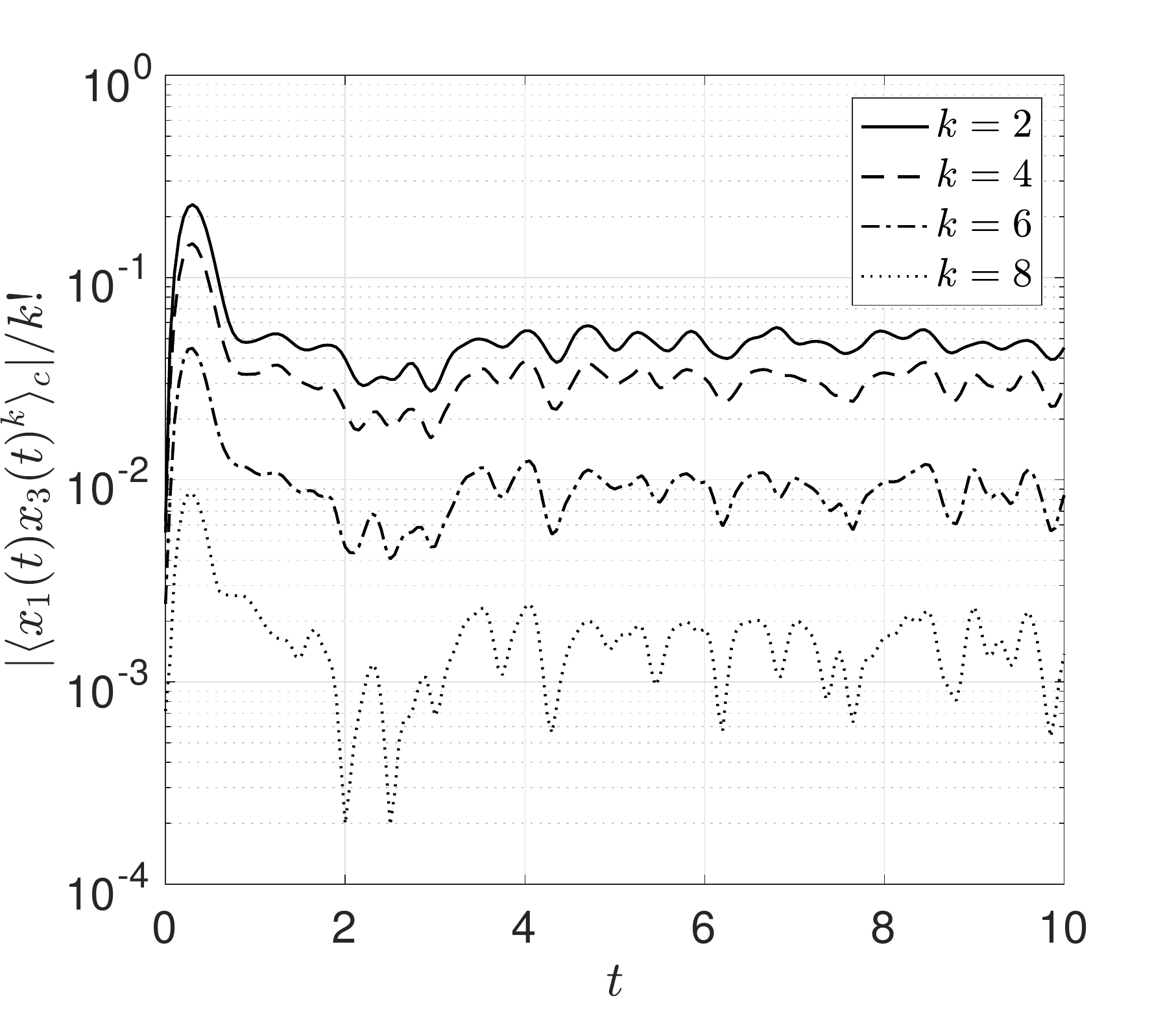}
}
\caption{Kraichnan-Orszag three mode problem. Absolute values 
of the first $8$ rescaled cumulants 
$\langle x_1(t) x_3(t)^k\rangle_c/k!$. The initial 
condition $x_i(0)$ ($i=1,2,3$) in \eqref{KO-dyn} 
is set to be i.i.d. Gaussian with mean and variance $1$. We 
estimated the cumulants numerically by using Monte Carlo ($50000$ 
sample paths) and then taking ensemble averages.
It is seen that the odd cumulants decay slowly with $k$, 
suggesting that the cumulant expansion \eqref{cumulantexp} 
cannot be truncated to a low order.  
This implies that any reasonably accurate approximation 
of the reduced-order equation \eqref{cumulantequation}  
involves high-order derivatives of 
$p(x_1,t)$ with respect to $x_1$.}
\label{fig:KOcumulants}
\end{figure}

{\remark $\,$} It can be shown that the innocent-looking 
 equation \eqref{BBGKY_marginal} is actually a PDE 
 involving derivatives of $p(x_1,t)$ up to order infinity 
 in the phase variable $x_1$. In fact, by using Kubo's cumulant 
expansion \cite{Kubo} of the joint characteristic function of 
$x_3(t)$ and $x_1(t)$, we can prove that  
\begin{equation}
\int_{-\infty}^\infty x_3 p(x_3,x_1,t)dx_3 = \mathbb{E}[x_3(t)] 
p(x_1,t) +\sum_{k=1}^{\infty} (-1)^{k+1}
\frac{\langle x_1(t) x_3(t)^k\rangle_c}{k!}
\frac{\partial^k p(x_1,t)}{\partial x_1^k},
\label{cumulantexp}
\end{equation}
where $\langle x_1(t) x_3(t)^k\rangle_c$ are 
classical cumulant averages\footnote{The cumulant
averages appearing in equation \eqref{cumulantexp} are defined as 
\begin{equation}
\langle x_1(t) x_3(t)^m\rangle_c= \mathbb{E}\left[x_1(t) x_3(t)^m\right]-
\mathbb{E}\left[ x_1(t)\right]\mathbb{E}\left[ x_3(t)^m \right].
\end{equation}}. 
A substitution of \eqref{cumulantexp} into 
\eqref{BBGKY_marginal} yields  the infinite-order PDE
\begin{equation}
\frac{\partial p(x_1, t)}{\partial t} = - \mathbb{E}[x_3(t)]
\frac{\partial \left(x_1 p(x_1,t) \right)}{\partial x_1}+
\sum_{k=1}^{\infty} (-1)^{k+1}\frac{\langle x_1(t) x_3(t)^k\rangle_c}{k!}
\frac{\partial^{k+1} \left(x_1p(x_1,t)\right)}{\partial x_1^{k+1}}.
\label{cumulantequation}
\end{equation}
\noindent
As shown in Figure \ref{fig:KOcumulants}, 
the rescaled cumulants $\langle x_1(t) x_3(t)^k\rangle_c/k!$ 
decay slowly with $k$, suggesting that the cumulant expansion 
\eqref{cumulantexp} cannot be truncated to a low-order. 
This implies that any reasonably accurate approximation of the 
reduced-order PDF equation \eqref{cumulantequation} 
involves high-order derivatives of $p(x_1,t)$ 
with respect to $x_1$.
The data-driven cumulant expansion approach we just described relies 
on computing sample paths of \eqref{KO-dyn}, estimating the 
cumulant averages $\langle x_1(t) x_3(t)^k\rangle_c$ using ensemble 
averaging, and then solving the PDE \eqref{cumulantequation} which 
potentially involves high-order derivatives of $p(x_1,t)$ with respect to 
$x_1$. Clearly this is not practical. 
A more effective approach relies on estimating the conditional 
expectation \eqref{condexpKO} directly from data and then solving 
the hyperbolic conservation law \eqref{BBGKY_marginal1} 
(first-order PDE).\vs

\noindent
More generally, if we are interested in the PDF of $k$-th component 
of the system \eqref{dyn}, then we need to express the right hand 
side of \eqref{marginal_general} in terms of conditional expectations, 
and estimate such expectations from data. If $G_k(\bm x)$ 
is in the form of a sum of separable functions, i.e.,
\begin{equation}
G_k(\bm x)=\sum_{l=1}^r \prod_{j=1}^N f^j_{kl}(x_j), 
\label{CP}
\end{equation}  
then we can explicitly write \eqref{marginal_general} as 
\begin{equation}
\frac{\partial p(x_k,t)}{\partial t}+\frac{\partial }{\partial x_k}\left( p(x_k,t) \sum_{l=1}^r f^k_{kl}(x_k)\mathbb{E}
\left[\left.f^1_{kl}(x_1)... f^{k-1}_{kl}(x_{k-1})f^{k+1}_{kl}(x_{k+1})... f^N_{kl}(x_N) \right| x_k\right]\right)=0.
\label{generalROPDFequation}
\end{equation}  
Well-known examples of nonlinear dynamical systems with 
velocity fields in the form \eqref{CP} are the Kraichnan-Orszag 
system \eqref{KO-dyn}, the 
Lorenz-63 \cite{Viswanath} and the Lorenz-96 \cite{Lorenz96,Karimi} systems,
and the the semi-discrete form of PDEs with quadratic 
nonlinearities such as the Burgers' equation and the Kuramoto-Sivashinsky equation. 
In the next Section, we discuss robust algorithms to compute conditional 
expectations from data, e.g., sample trajectories of \eqref{dyn}. 


\section{Estimating conditional expectations from data}
\label{sec:estimating_conditional_expectations}
Computing conditional expectations from data or sample 
trajectories is a key step in determining accurate closure 
approximations of reduced-order PDF 
equations. A major challenge to fitting a conditional 
expectation is ensuring accuracy and stability. More 
importantly, the estimator must be flexible and 
effective for a wide range of numerical applications.
Let us briefly recall what conditional expectations are 
and, more importantly, how to compute them 
efficiently based on sample paths. To this end, 
let us provide a simple example.
\begin{figure}[t]
\centerline{\hspace{0.2cm}(a)\hspace{7cm} (b)}
	\centering
	\begin{subfigure}[b]{0.45\textwidth}
		\includegraphics[width=\textwidth]{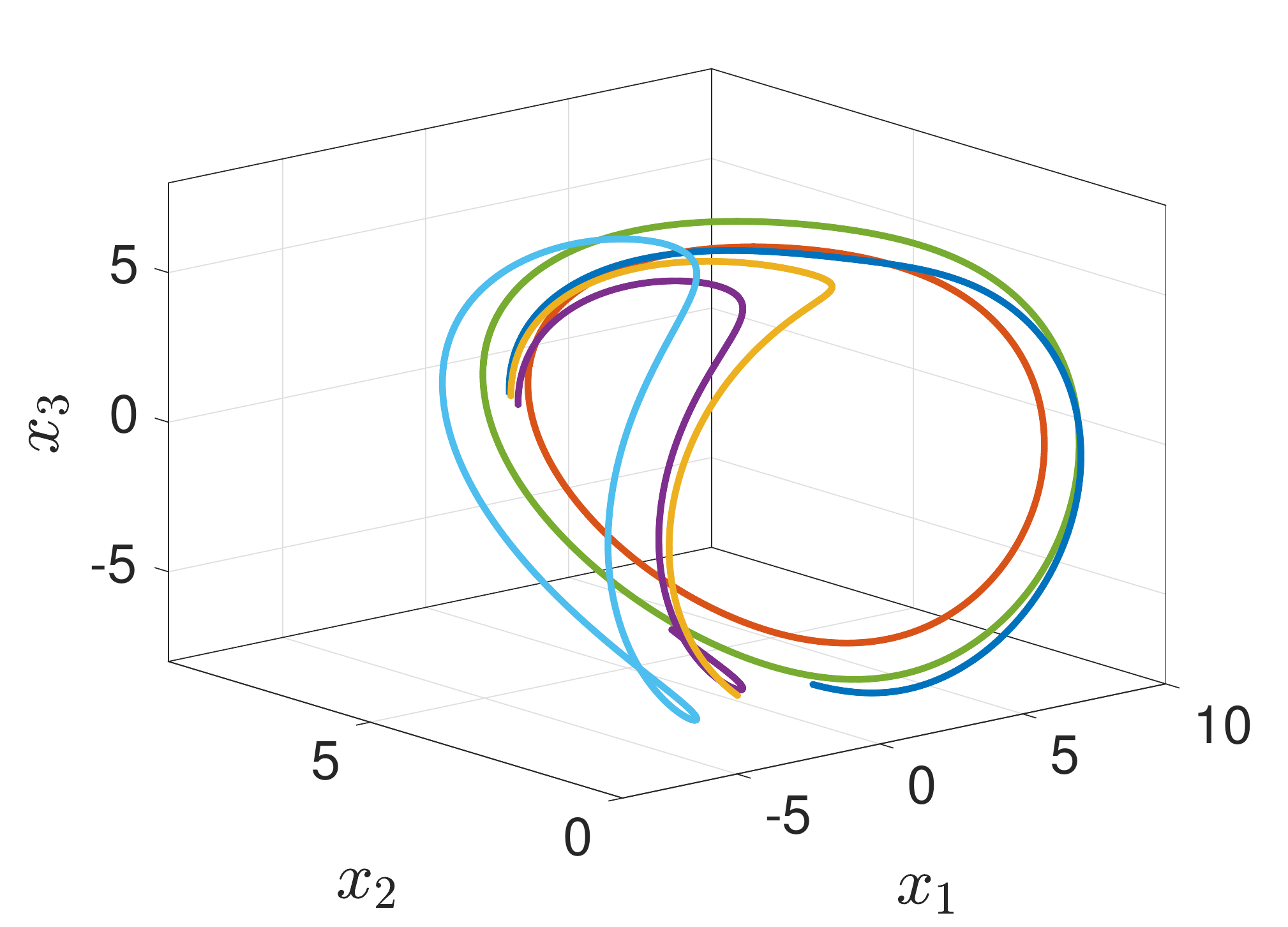}
	\end{subfigure}
	~
	\begin{subfigure}[b]{0.44\textwidth}
		\includegraphics[width=\textwidth]{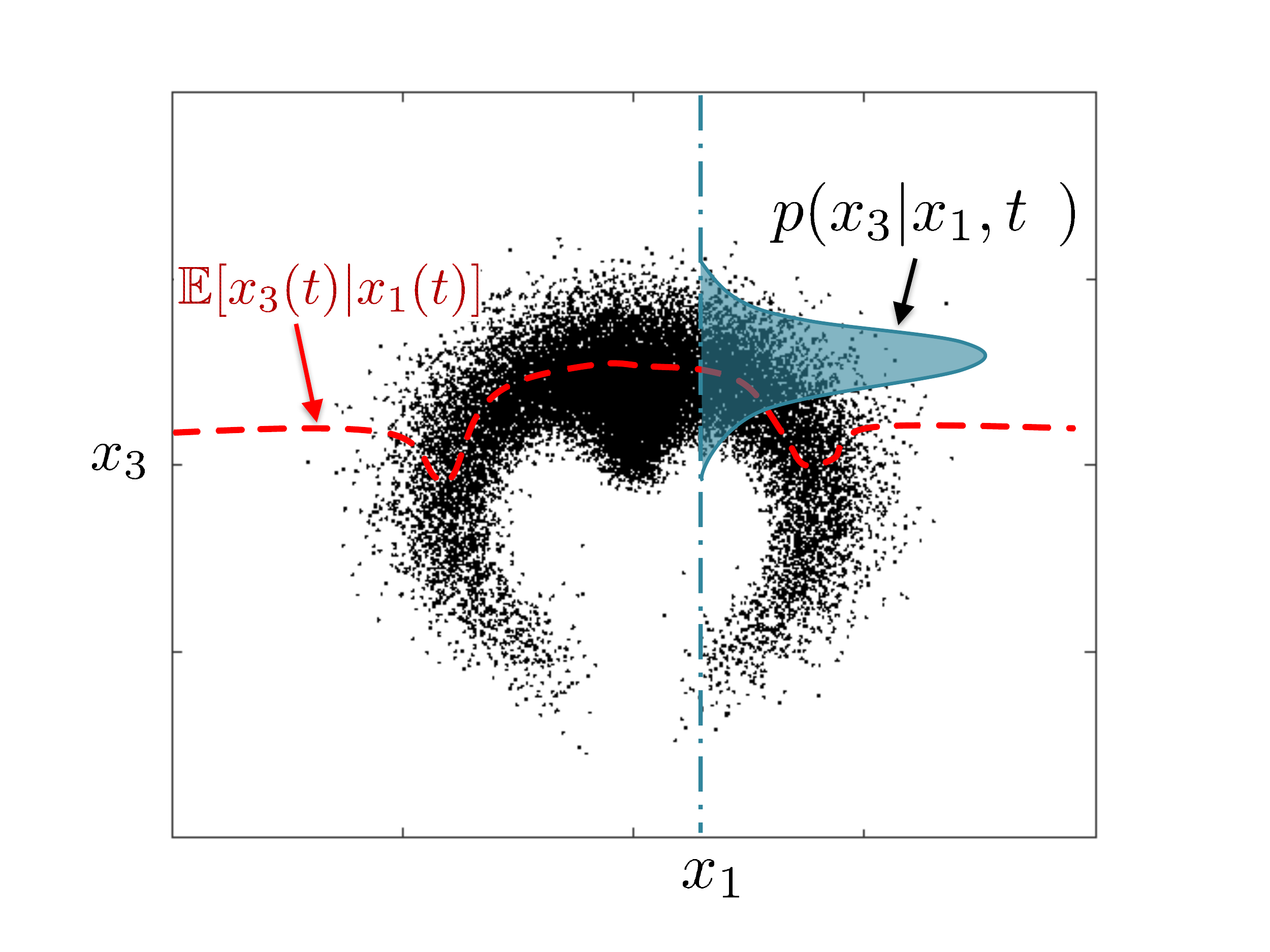}
	\end{subfigure}
\caption{Kraichnan-Orszag three mode problem. (a) Sample 
trajectories of \eqref{KO-dyn} corresponding to random samples 
the initial condition $\bm x_0$. (b) Solution samples projected 
into the plane $(x_1,x_3)$ at time $t$. For each value of $x_1$, the 
conditional PDF $p(x_3|x_1,t)$ can be estimated based on 
points sitting on or lying nearby the vertical dashed line. The 
conditional expectation $\mathbb{E}[x_3(t)|x_1(t)]$ is 
the mean of such conditional PDF.}
\label{fig:conditional}
\end{figure}
\begin{figure}[t]
\centerline{\footnotesize\hspace{0.6cm}$10$ samples \hspace{3.6cm} $100$ samples \hspace{3.4cm} $1000$ samples}
	\centerline{
	\rotatebox{90}{\hspace{1.3cm} \footnotesize moving average}
		\includegraphics[width=0.3\textwidth]{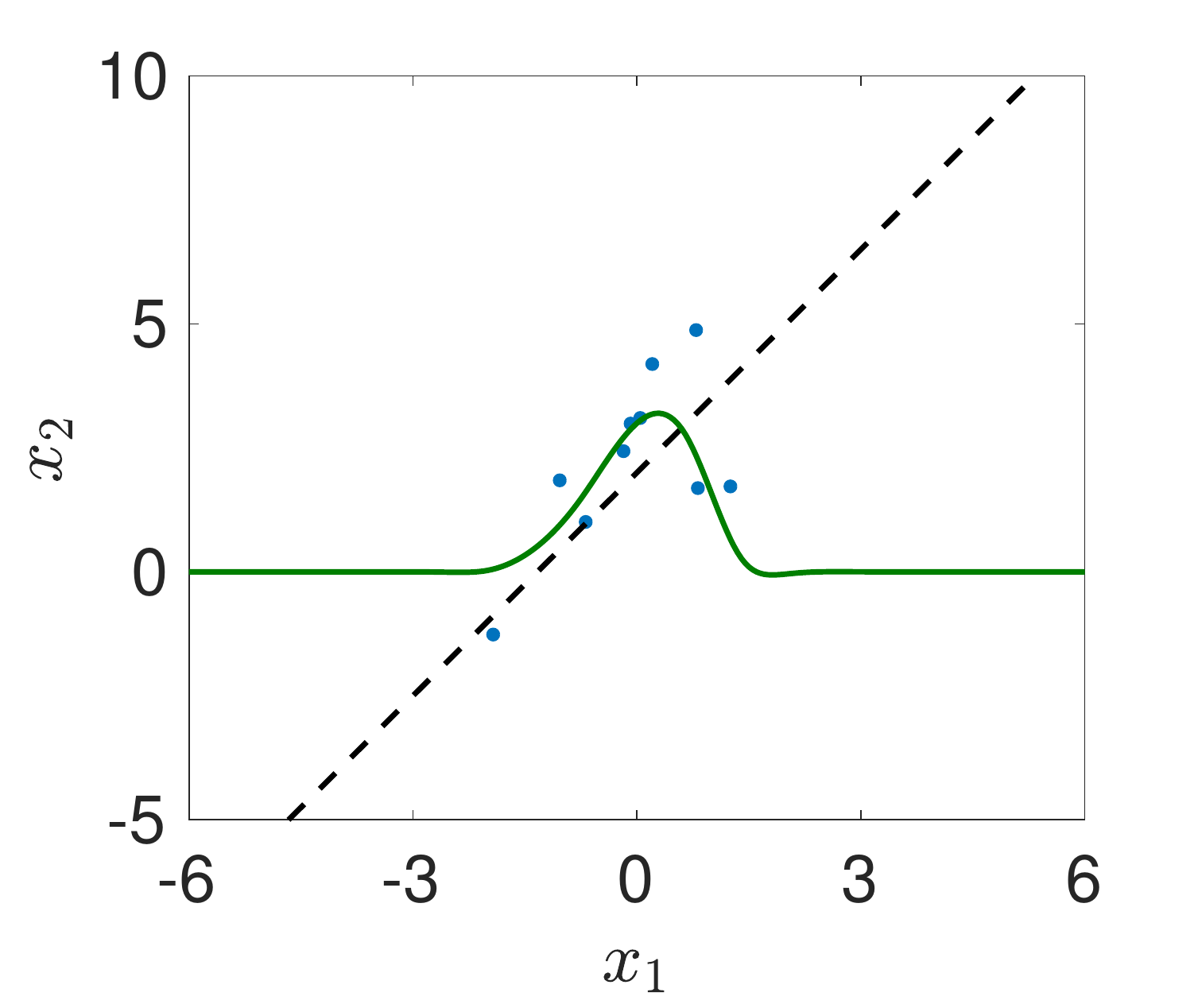}
		\includegraphics[width=0.3\textwidth]{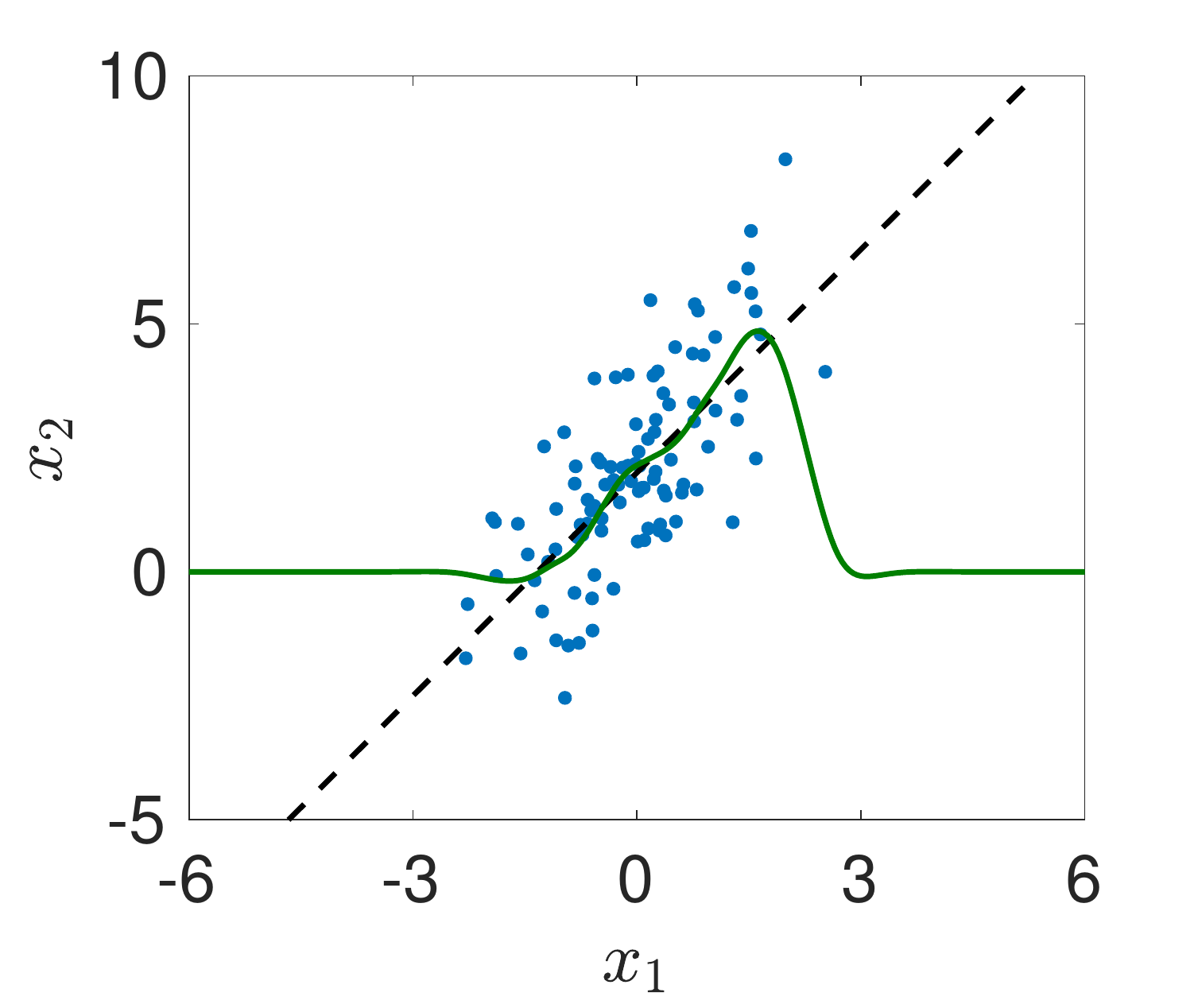}
		\includegraphics[width=0.3\textwidth]{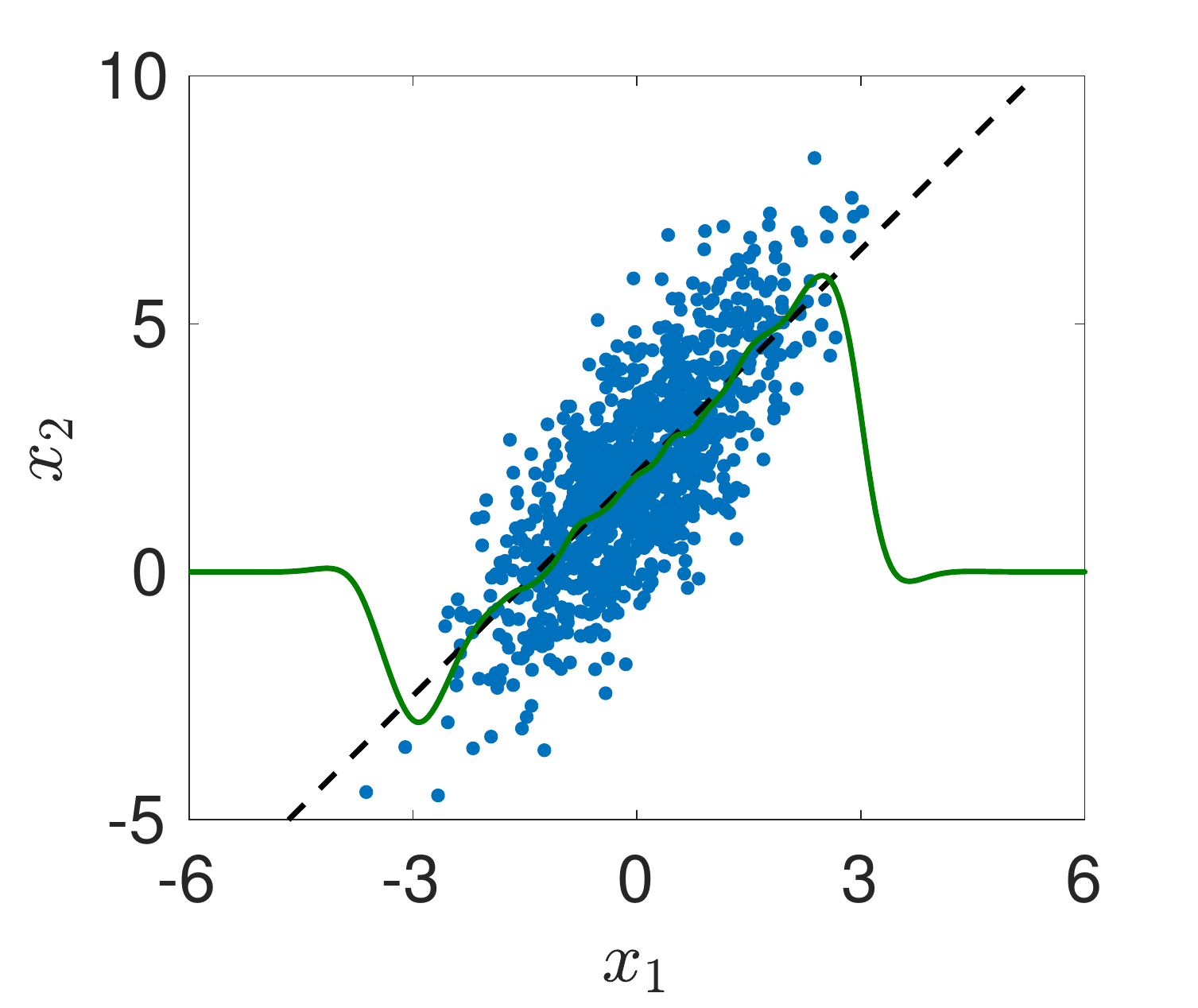}	
}
	
	\centerline{
	\rotatebox{90}{\hspace{1.3cm} \footnotesize smoothing splines}
       \includegraphics[width=0.3\textwidth]{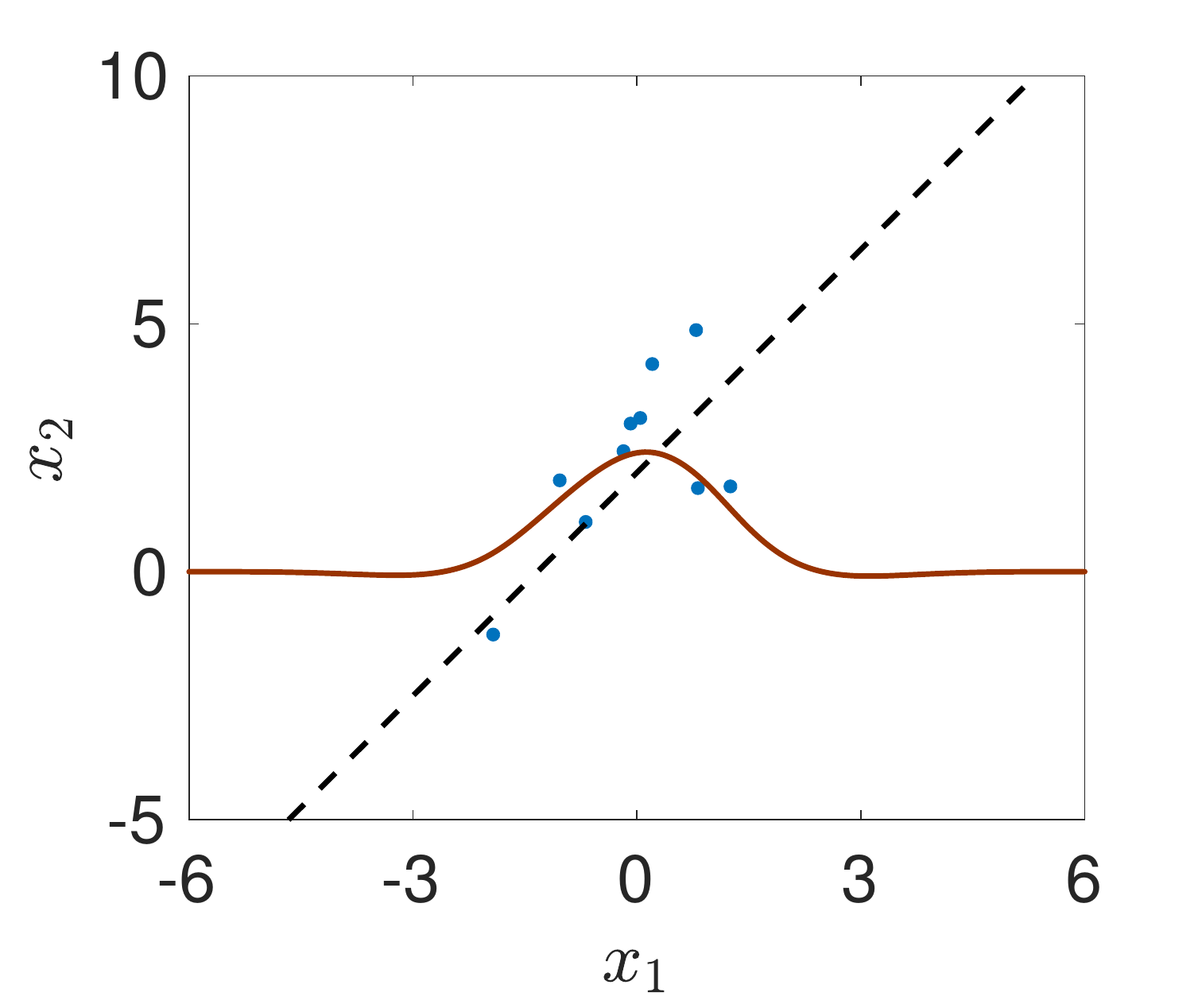}
	    \includegraphics[width=0.3\textwidth]{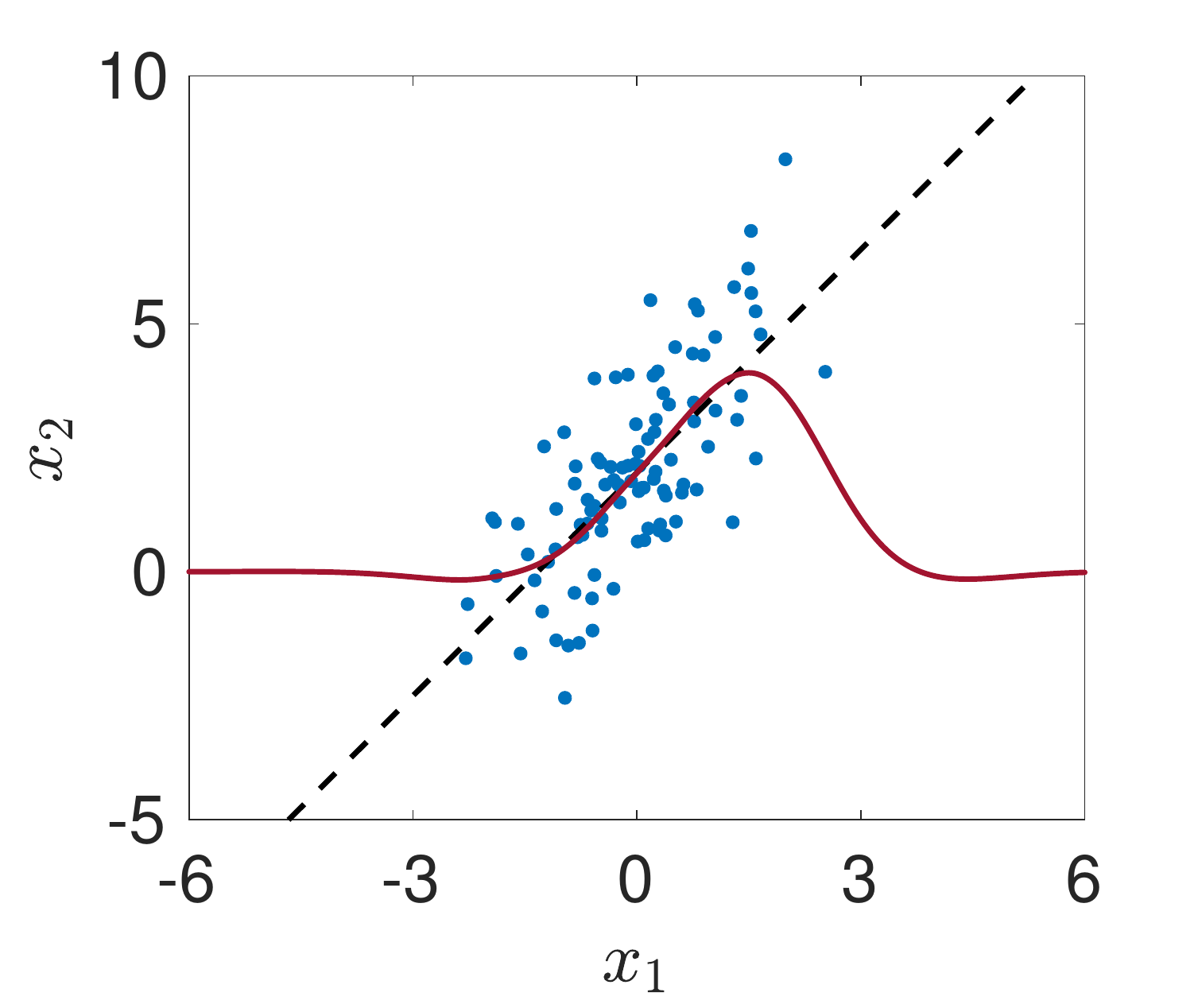}
	    \includegraphics[width=0.3\textwidth]{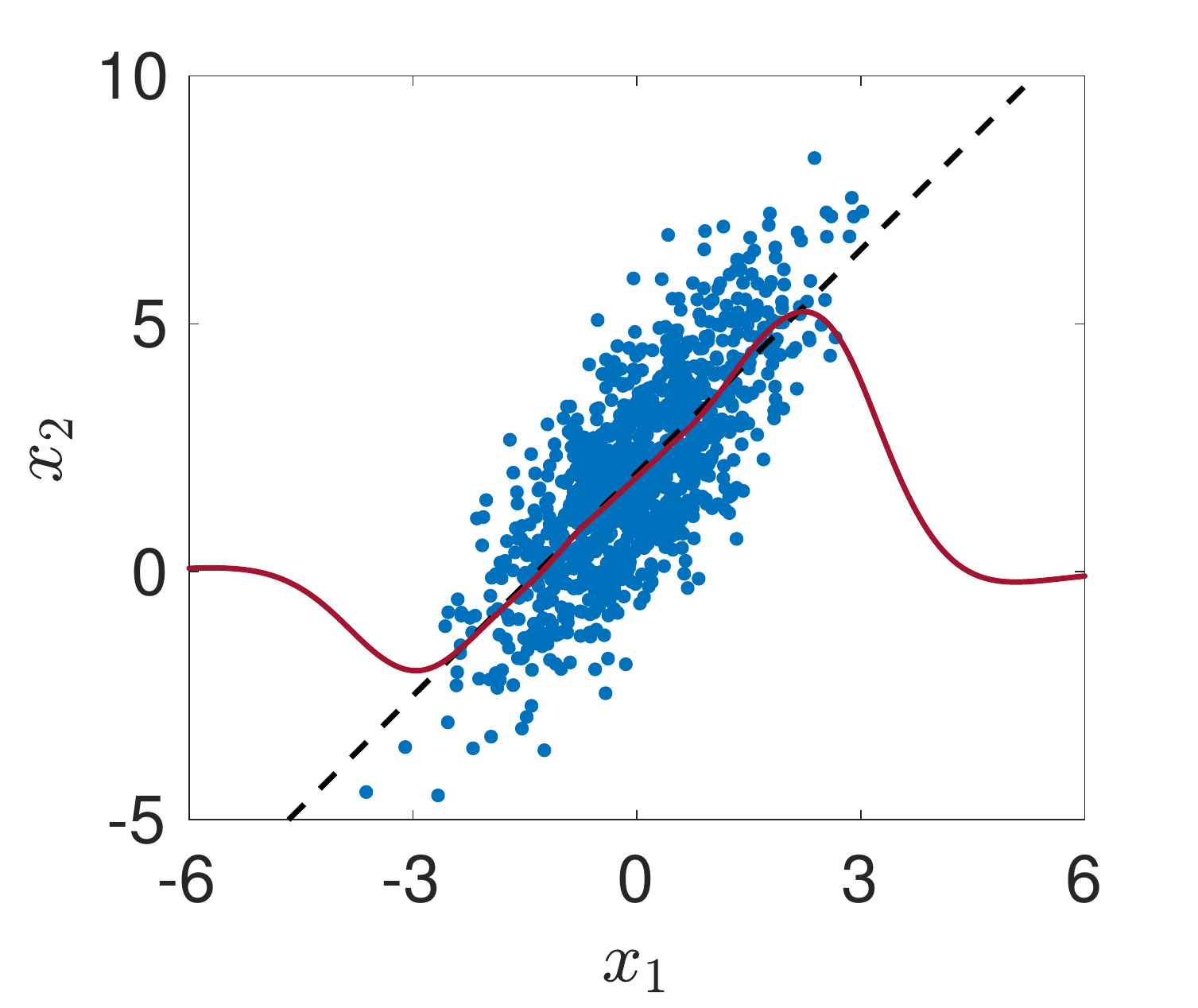}
}
\caption{Numerical estimation of the conditional expectation 
\eqref{condexp1} for different number of samples 
of \eqref{jointPDFG}. Shown are results obtained with moving 
averages (first row) and cubic smoothing splines (second row). It is seen that 
both methods converge to the correct conditional expectation in 
the active region as we increase the number of samples.
}
\label{fig:estimators}
\end{figure}
{\example$\,$} Consider the random 
processes $x_1(t)$ and $x_3(t)$ defined by the dynamical 
system \eqref{KO-dyn} with random initial state. 
The conditional expectation of $x_3(t)$ given $x_1(t)$ 
is defined in \eqref{condexpKO}. The geometric meaning 
of such conditional expectation is illustrated in 
Figure \ref{fig:conditional}. We first compute sample 
trajectories of \eqref{KO-dyn} -- see 
Figure \ref{fig:conditional}(a) -- by sampling the initial condition 
and evolving it in time . We then project the solution samples we 
obtain at time $t$ into the plane $(x_1,x_3)$, to obtain the scatter plot 
in Figure \ref{fig:conditional}(b). For each value of $x_1$, the 
conditional PDF $p(x_3|x_1,t)$ can be estimated based on all 
points sitting on or lying nearby the vertical dashed line. The 
conditional expectation $\mathbb{E}[x_3(t)|x_1(t)]$ is 
the mean of such conditional PDF.
\vs

\noindent Hereafter we present two different approaches 
to estimate conditional expectations from data based on 
moving averages and smoothing splines. 
The moving average estimate is obtained by first sorting the data into 
bins and then computing the average within each bin. With such 
averages available, we can construct a smooth interpolant 
using the average value within each bin. Some factors that affect 
the bin average approximation are the bin size (the number of samples in 
each bin) and the interpolation method used in the final step.
Another approach to estimate conditional expectations uses 
smoothing splines.  This approach seeks to minimize a penalized 
sum of squares. A smoothing parameter determines the balance 
between smoothness and goodness-of-fit in the least-squares 
sense \cite{Craven}. The choice of smoothing parameter is 
critical to the accuracy of the results. Specifying 
the smoothing parameter a priori is generally yields poor 
estimates \cite{Pope}. Instead, cross-validation and 
maximum likelihood estimators can guide the choice the optimal 
smoothing value for the data set \cite{Wahba}. Such methods 
can be computationally intensive, especially when the spline 
estimate is performed at each time step. Other techniques 
to compute conditional expectations can leverage on recent 
developments on deep learning \cite{Graupe2016}.

In Figure \ref{fig:estimators} we compare the performance 
of the moving average and smoothing splines approaches 
in approximating the conditional expectation of two jointly 
Gaussian random variables. Specifically, we consider the joint distribution
\begin{equation}
p(x_1,x_2)=\frac{1}{2\pi\sigma_1\sigma_2\sqrt{1-\rho^2}}
\exp\left(-\frac{1}{2(1-\rho^2)}\left[\frac{(x_1-\mu_1)}{\sigma_1^2}\frac{(x_2-\mu_2)}{\sigma_2^2}-
\frac{2\rho(x_1-\mu_1)(x_2-\mu_2)}{\sigma_1\sigma_2)}
\right] \right)
\label{jointPDFG}
\end{equation}
with parameters $\rho=3/4$, $\mu_1=0$, $\mu_2=2$, 
$\sigma_1=1$ $\sigma_2=2$.
As is well known \cite{Papoulis}, the conditional 
expectation of $x_2$ given $x_1$ can be expressed 
as\footnote{Given two random variables with joint 
PDF $p(x_1,x_2)$,  the conditional expectation of $x_2$ given $x_1$ is defined as 
\begin{equation}
\mathbb{E}[x_2|x_1]=\int_{-\infty}^\infty x_2 p(x_2|x_1)dx_2=\frac{1}{p(x_1)}
\int_{-\infty}^\infty x_2 p(x_1,x_2) dx_2, 
\end{equation}
where $p(x_1)$ is the marginal of $p(x_1,x_2)$ with 
respect to $x_2$. 
}
\begin{align}
\mathbb{E}[x_2|x_1]=\mu_2+\rho\frac{\sigma_2}{\sigma_1}(x_1-\mu_1)
=2+\frac{3}{2}x_1.
\label{condexp1}
\end{align}
Such conditional expectation is plotted in Figure \ref{fig:estimators} 
(dashed line), together with the plots of the conditional average
estimates we obtain with the moving average and the smoothing 
spline approaches for different numbers of samples. 
It is seen that both methods converge to the correct conditional 
expectation as we increase the number of samples. Note, however, that 
convergence is achieved in regions where the PDF \eqref{jointPDFG} is 
not small (see the subsequent Remark  \ref{rmk:3.2}). 
Both estimators require setting suitable parameters to compute 
expectations, e.g., the width of the moving average window
in the moving average approach, or the smoothing parameter in 
the cubic spline approximant.

{\remark $\,$} If the joint PDF of $x_1$ and $x_2$ is not 
compactly supported, then the conditional expectation is 
defined in the whole real line.  It is computationally 
challenging to estimate the expectation \eqref{condexp1} 
in regions where the joint PDF is very small \cite{Casella-Berger}. 
At the same time, if we are not interested in rare events (i.e., tails of 
probability densities), then resolving the dynamics in such 
regions of small probability is not needed. This means that if we 
have available a sufficient number of sample trajectories
\footnote{In Section \ref{sec:information_content}, we will
quantify what a sufficient number of trajectories is, and propose
a new way to measure the information content of data 
based on PDEs.} then we can identify the active regions where 
the dynamics are happening with high probability, and 
approximate the conditional expectation only within such 
regions \cite{Chorin2009}. Outside the active regions, 
we set the expectation equal to zero. 

{\remark $\,$ \label{rmk:3.2}} If the joint PDF of $x_1$ and $x_2$ 
is compactly supported, e.g. uniform in the square $[0,1]^2$, 
then the conditional expectation is undefined outside the 
support of the joint PDF. 
This means, in principle, that we are not allowed to set any 
value for the conditional expectation outside the domain where 
it exists. However, a quick look at the structure of the 
reduced-order PDF equations we are considering in this paper, 
e.g., equation \eqref{generalROPDFequation}, suggests 
that the conditional expectation plays the role of a velocity field 
advecting the reduced-order PDF. Therefore, setting such 
velocity vector equal to zero in the regions where 
the reduced order PDF is very small or even undefined, 
does not affect the PDF propagation process. 
On the other hand, setting the conditional expectation 
equal to zero in low- or zero-probability regions greatly 
simplifies the mathematical discretization of PDEs in the 
form \eqref{generalROPDFequation}.
\vs

\begin{figure}[t]
\centerline{
    \rotatebox{90}{\hspace{0.7cm} \footnotesize Kraichnan-Orszag 
    \cite{Orszag}}\hspace{0.5cm}
		\includegraphics[width=0.3\textwidth]{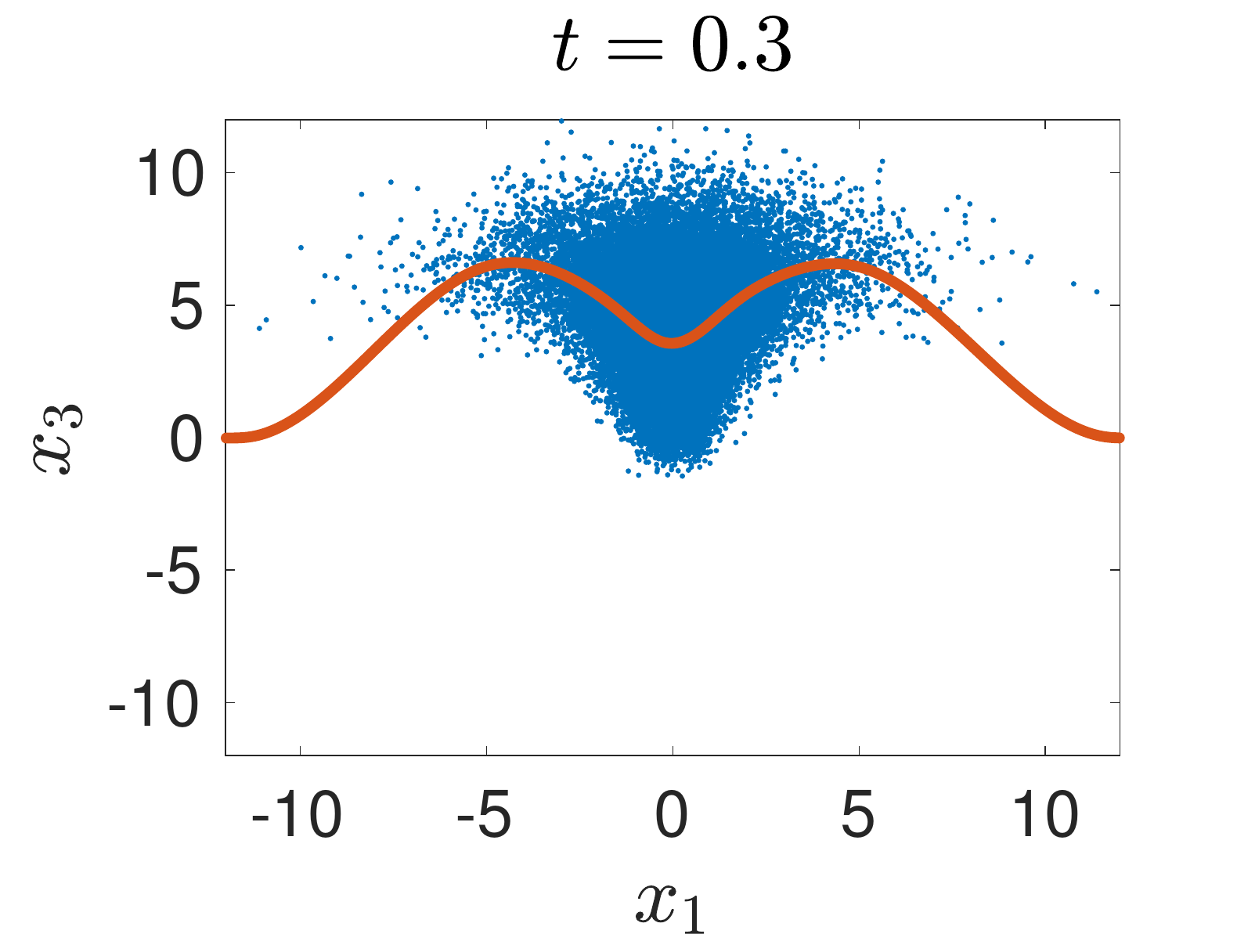}
		\includegraphics[width=0.3\textwidth]{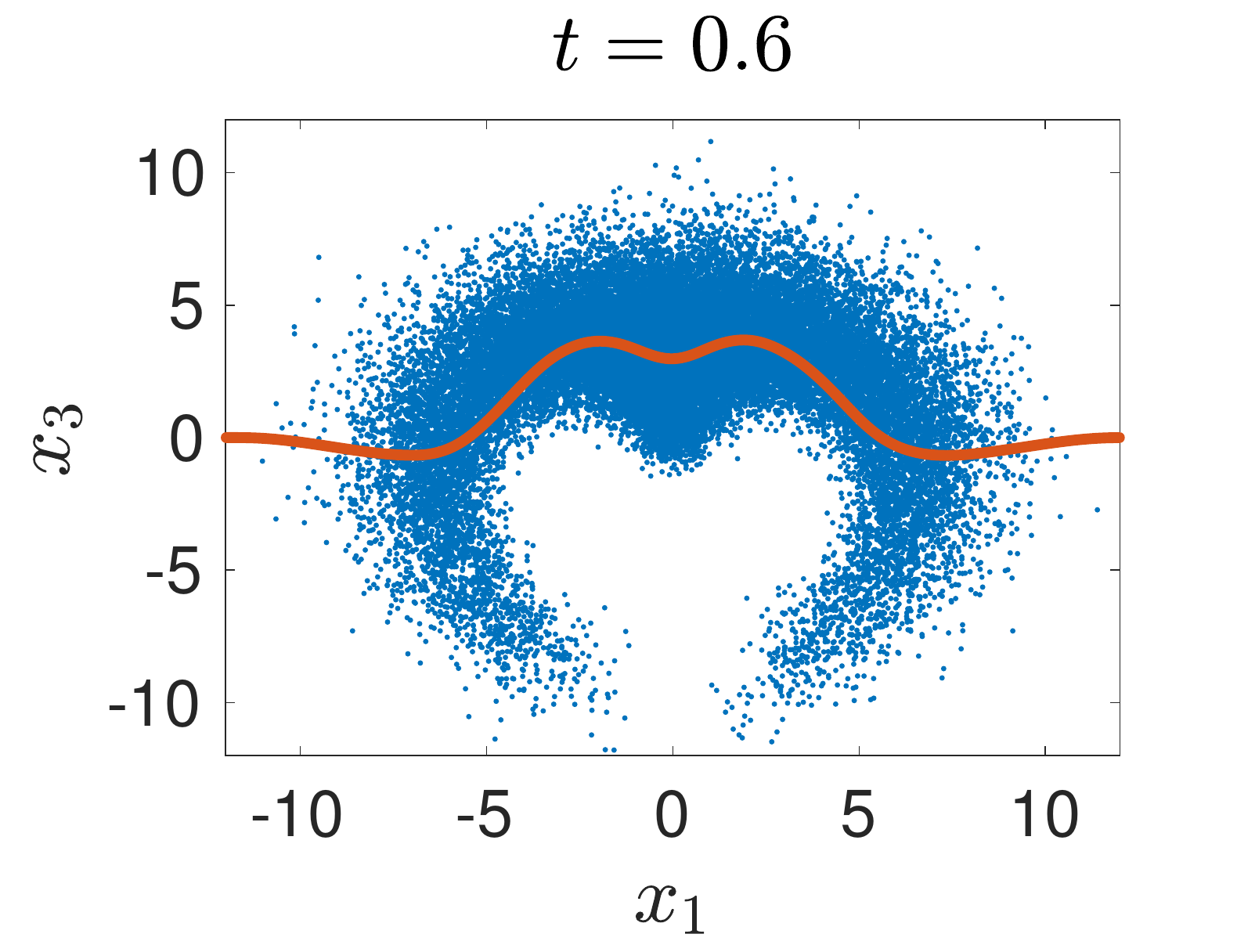}
		\includegraphics[width=0.3\textwidth]{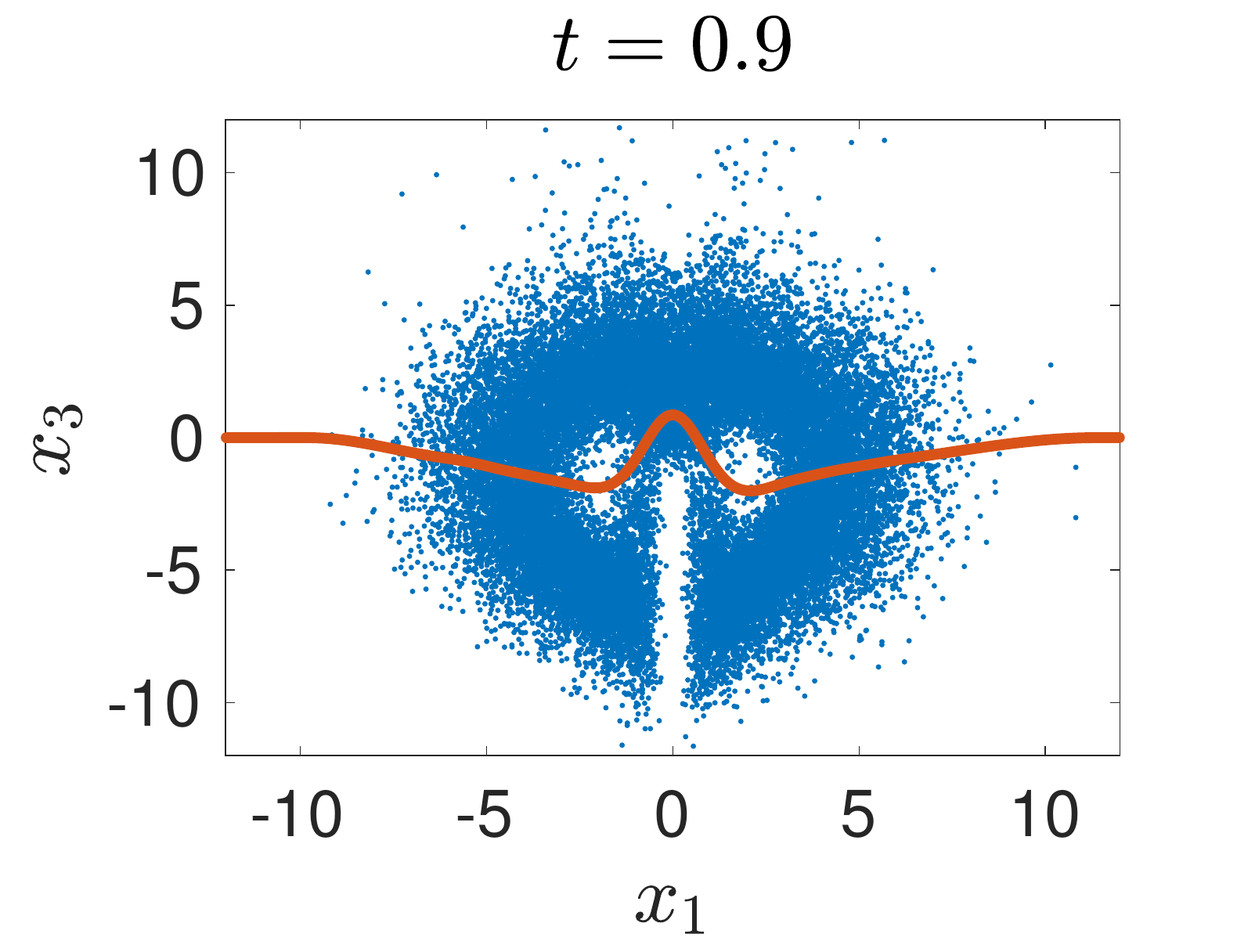}
	}
	
\centerline{ 
    \rotatebox{90}{\hspace{0.5cm} \footnotesize High-dimensional system}
        \hspace{0.5cm}
	  \includegraphics[width=0.3\textwidth]{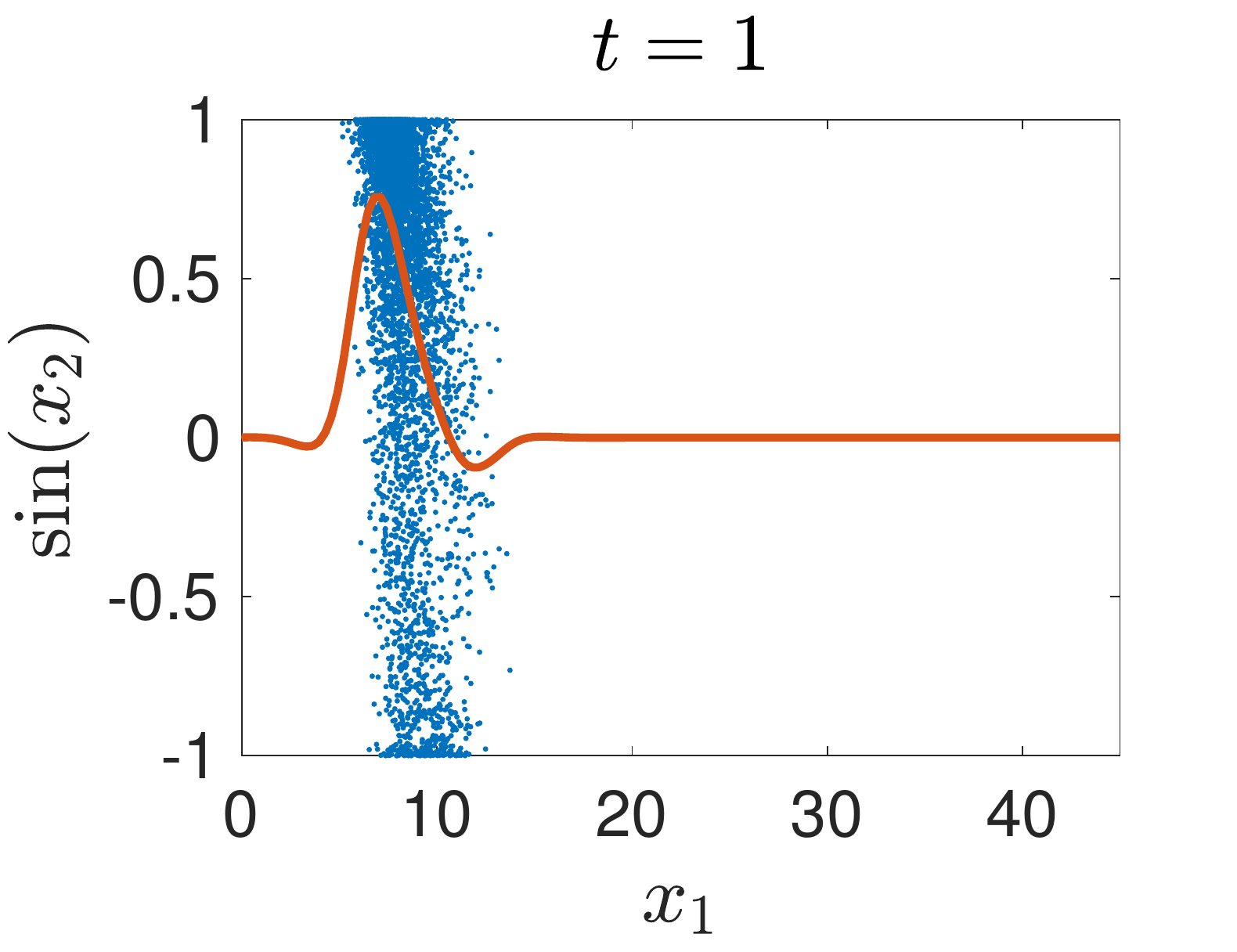}
	  \includegraphics[width=0.3\textwidth]{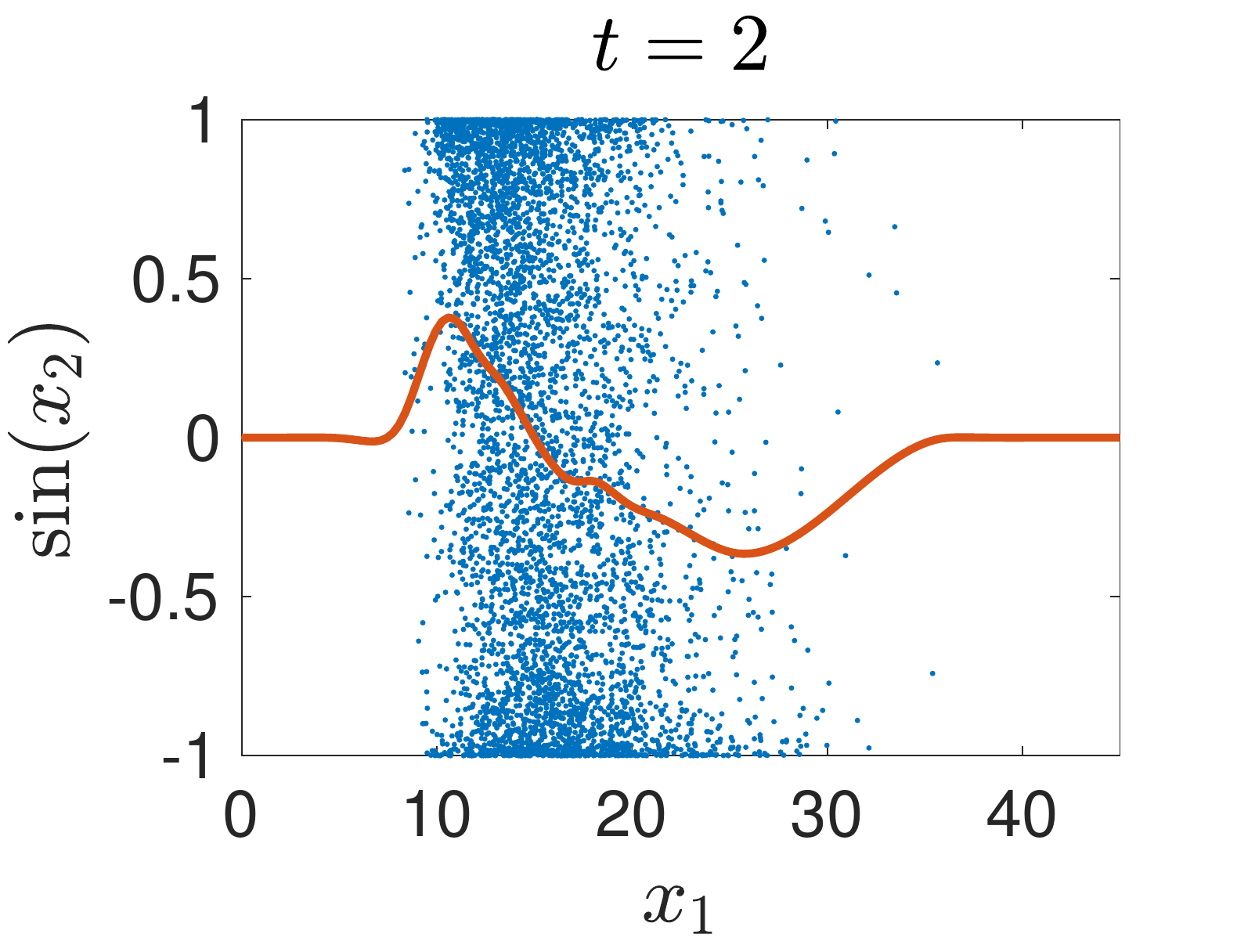}
	  \includegraphics[width=0.3\textwidth]{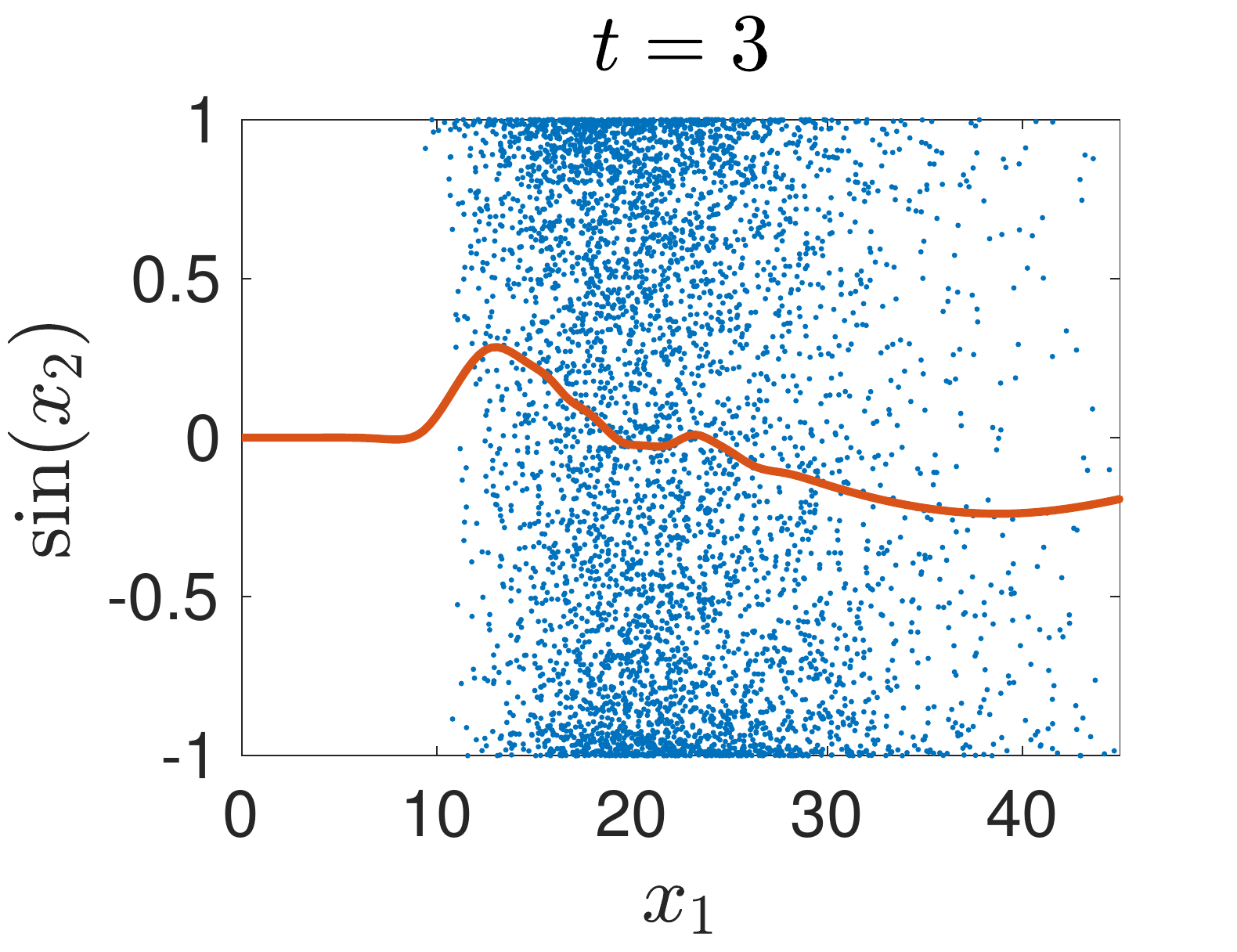}
}

\centerline{
    \rotatebox{90}{\hspace{0.8cm} \footnotesize Malaria model \cite{Meara}}
    \hspace{0.5cm}
	\includegraphics[width=0.3\textwidth]{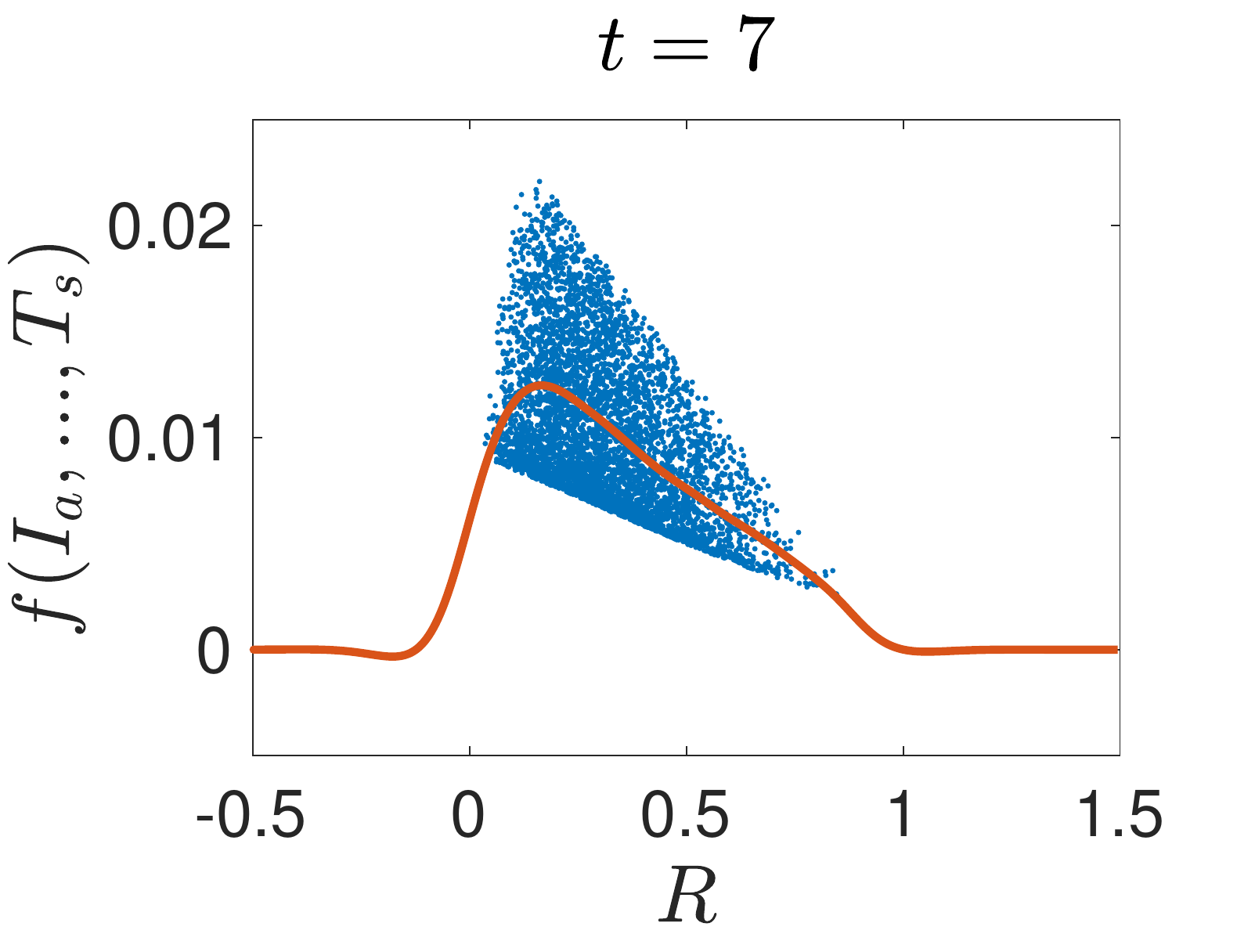}
	\includegraphics[width=0.3\textwidth]{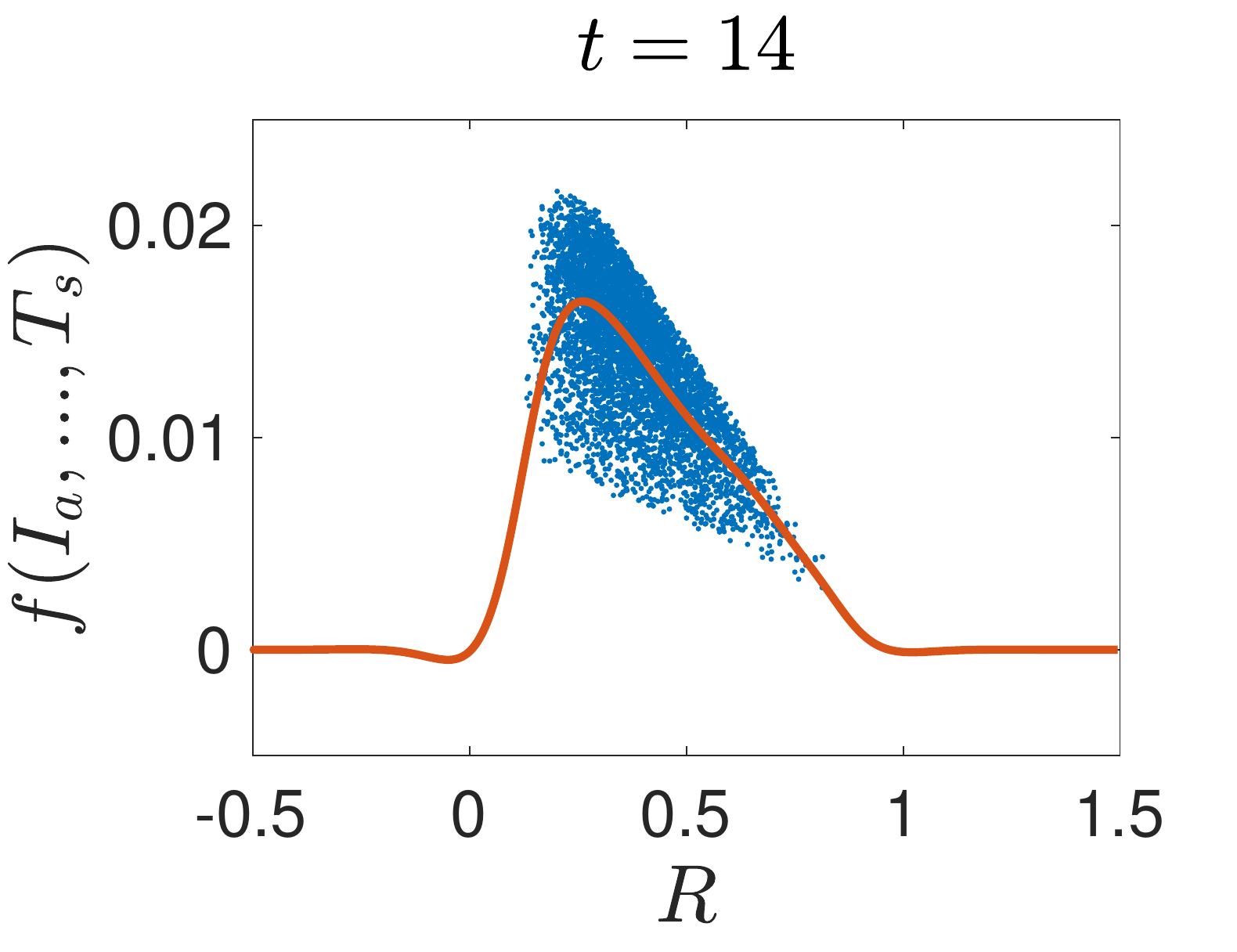}
	\includegraphics[width=0.3\textwidth]{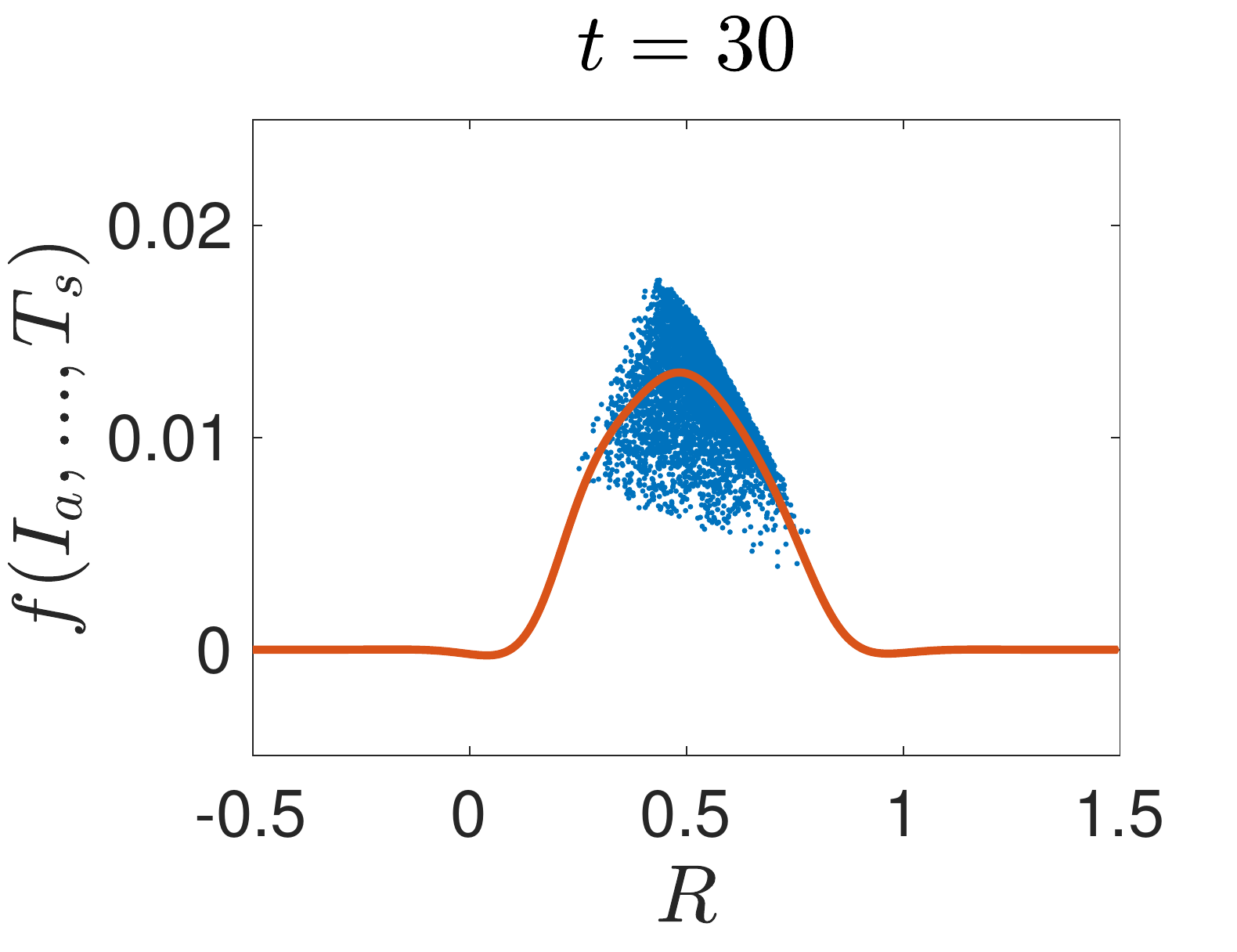}
}

\caption{Data-driven smoothing spline estimation of 
the conditional expectations arising in the study of the 
dynamical systems we study in Section \ref{sec:numerics}. }
	\label{fig:expectations}
\end{figure}

\noindent
In Figure \ref{fig:expectations}, we summarize the results we 
obtain by applying the smoothing spline conditional expectation 
estimator to several dynamical systems studied 
in detail in Section \ref{sec:numerics}. Such systems have 
polynomial-type nonlinearities. The quantity of interest 
in each case is indicated in the $x$-axis of the plots. 
When using this approach, we must be careful to provide 
enough samples for the estimator to adequately capture 
the support of the underlying PDF. If we do not have enough 
samples, the estimator will not be consistent with the true 
conditional expectation. 

\section{Measuring the information content of data with PDEs}
\label{sec:information_content}

In this section, we address the important question of 
whether enough useful data is being injected into the 
reduced-order PDF equation for the purpose of computing an accurate 
numerical solution.
To this end, we develop a new framework based on hyperbolic 
systems that allows us to measure a posteriori the information 
content of data. The key idea relies on developing evolution 
equations for the unclosed terms (conditional expectations) 
appearing in the reduced-order PDF equations. 
As we will see, such equations are generally unclosed, i.e., they 
require data to be solved, but they have the important property 
that their solution can be compared with data. This allows us 
to measure quantitatively whether we have enough 
data to reliably compute the closure approximation.  
In other words, we can measure the information content 
of data by solving a hyperbolic system of PDEs.

\begin{figure}[t]
\centerline{\hspace{0.0cm}
$\mathbb{E}\left[x^3_2(t)|x_1(t)\right]$ \hspace{3.1cm} 
$\mathbb{E}\left[x^2_2(t)|x_1(t)\right]$ \hspace{3.3cm}
$\mathbb{E}[x_3(t)|x_1(t)]$
}
\centerline{
\includegraphics[width=0.33\textwidth]{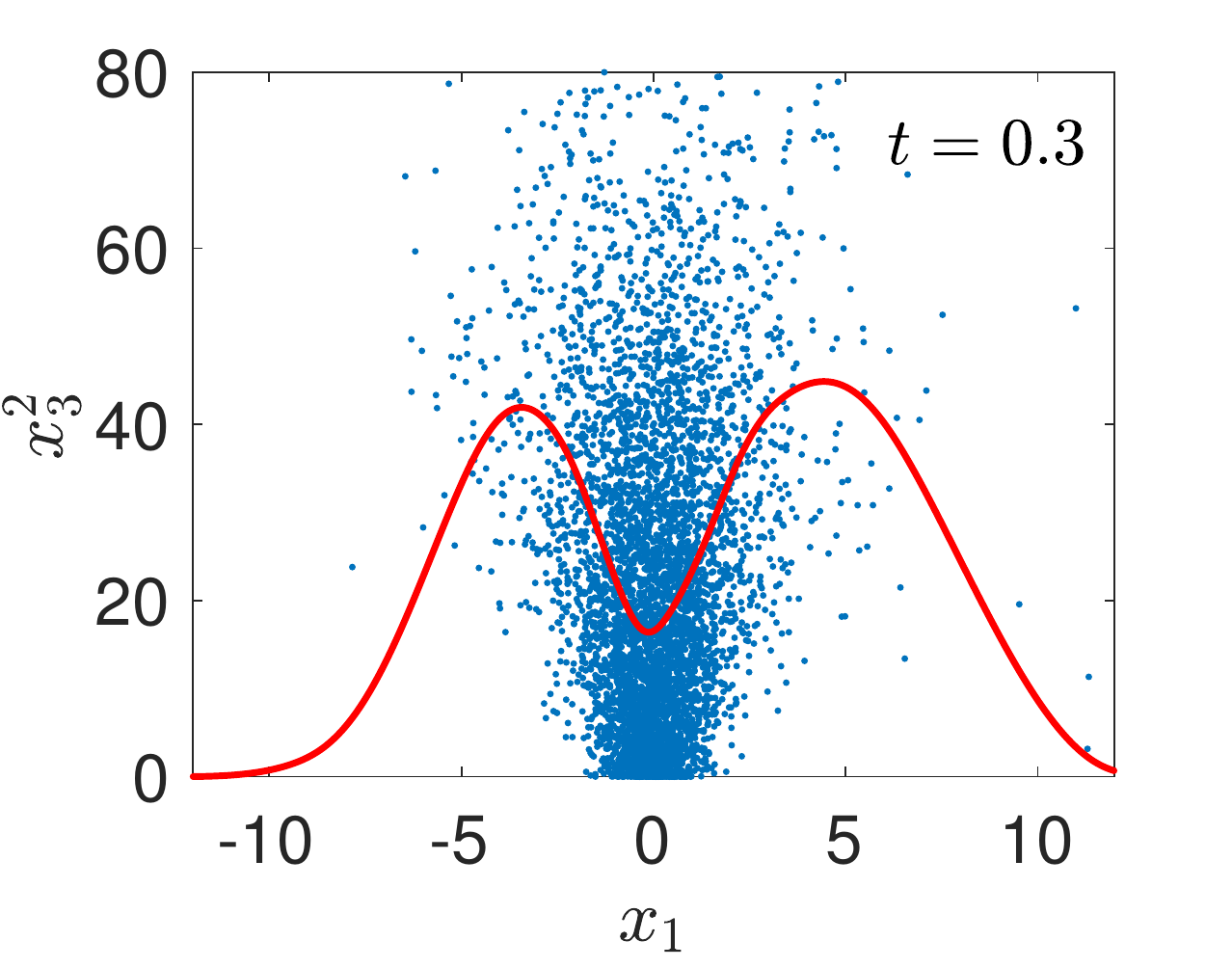}
\includegraphics[width=0.33\textwidth]{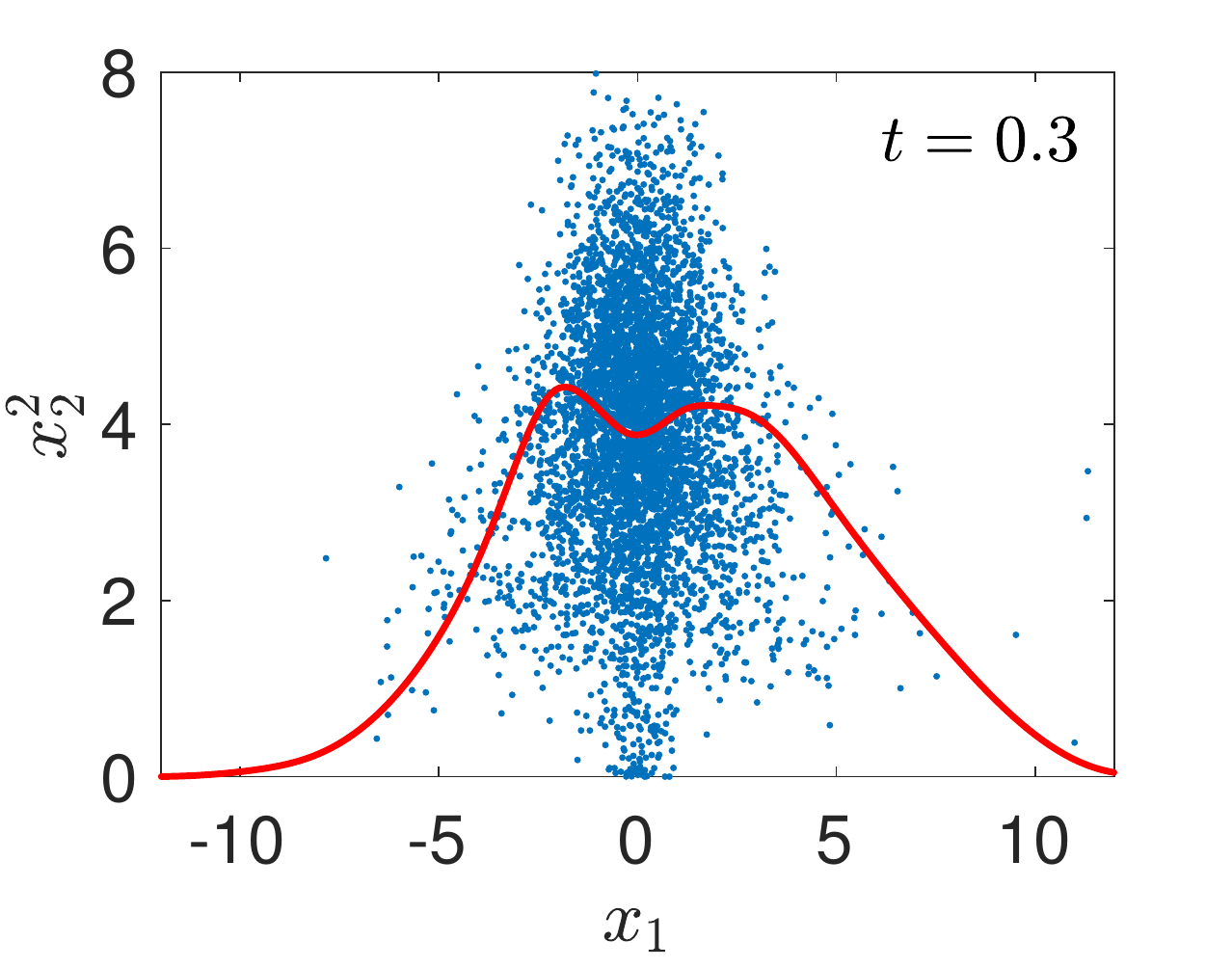}
\includegraphics[width=0.33\textwidth]{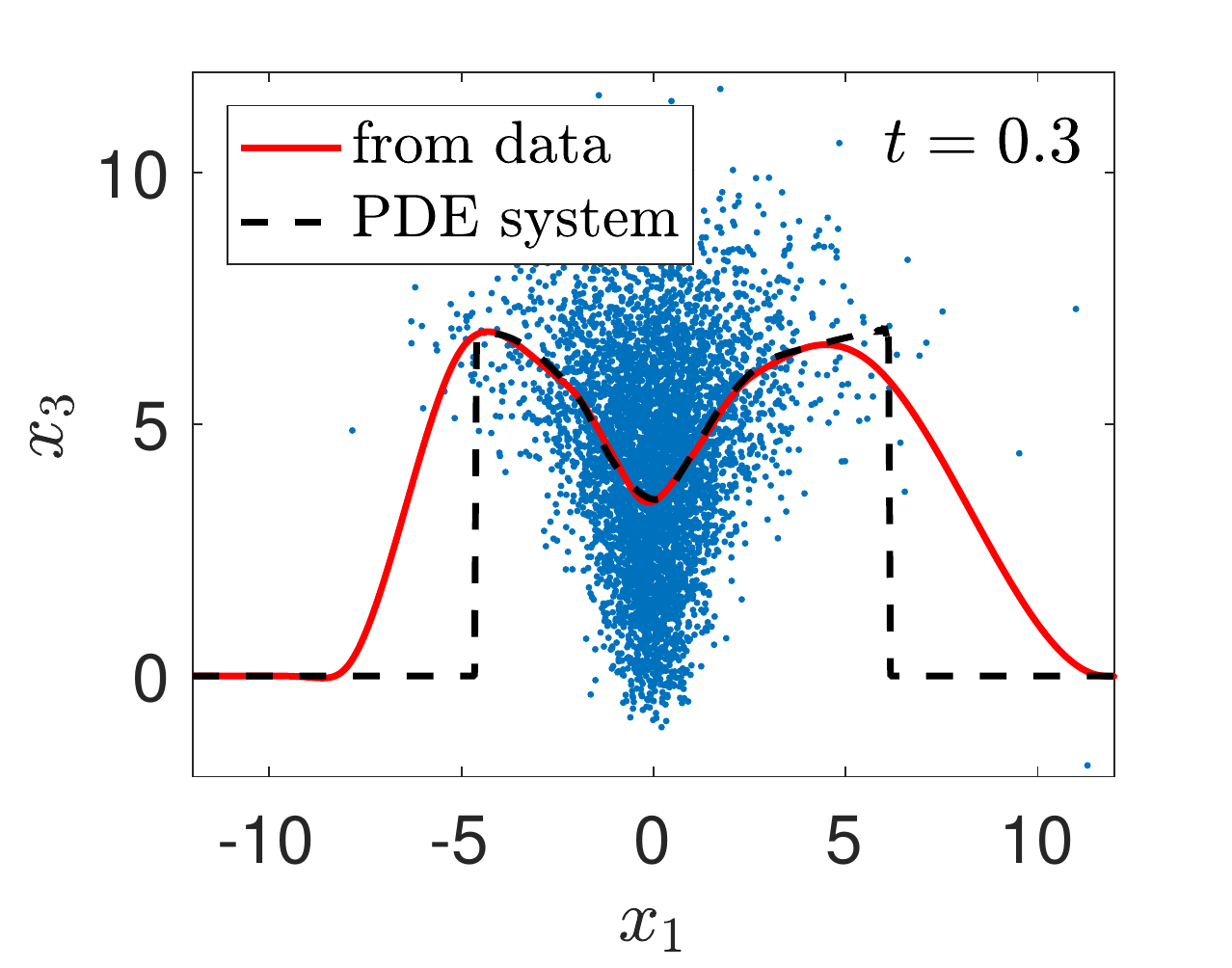}
}
\centerline{
\includegraphics[width=0.33\textwidth]{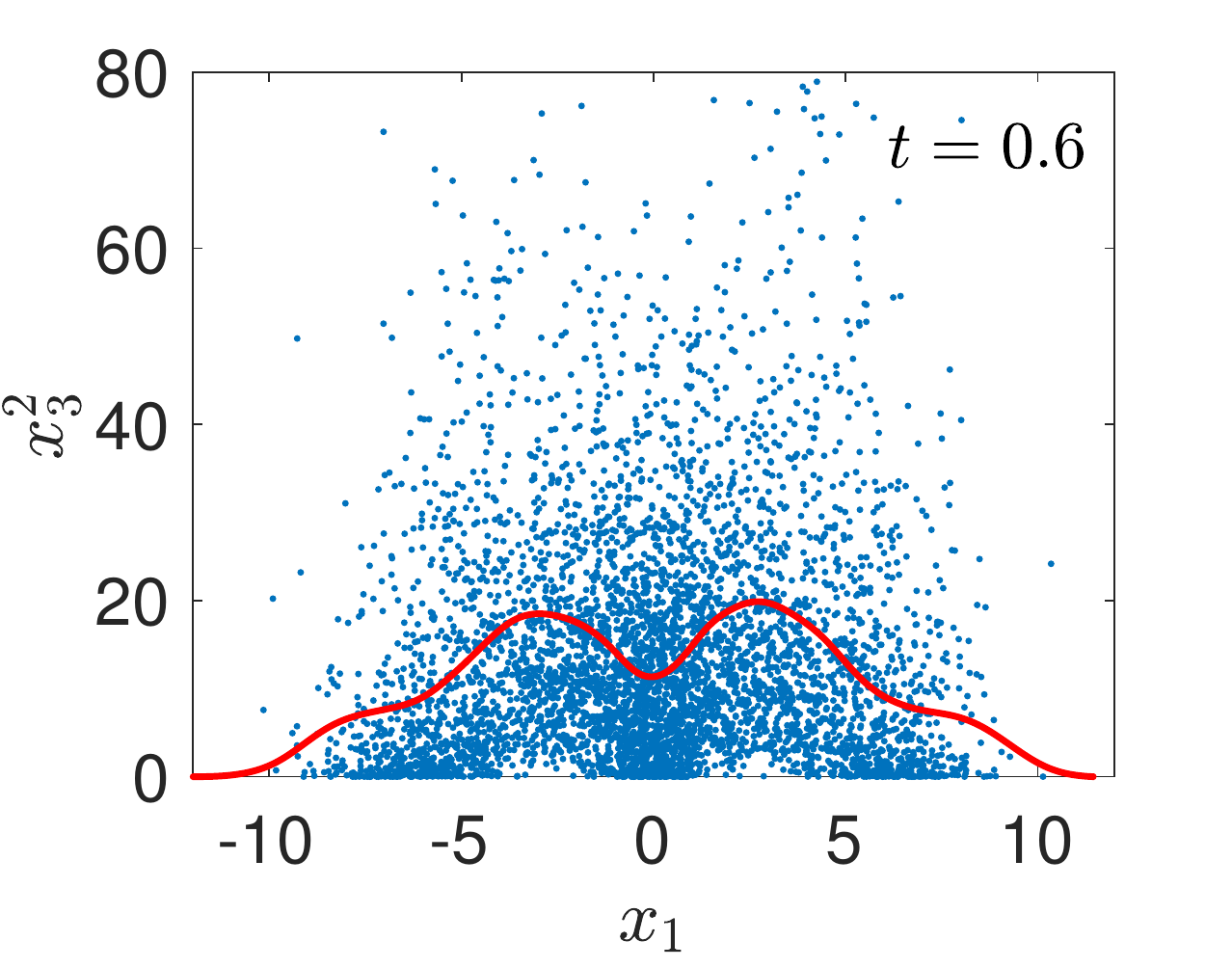}
\includegraphics[width=0.33\textwidth]{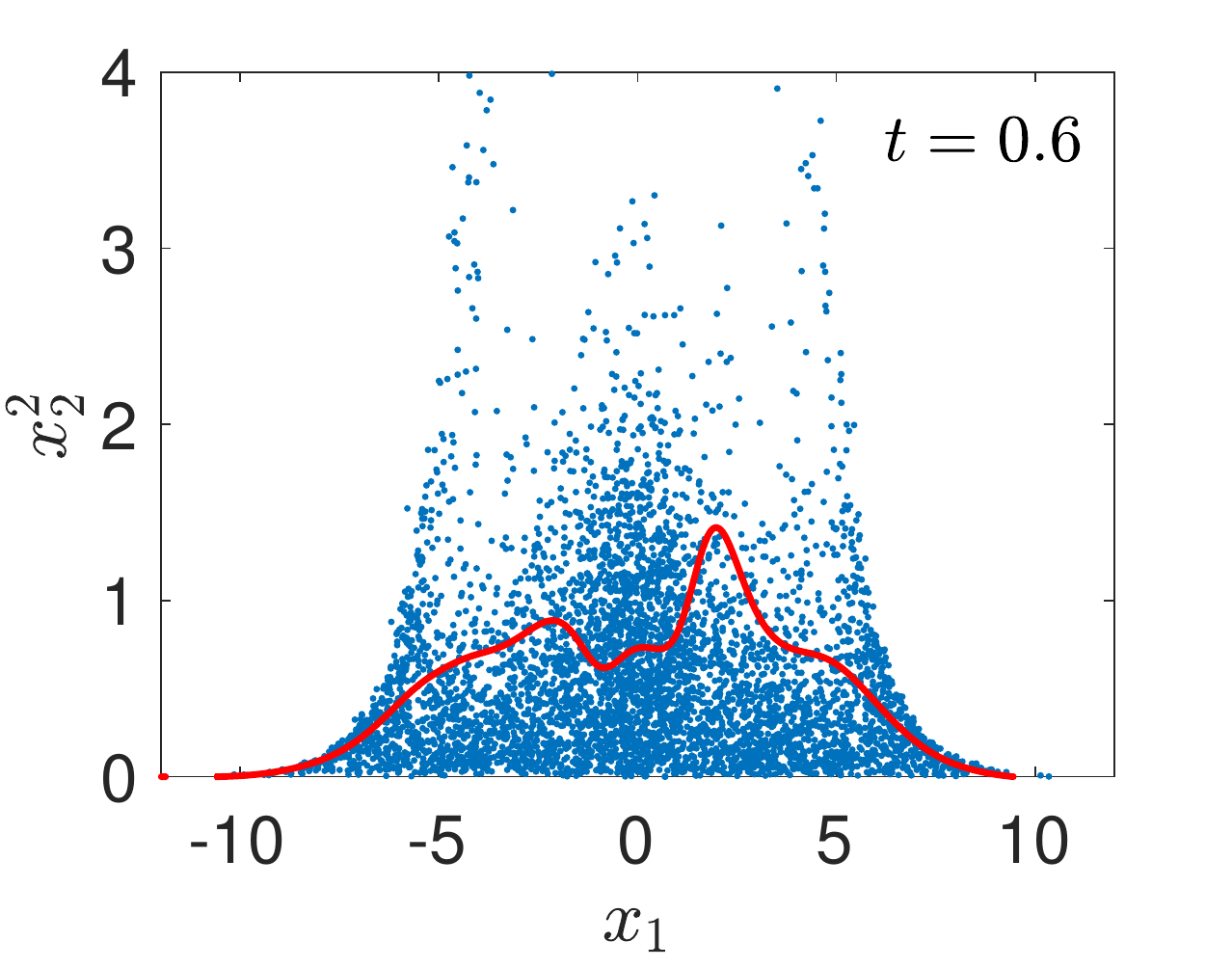}
\includegraphics[width=0.33\textwidth]{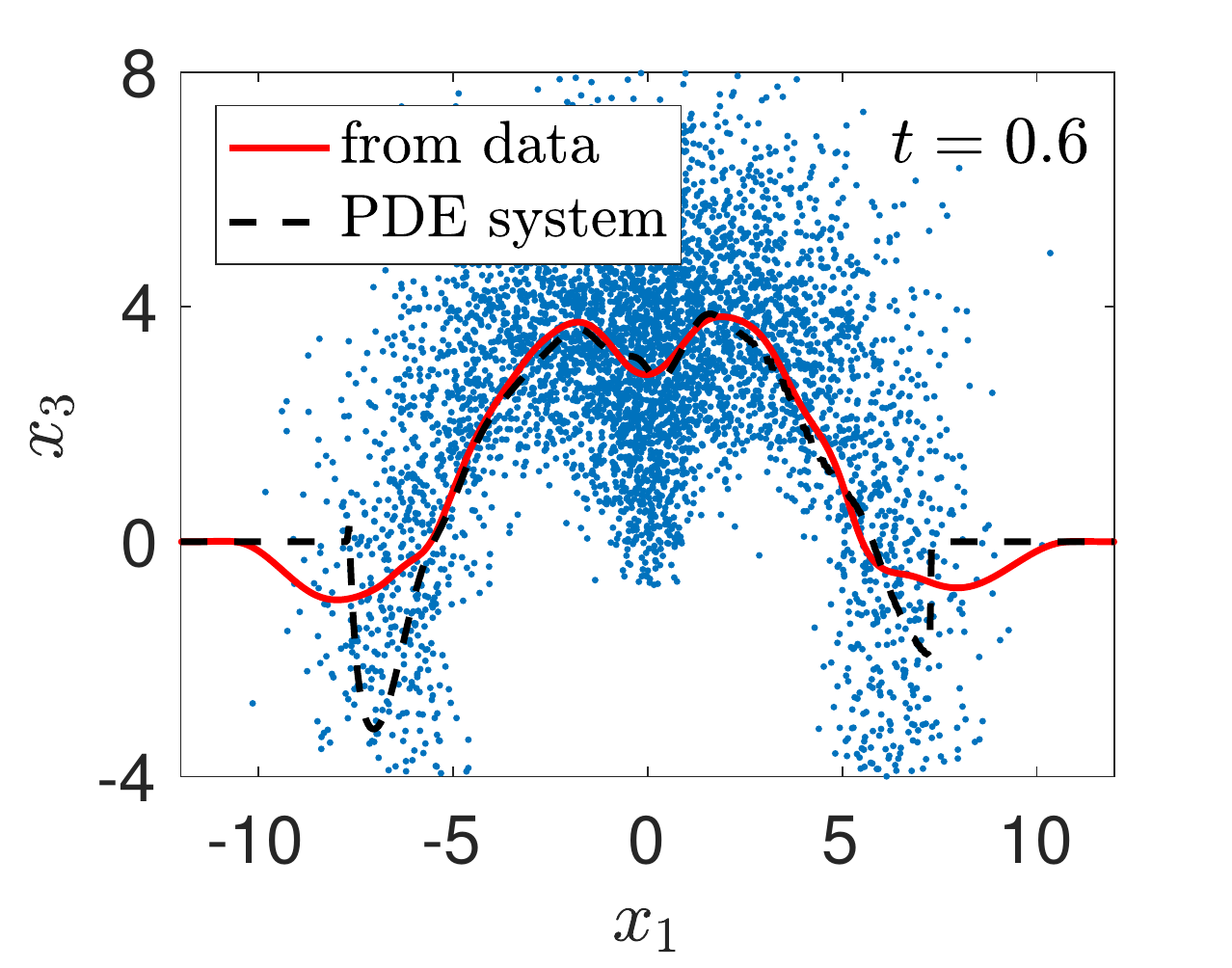}
}
\caption{Kraichnan-Orszag system. 
Data-driven estimates of the conditional expectations $\mathbb{E}[x^2_3(t)|x_1(t)]$ (first column) and $\mathbb{E}[x^2_2(t)|x_1(t)]$ (second column) at different times based on 5000 sample trajectories. 
In the third column we compare the conditional expectation  
$\mathbb{E}[x_3(t)|x_1(t)]$ we obtain from data with the one computed 
by solving the nonlinear PDE system \eqref{BBGKY_marginal1}-\eqref{condexpKOequation} with $\mathbb{E}[x^2_3(t)|x_1(t)]$ and $\mathbb{E}[x^2_2(t)|x_1(t)]$ estimated from data. It is seen that the two 
conditional expectations coincide within the active region of 
the reduced-order PDF $p(x_1,t)$. 
The error in $\mathbb{E}[x_3(t)|x_1(t)]$ provides 
a measure of the information content of the sample set obtained from \eqref{KO-dyn}.}
\label{fig:information_content}
\end{figure}

To describe the method, let us consider again the 
Kraichnan-Orszag system \eqref{KO-dyn}. As before, 
suppose we are interested in the phase space 
function $u(\bm x)=x_1$ (first-component of the system). 
We have seen that the dynamics of the PDF of $x_1$ is 
governed by the unclosed transport equation 
\eqref{BBGKY_marginal1}, where the conditional expectation 
$\mathbb{E}[x_3(t)|x_1(t)]$ can be estimated directly 
from sample trajectories. Recall that such conditional 
expectation is defined in \eqref{condexpKO}. 
By differentiating such expression with respect to 
time we obtain
\begin{equation}
\frac{\partial \mathbb{E}[x_3(t)|x_1(t)]}{\partial t}= 
-\frac{1}{p(x_1,t)}\frac{\partial p(x_1,t)}{\partial t}
\mathbb{E}[x_3(t)|x_1(t)] + \frac{1}{p(x_1,t)}
\int_{-\infty}^{\infty} x_3\frac{\partial p(x_1,x_3,t)}{\partial t}dx_3.
\label{co1}
\end{equation}
By substituting \eqref{BBGKY_marginal1} and the 
reduced-order PDF equation for $p(x_1,x_3,t)$ obtained by integrating 
the Liouville equation \eqref{KOLiouville} with respect to $x_2$ into 
\eqref{co1}, we find 
\begin{align}
\frac{\partial\mathbb{E}[x_3(t)|x_1(t)] }{\partial t}= &
\frac{\mathbb{E}[x_3(t)|x_1(t)]}{p(x_1,t)}
\frac{\partial }{\partial x_1}\left(x_1p(x_1,t)\mathbb{E}[x_3(t)|x_1(t)]\right) -\cdots \nonumber \\
& \frac{1}{p(x_1,t)}
\frac{\partial }{\partial x_1}\left(x_1p(x_1,t) \mathbb{E}[x^2_3(t)|x_1(t)]\right) - x_1^2 + \mathbb{E}[x^2_2(t)|x_1(t)].
\label{condexpKOequation}
\end{align}
This is the formally exact evolution equation for the conditional
expectation of $x_3(t)$ given $x_1(t)$ in the Kraichnan-Orszag 
system.  The solution to the nonlinear PDE system
\eqref{BBGKY_marginal1}-\eqref{condexpKOequation} can be 
computed in a data-driven setting by estimating  
$\mathbb{E}[x^2_3(t)|x_1(t)]$ and $\mathbb{E}[x^2_2(t)|x_1(t)]$ 
from sample trajectories of \eqref{KO-dyn} as we discussed in Section 
\ref{sec:estimating_conditional_expectations}. 
The conditional expectation $\mathbb{E}[x_3(t)|x_1(t)]$ 
we obtain by solving the system 
\eqref{BBGKY_marginal1}-\eqref{condexpKOequation} 
can be then compared with its data-driven estimate. 
This provides an indication of whether 
we have sufficient data to compute an accurate closure 
of \eqref{BBGKY_marginal1}. This procedure is illustrated 
in Figure \ref{fig:information_content}.
To avoid numerical issues, we solve 
equation \eqref{condexpKOequation} only in the regions where 
the PDF $p(x_1,t)$ is larger than a threshold, hence the jump 
in the black dashed line observed in the rightmost plots
of Figure \ref{fig:information_content}. As we will see in the 
next Section, this issue can be avoided if we change the 
coordinates appropriately.

\vs
\noindent  
The methodology we just described for the Kraichnan-Orszag 
system can be easily generalized and applied to arbitrary nonlinear 
systems in the form \eqref{dyn}, and any observable \eqref{observable}. 
In Section \ref{sec:numerics}, we study two examples.

\subsection{Computational aspects}
\begin{figure}
\centerline{\hspace{0.3cm}(a)\hspace{8cm}(b)}
\centerline{
\includegraphics[width=0.5\textwidth]{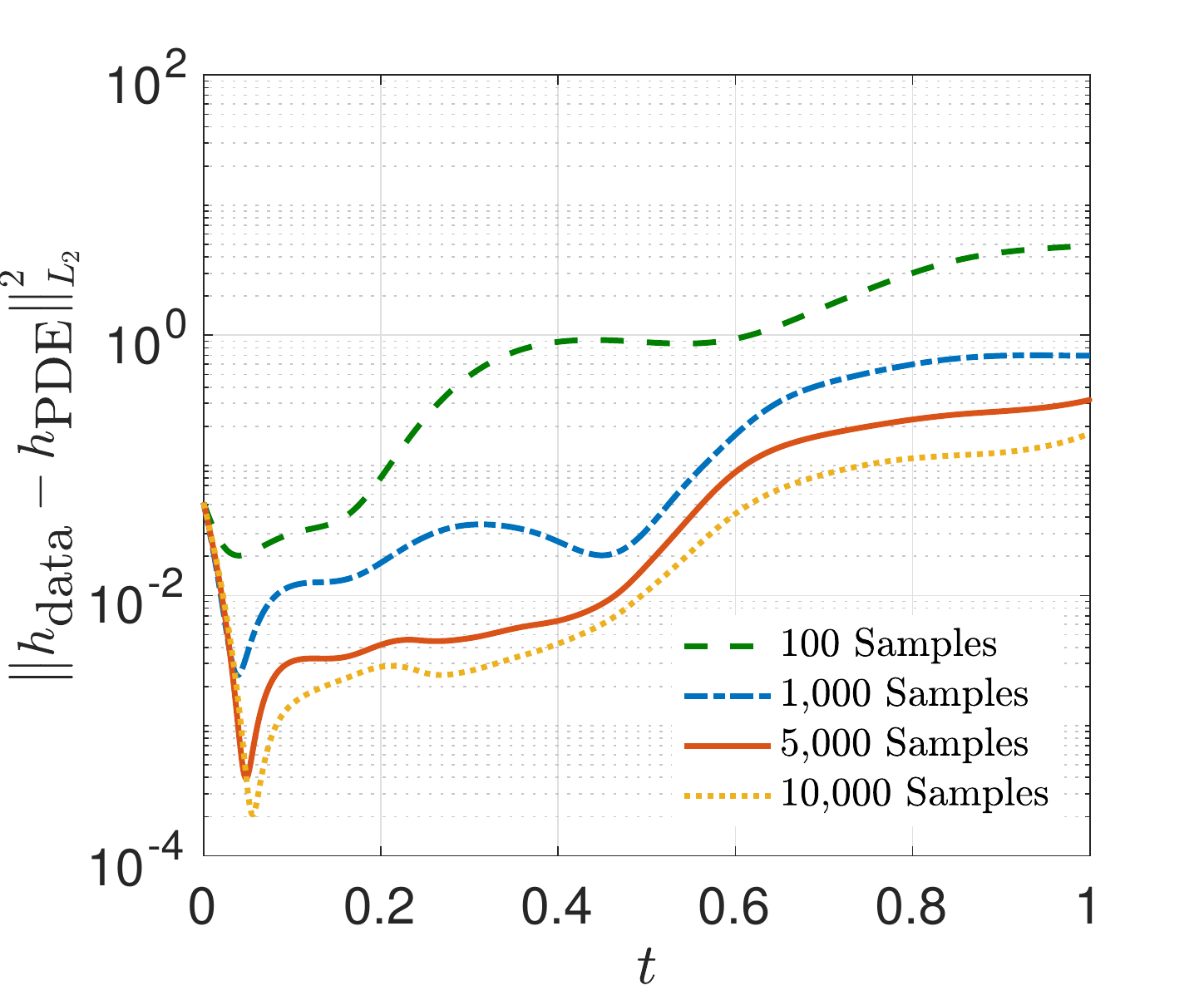}
\includegraphics[width=0.5\textwidth]{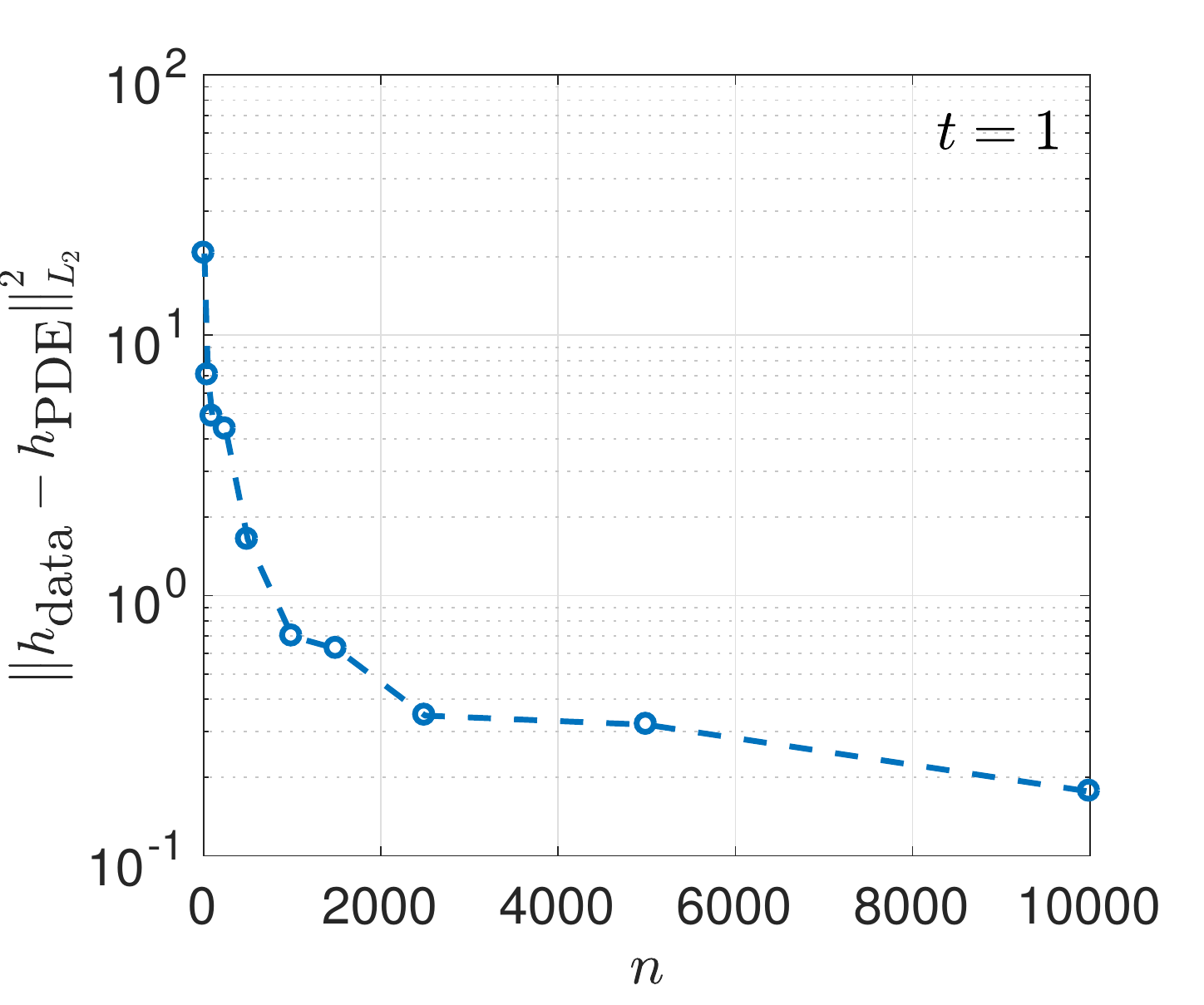}
}

\caption{Kraichnan-Orszag three-mode problem. (a) Time-dependent 
errors in the function \eqref{hdef} plotted for a variety of sample 
sizes. (b) Decay of the error at $t=1$ as a function of 
the number of sample trajectories.}
\label{fig:information_content_h}
\end{figure}

Solving the nonlinear PDE system \eqref{BBGKY_marginal1}-\eqref{condexpKOequation} numerically is not straightforward. 
There are indeed several subtleties and difficulties that are 
hard to overcome. Perhaps, the most relevant 
 is related to the presence of terms at the right hand side of 
\eqref{condexpKOequation} multiplied by $1/p(x_1,t)$. After 
simplification, these terms can be written as 
 $\partial \log(p(x_1,t))\partial x_1$. For instance,
we have  
\begin{align}
\frac{\mathbb{E}[x_3(t)|x_1(t)]}{p(x_1,t)}
\frac{\partial }{\partial x_1}\left(x_1p(x_1,t)\mathbb{E}[x_3(t)|x_1(t)]\right)=&\mathbb{E}[x_3(t)|x_1(t)]
\frac{\partial }{\partial x_1}\left(x_1\mathbb{E}[x_3(t)|x_1(t)]\right)+\cdots \nonumber\\
& x_1\mathbb{E}[x_3(t)|x_1(t)]^2 \frac{\partial \log(p(x_1,t))}{\partial x_1}.
\end{align}
Clearly, such terms can easily yield numerical overflow in regions 
where $p(x_1,t)$ is very small. This problem can be mitigated by 
using adaptive algorithms that can track the support of the PDF 
$p(x_1,t)$ (see, e.g., \cite{Heyrim}). An alternative approach
relies on coordinate transformation. In particular, rather than 
solving equation \eqref{condexpKOequation} for 
$\mathbb{E}[x_3(t),x_1(t)]$, we can solve it for the product 
of $\mathbb{E}[x_3(t),x_1(t)]$ and $p(x_1,t)$. This product 
represents the integral of $x_3p(x_1,x_3,t)$ with respect 
to $x_3$.
Let us define\footnote{Recall that 
\begin{equation}
\mathbb{E}[x_3(t)|x_1(t)] =\int_{-\infty}^\infty x_3 p(x_3|x_1,t)dx_3=
\frac{1}{p(x_1,t)}\int_{-\infty}^\infty x_3 p(x_1,x_3,t) dx_3.
\end{equation}
 Therefore, 
\begin{equation}
\int_{-\infty}^\infty x_3 p(x_1,x_3,t) dx_3 = p(x_1,t) 
\mathbb{E}[x_3(t)|x_1(t)]. 
\end{equation}}
\begin{align}
h(x_1,t) = &\int_{-\infty}^\infty x_3 p(x_3,x_1,t)dx_3
= p(x_1,t) \mathbb{E}[x_3(t)|x_1(t)].
\label{hdef}
\end{align}
The evolution equation for $h(x_1,t)$ can 
be obtained by following the same steps we followed to 
derive equation \eqref{condexpKOequation}. This yields,  
\begin{equation}
\frac{\partial h(x_1,t)}{\partial t} =-\frac{\partial }{\partial x_1}\left(x_1p(x_1,t)\mathbb{E}\left[\left.x_3(t)^2\right|x_1(t)\right]\right)-x_1^2p(x_1,t) +
\mathbb{E}\left[\left.x_2(t)^2\right|x_1(t)\right]p(x_1,t).
\label{KOu}
\end{equation} 
On the other hand, equation \eqref{BBGKY_marginal1} can be written 
in terms of $h(x_1,t)$ as 
\begin{equation}
\frac{\partial p(x_1,t)}{\partial t} = -\frac{\partial }{\partial x_1}
\left(x_1 h(x_1,t)\right). 
\label{KOPDF_Eq1}
\end{equation}
The hyperbolic system \eqref{KOu}-\eqref{KOPDF_Eq1} is linear, 
with two unclosed terms represented by the conditional expectations 
$\mathbb{E}[x_3(t)^2|x_1(t)]$ and $\mathbb{E}[x_2(t)^2|x_1(t)]$.
The solution can be computed in a data-driven setting by estimating  
these conditional expectations from sample trajectories by 
using the methods of Section \ref{sec:estimating_conditional_expectations}. 
The advantages of solving  \eqref{KOu}-\eqref{KOPDF_Eq1} over solving 
directly \eqref{BBGKY_marginal1} rely on the fact that once $h(x_1,t)$ 
and $p(x_1,t)$ are available, then we can immediately compute 
the conditional expectation $\mathbb{E}\left[x_3(t)|x_1(t)\right])$
using \eqref{hdef}, and compare it with the data-driven 
estimate. This provides a measure of the information content of 
sample trajectories for the particular phase space function 
we are interested in, i.e., $u(\bm x(t))=x_1(t)$ in this case. 
In Figure \ref{fig:information_content_h}, we 
plot the time-dependent $L_2$ error between 
the function \eqref{hdef} we obtain from data and the one we obtain 
by solving the PDE system  \eqref{KOu}-\eqref{KOPDF_Eq1}. The 
 benchmark solution of $h(x_1,t)$ is obtained by estimating 
$\mathbb{E}\left[x_3(t)|x_1(t)\right])$ from data and then 
multiplying it by an accurate kernel density estimate of $p(x_1,t)$. 

{\remark $\,$} Equation \eqref{KOu} can be immediately 
integrated in time to obtain
\begin{align}
h(x_1,t) =&  h(x_1,0) - \frac{\partial }{\partial x_1}
\left(x_1\int_{0}^tp(x_1,\tau)\mathbb{E}\left[\left.x_3(\tau)^2\right|x_1(\tau)\right]d\tau \right)- \cdots\nonumber\\
& x_1^2\int_{0}^t p(x_1,\tau)d\tau + 
\int_0^t\mathbb{E}\left[\left.x_2(\tau)^2\right| x_1(\tau)\right]p(x_1,\tau)d\tau.
\label{KOu1}
\end{align}
A substitution of \eqref{KOu1} into \eqref{KOPDF_Eq1} yields 
a rather complicated integro-differential PDF equation for $p(x_1,t)$. 
Such equation has exactly the same solution as  
equation \eqref{BBGKY_marginal1}, although the unclosed 
terms (conditional expectations) that need to be estimated 
from sample trajectories are different.
\vs
 
\begin{figure}[t]
\centerline{\hspace{-0.5cm}(a)\hspace{6.4cm} (b)}
\centerline{
		\includegraphics[width=0.41\textwidth]{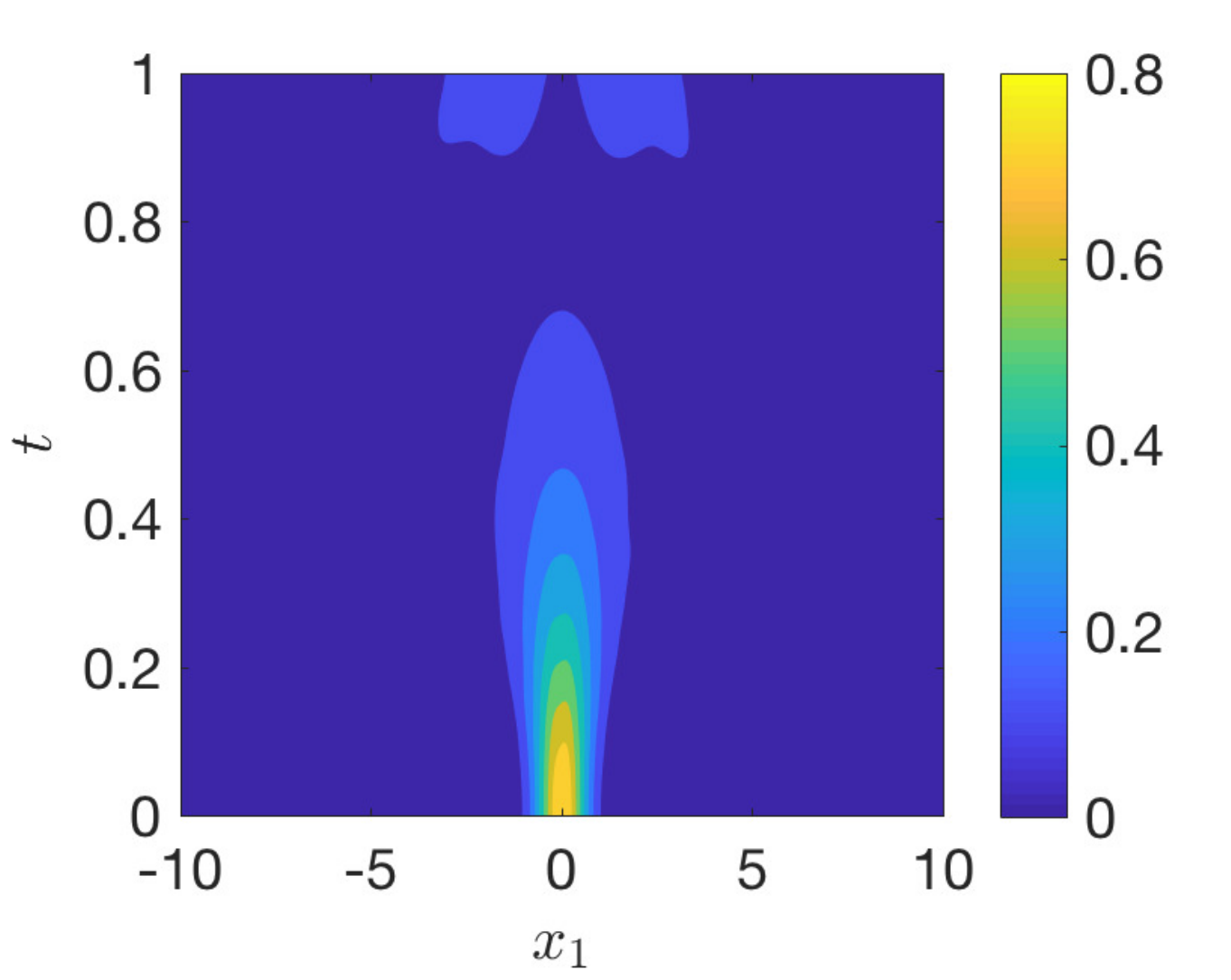}
		\includegraphics[width=0.41\textwidth]{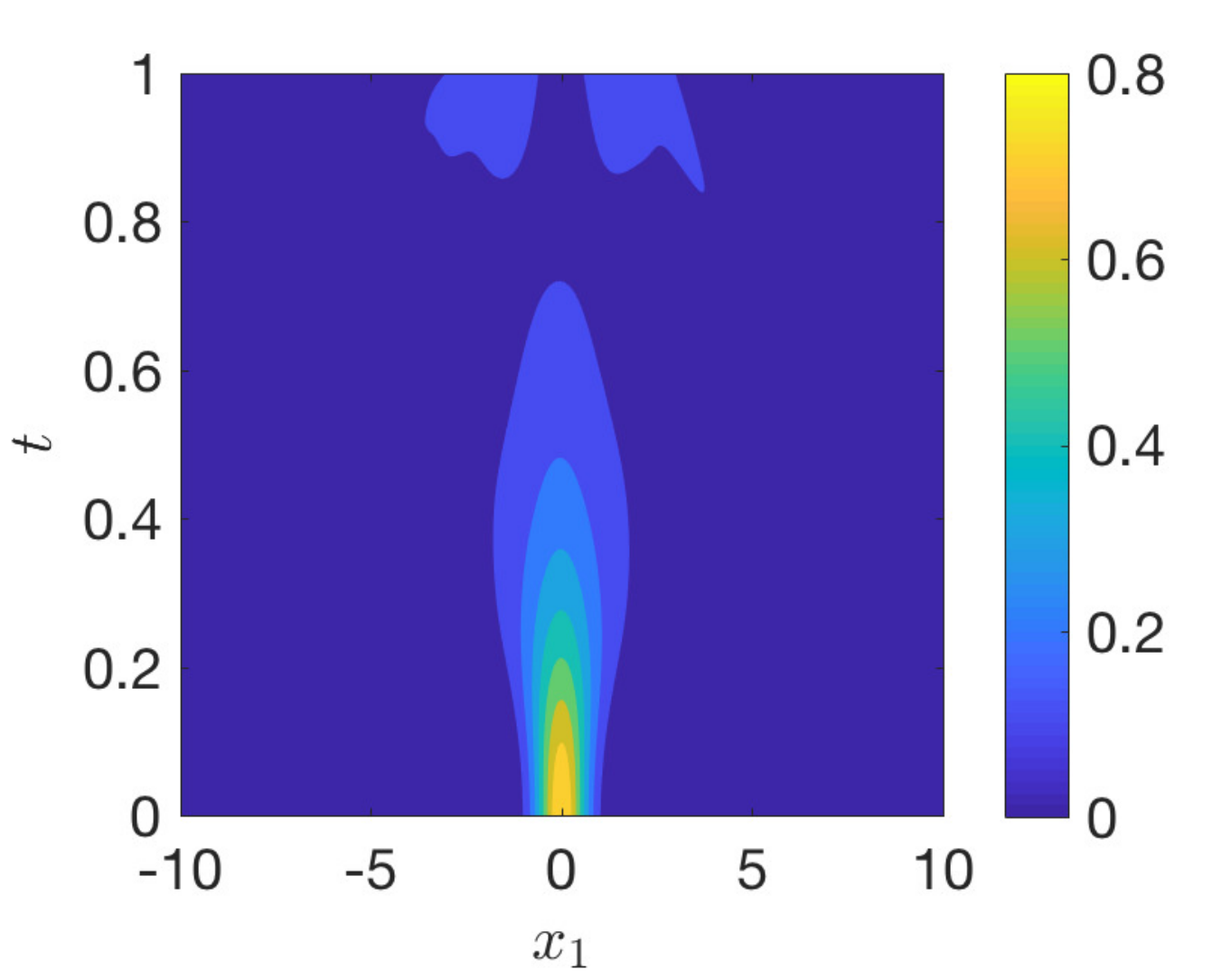}}
\centerline{
\includegraphics[width=0.28\textwidth]{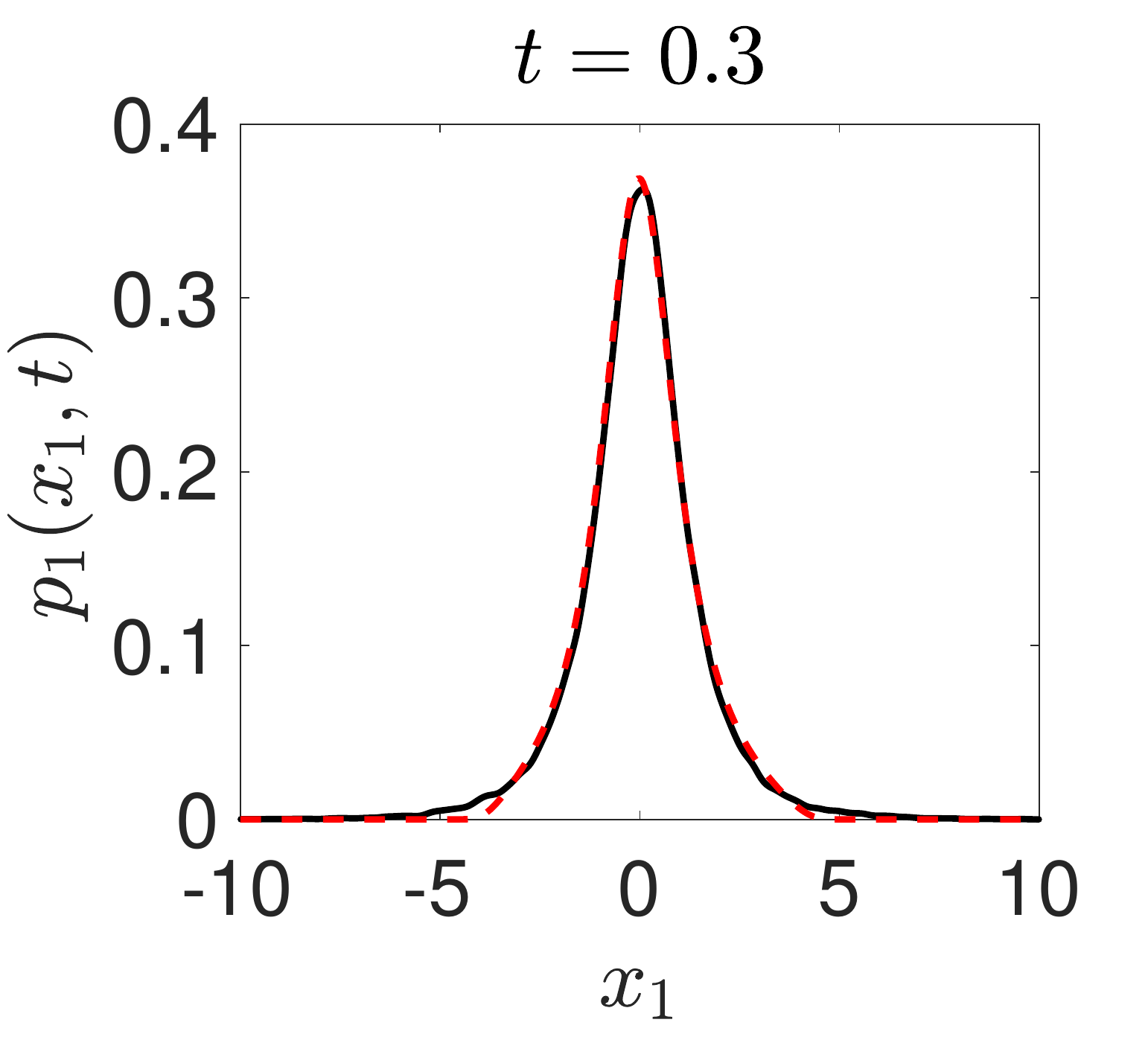}
\includegraphics[width=0.28\textwidth]{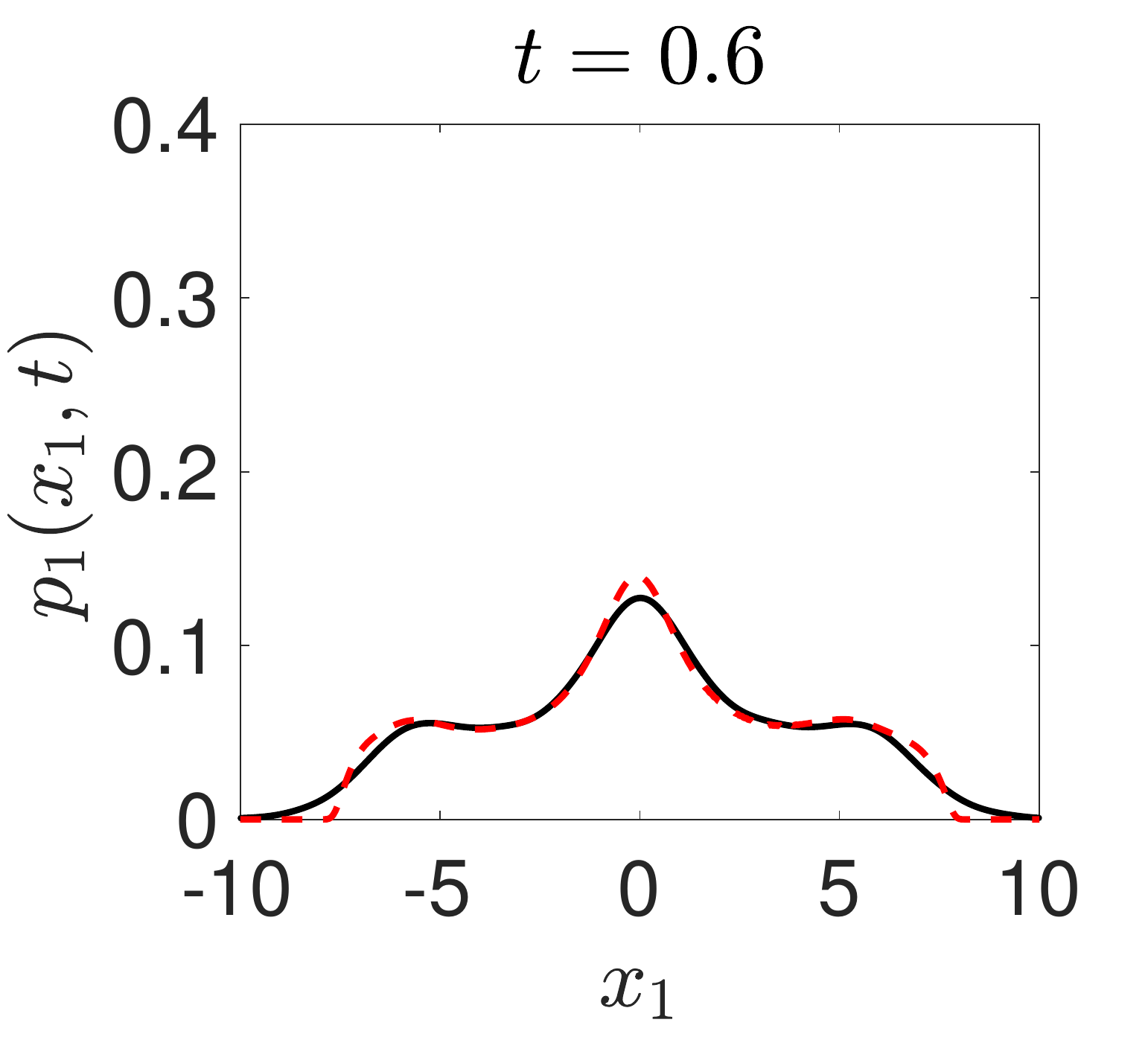}
\includegraphics[width=0.28\textwidth]{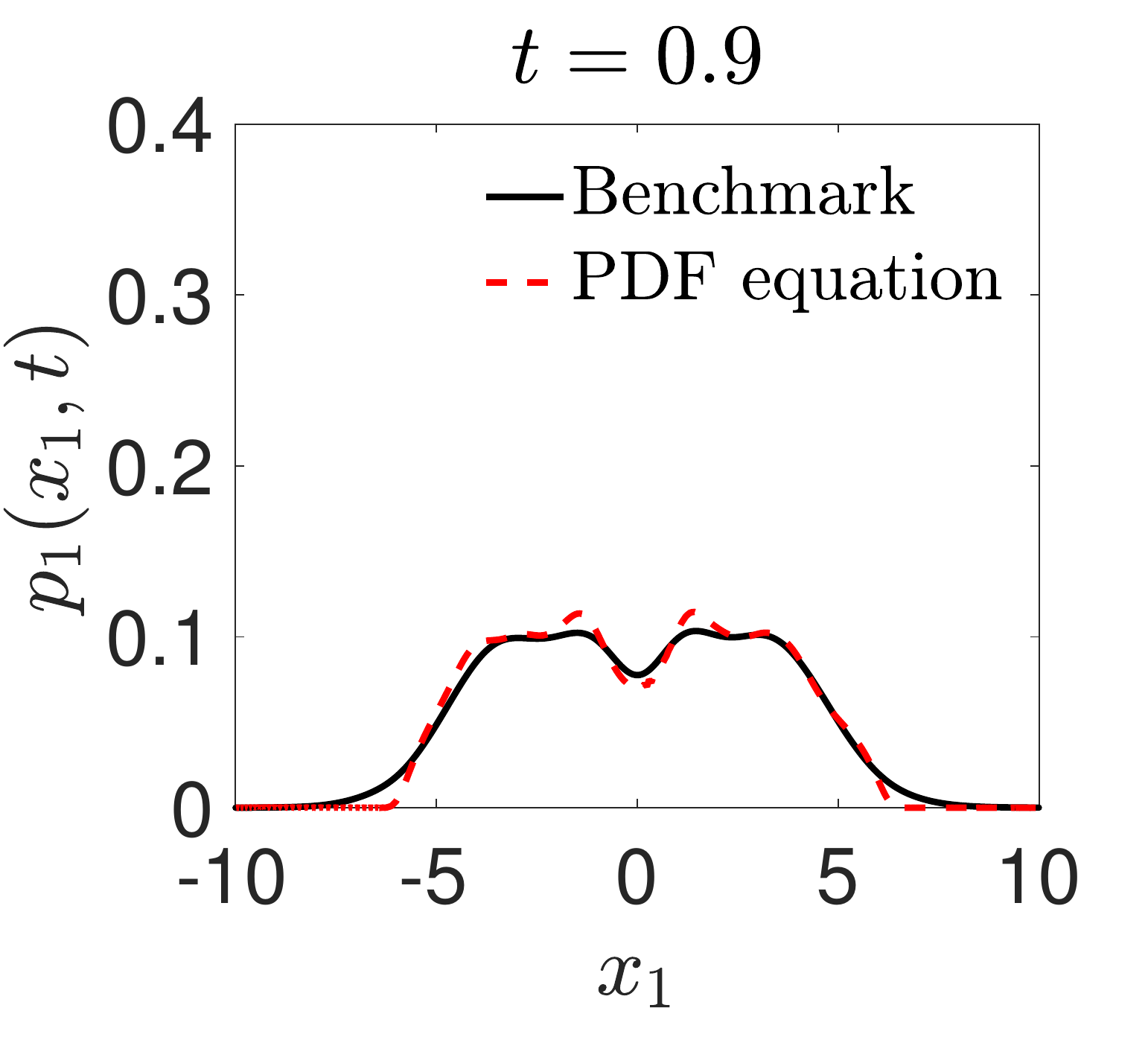}
}
\caption{Kraichnan-Orszag three-mode problem. (a)
Accurate kernel density estimate of $p_1(x_1, t)$ 
based on $30000$ sample trajectories. (b)  Numerical solution of \eqref{BBGKY_marginal1} obtained by estimating $\mathbb{E}[x_3(t)|x_1(t)]$
with $5000$ sample trajectories.}
\label{fig:KO_PDF}
\end{figure}

\begin{figure}[t]
\centerline{\hspace{-0.5cm}(a)\hspace{7cm} (b)}
\centerline{
		\includegraphics[width=0.45\textwidth]{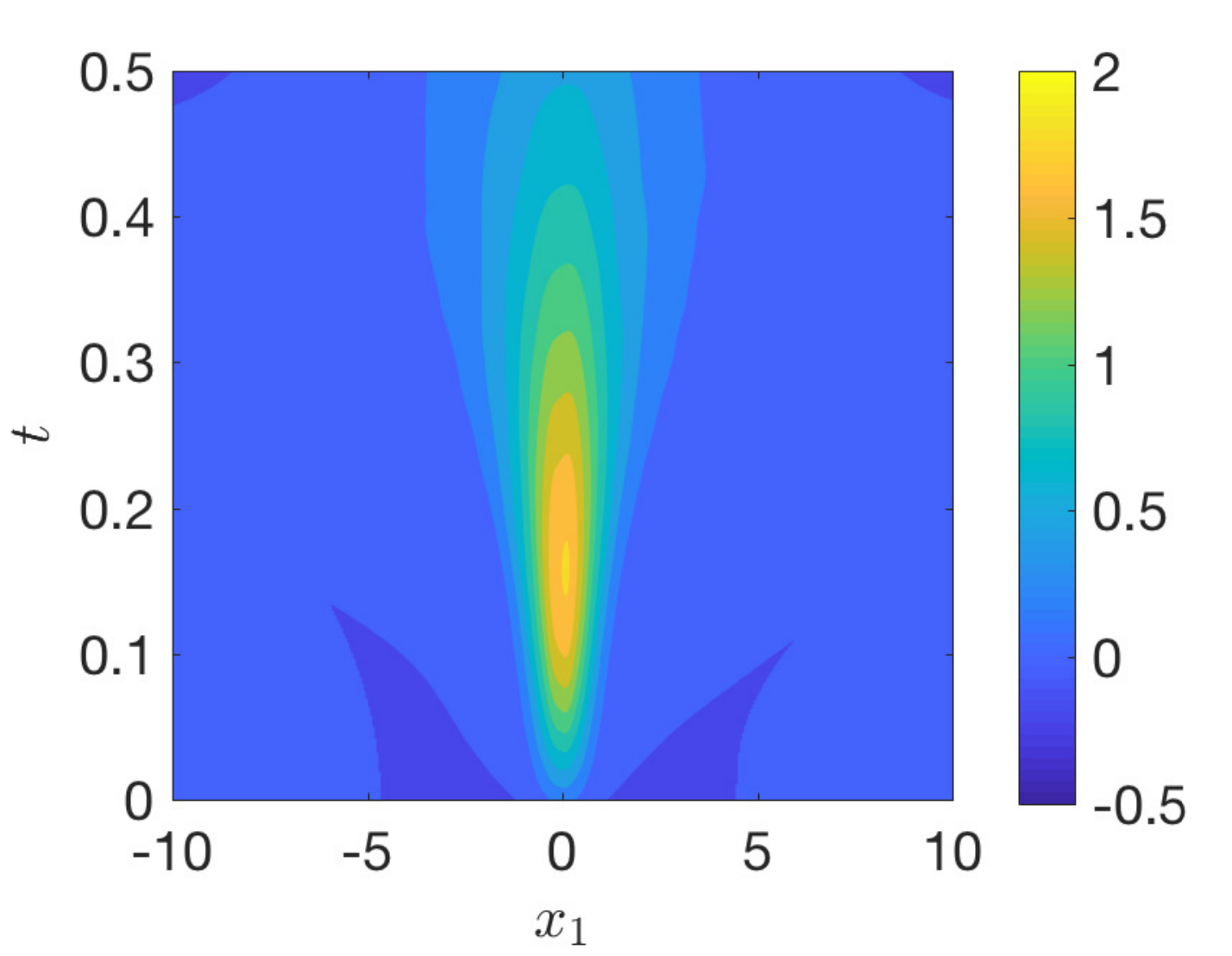}
		\includegraphics[width=0.45\textwidth]{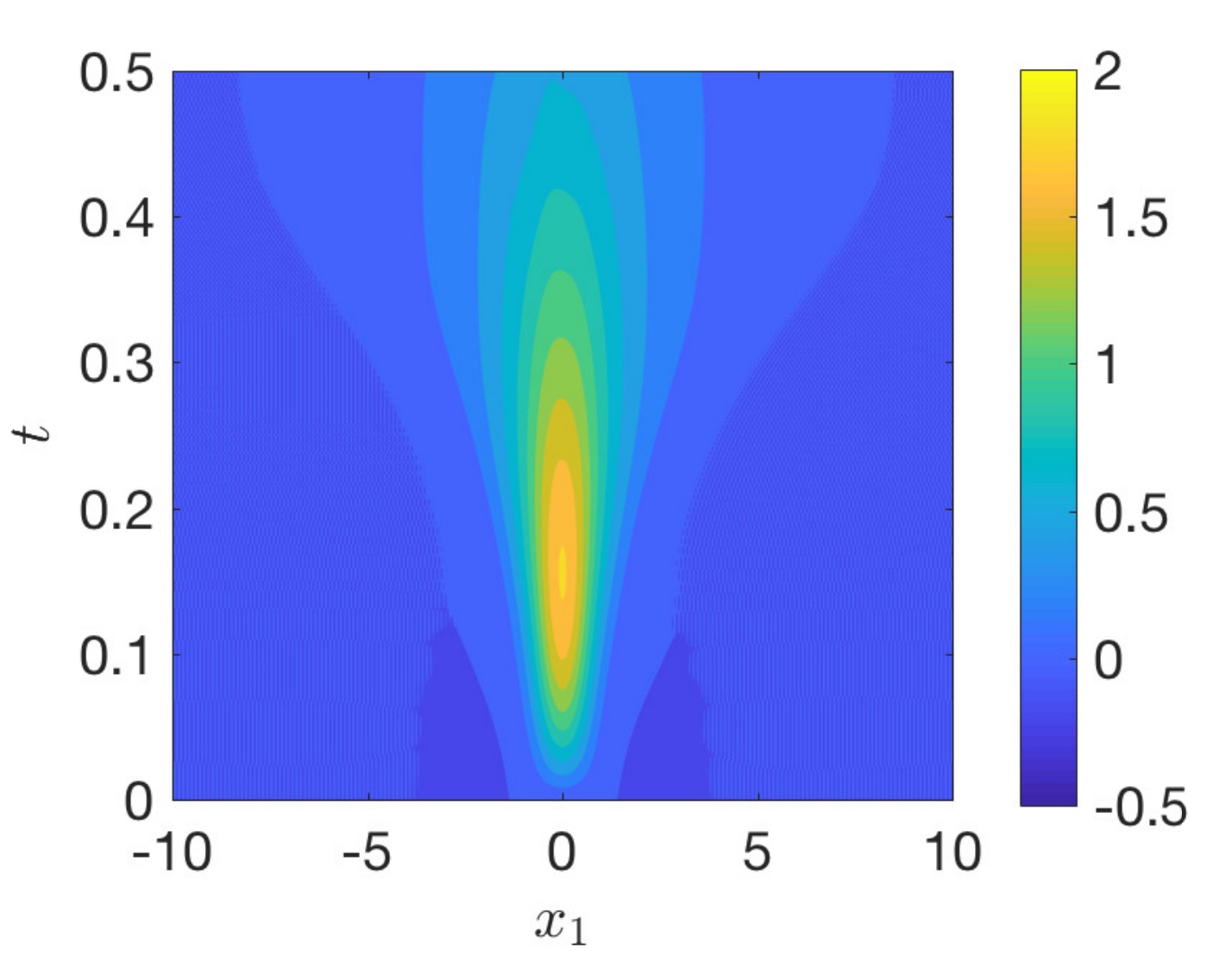}}
\caption{Kraichnan-Orszag three-mode problem.  Comparsion between 
the time evolution of the function \eqref{hdef} we obtain from data 
($30000$ sample paths) (a),  
and from the numerical solution of the system of equations \eqref{KOu}-\eqref{KOPDF_Eq1} (b). In the latter case, the unclosed conditional 
expectations in \eqref{KOu} were estimated with smoothing splines 
and $5000$ samples.}
\label{fig:h}
\end{figure}

\noindent
In Figure \ref{fig:KO_PDF} we plot the PDF dynamics 
we obtain by solving \eqref{BBGKY_marginal1} with an 
accurate Fourier spectral method. The conditional expectation is 
estimated based on $5000$ sample trajectories. 
In Figure \eqref{fig:h} we compare the function 
\eqref{hdef} we obtain from data ($30000$ sample trajectories) 
to the numerical solution of the system of equations \eqref{KOu}-\eqref{KOPDF_Eq1}. The unclosed terms were estimated 
with $5000$ sample paths. 

\section{Data-driven estimation of the Mori-Zwanzig memory integral} 
\label{sec:data-drivenMZ}
The Mori-Zwanzig (MZ) formulation is a technique of irreversible 
statistical mechanics that allows us to formally integrate 
out an arbitrary number of phase variables in nonlinear 
dynamical systems. In doing so, we obtain exact evolution equations 
for quantities of interest such as macroscopic observables 
in high-dimensional phase spaces \cite{VenturiBook,Venturi_MZ,Chorin,Dominy2017}. 
To describe the method, consider the dynamical 
system \eqref{dyn} and assume that the components of 
the random initial state $\bm x_0(\omega)$ are statistically 
independent\footnote{The general case where the 
components of $\bm x_0(\omega)$ are not statistically 
independent can be treated similarly.}. 
Furthermore, suppose we are interested in the PDF of the 
first component of the system, i.e., $p(x_1,t)$. 
To study the dynamics of such PDF, we define the 
following projection operator  
\begin{equation}
P f(\bm x) = \left<f(\bm x)\right>\prod_{j=2}^n p(x_j,0) \qquad 
\left<f\right>=\int_{-\infty}^\infty\cdots\int_{-\infty}^\infty f(\bm x) 
\prod_{j=2}^n dx_j 
\end{equation}
as well as the complementary projection $Q=I-P$. Note that $P$ sends 
the joint PDF $p(\bm x,t)$ into the separated state 
\begin{equation}
P p(\bm x,t)=p(x_1,t)p(x_2,0)\cdots p(x_N,t).
\end{equation}
Moreover, $p(x_1,t)=\left<P p(\bm x,t)\right>$.
Applying $P$ to the Liouville 
equation (\ref{Liouville}) and formally integrating 
out the orthogonal dynamics $Qp$ yields the Mori-Zwanzig (MZ) 
equation 
\cite{Venturi_MZ,VenturiBook} 
\begin{equation}
\frac{\partial p(x_1,t)}{\partial t}=
\left<PLP\right>p(x_1,t)+\int_{0}^t\left<PLe^{(t-s)QL}QL\right>p(x_1,s)ds.
\label{MZ}
\end{equation}
As demonstrated in \cite{Dominy2017}, there exists 
a duality between the MZ-PDF equation \eqref{MZ} and the 
more classical MZ formulation in the space of observables 
(see \cite{Gouasmi,Chorin}). Such duality is the same that 
pairs the Frobenious-Perron and the Koopman operators 
we discussed in Section \ref{sec:intro}. Equation \eqref{MZ} is 
very challenging  to solve. One of the main mathematical 
difficulties is the evaluation of the memory 
integral (second term at the right hand side). 
This term arises from purely formal mathematical 
manipulations, i.e., by using the variation 
of constants formula or the Dyson 
formula (see \cite{Yuan1,Venturi_MZ, Chorin}). 
Hence, the memory integral does not incorporate 
any information about the {\em structure} 
of the dynamical system \eqref{dyn}, i.e., it holds for 
 any Liouville operator $L$.

At this point, we notice that a comparison between 
equations \eqref{MZ} and \eqref{marginal_general} 
allows us to write the MZ memory integral 
as\footnote{This key observation 
allows us to represent the MZ memory integral in a way that 
depends on the specific dynamical system under consideration. In particular, 
if we have available an estimate of $p(x_1,t)$ and the conditional expectation 
$\mathbb{E}[G_1(\bm x(t))|x_1(t)]$ (see equation \eqref{MZapprox}), e.g., from sample paths, then we can easily estimate the MZ memory.}
\begin{equation}
\int_{0}^t\left<PLe^{(t-s)QL}QL\right>p(x_1,s)ds = -\left<PLP\right>p(x_1,t) -
\frac{\partial }{\partial x_1} \int_{-\infty}^\infty \cdots \int_{-\infty}^\infty 
\left(G_1(\bm x)p(\bm x,t)\right)dx_2\cdots dx_N.
\end{equation}
The streaming term $\left<PLP\right>p(x_1,t)$, 
can be generally expressed as  
\begin{align}
\left<PLP\right>p(x_1,t) =  
-\frac{\partial }{\partial x_1}\left( p(x_1,t) 
\int_{-\infty}^\infty \cdots \int_{-\infty}^\infty G_1(\bm{x}) 
\prod_{j=2}^N p(x_j,0)dx_j\right).
\end{align}
This implies that the MZ memory integral \eqref{MZ} can be 
expressed as
\begin{equation}
\int_{0}^t\left<PLe^{(t-s)QL}QL\right>p(x_1,s)ds  = 
-\frac{\partial }{\partial x_1}\left(p(x_1,t)M(x_1,t)\right),
\label{MZapprox}
\end{equation}
where the function $M(x_1,t)$ is
\begin{equation}
M(x_1,t)=\mathbb{E}\left[G_1(\bm x(t))|x_1(t)\right]-
\mathbb{E}\left[G_1(\bm x(0))|x_1(0)\right].
\label{M}
\end{equation}
Clearly, $\mathbb{E}\left[G_1(\bm x(0))|x_1(0)\right]$ is known, since  
the statistical properties of the initial state are assumed to be known. 
On the other hand, the term  $\mathbb{E}\left[G_1(\bm x(t))|x_1(t)\right]$, 
i.e., the conditional expectation of  $G_1(\bm x(t))$ given $x_1(t)$ 
is usually not known, and it cannot be computed based on the 
PDF $p(x_1,t)$ alone. However, it can be estimated from sample 
trajectories of \eqref{dyn} as we discussed in Section \ref{sec:estimating_conditional_expectations}.
The specific form of $M(x_1,t)$ depends on the structure of the 
dynamical system, in particular on the form of $G_1(\bm x)$. Let us provide
a simple example. 

\begin{figure}[t]
\centerline{
		\includegraphics[width=0.455\textwidth]{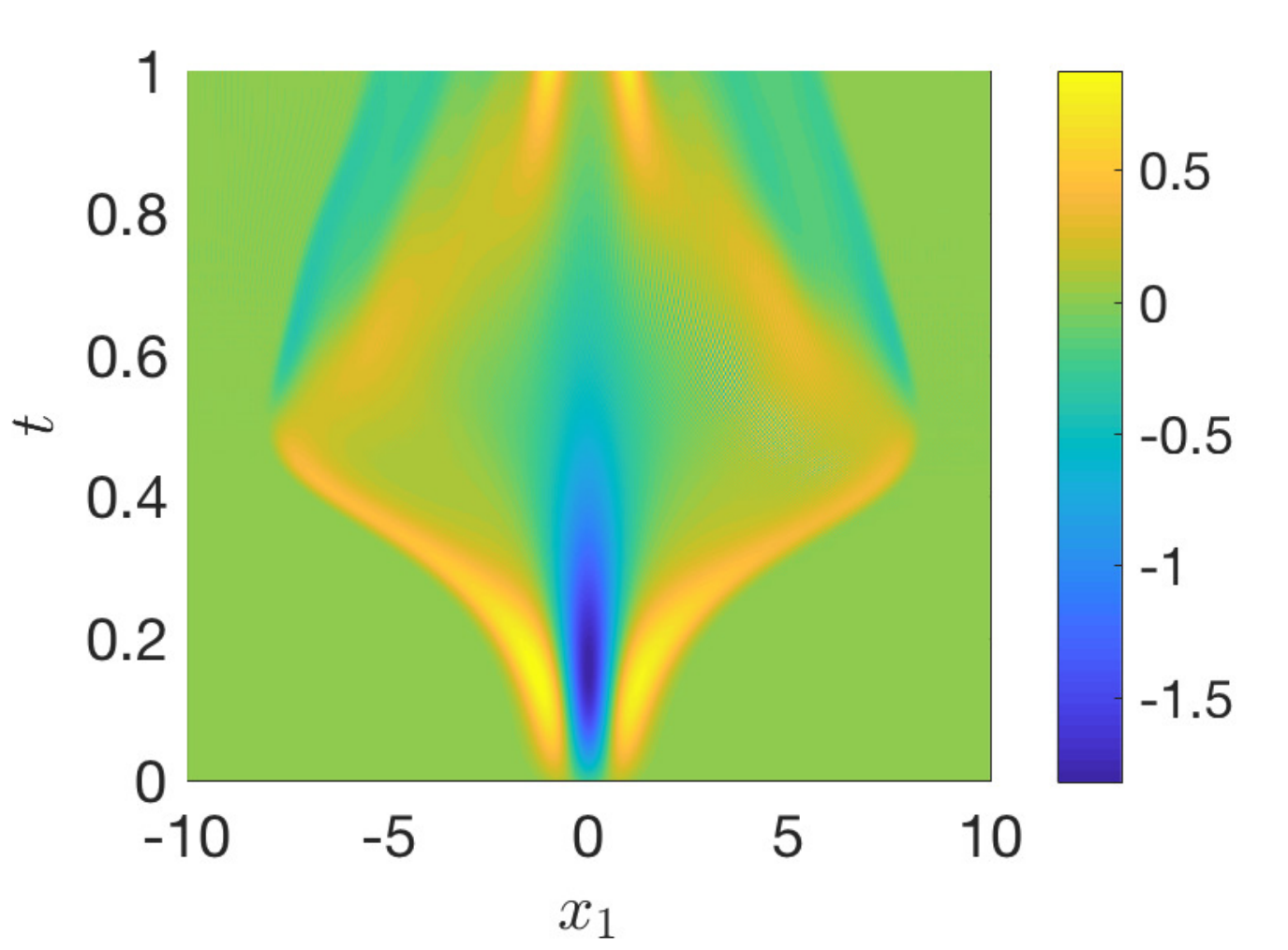}
		\includegraphics[width=0.4\textwidth]{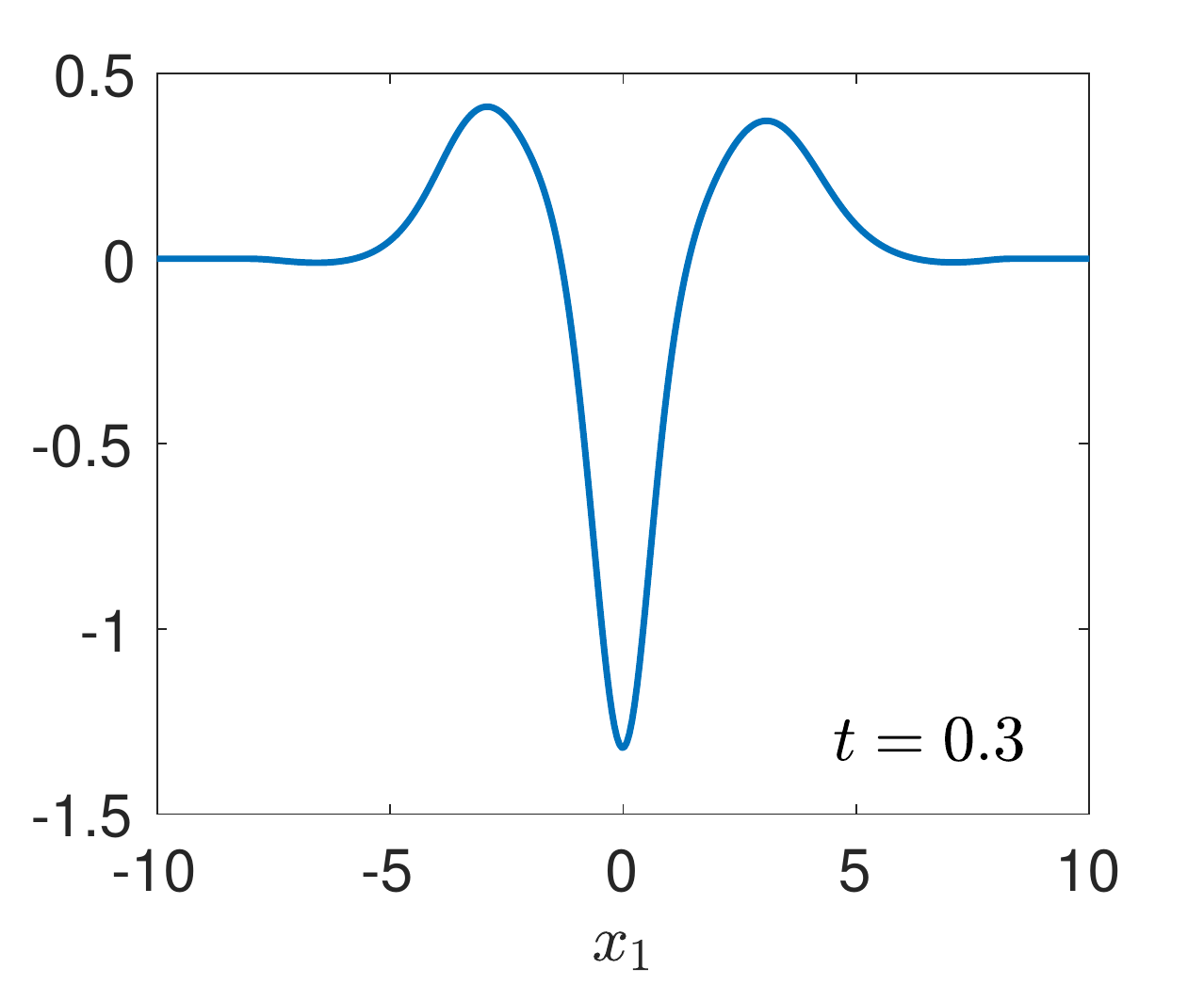}}
\caption{Kraichnan-Orszag three mode problem. 
Mori-Zwanzig memory integral appearing in the MZ-PDF 
equation \eqref{MZ}. The memory is computed based 
on $5000$ sample trajectories and equations 
\eqref{MZapprox}-\eqref{MZKO}. }
\label{fig:KO_MZ}
\end{figure}

{\example $\,$} Consider the Kraichnan-Orszag 
system \eqref{KO-dyn} evolving from a random initial 
state $\bm x(0)$ with i.i.d. components. A 
simple  calculation shows that the function \eqref{M} in 
this case takes the form 
\begin{equation}
M(x_1,t)=x_1\left(\mathbb{E}[x_3(t)|x_1(t)]-\mathbb{E}[x_3(0)]\right).
\label{MZKO}
\end{equation}
The conditional expectation $\mathbb{E}[x_3(t)|x_1(t)]$ 
can be estimated from sample trajectories  by using the methods 
we discussed in Section \ref{sec:estimating_conditional_expectations}. 
With estimates of $p(x_1,t)$ and $\mathbb{E}[x_3(t)|x_1(t)]$ available, 
it is easy to compute the Mori-Zwanzig memory integral \eqref{MZapprox}-\eqref{MZKO}. In Figure \ref{fig:KO_MZ} we plot the results we obtain with 
$5000$ sample trajectories.  

\section{Numerical examples}
\label{sec:numerics}
In this section we apply the mathematical methods 
presented in Sections \ref{sec:intro}-\ref{sec:data-drivenMZ}
to a high-dimensional nonlinear dynamical system and a drug-resistant 
malaria propagation model \cite{Meara,Ewungkem}.

\subsection{High-dimensional nonlinear dynamics} 
\label{venturi-16}
\begin{figure}[t]
\centerline{\hspace{-.2cm}(a)\hspace{6.6cm}(b)}
	\centering{
		\includegraphics[width=0.42\textwidth]{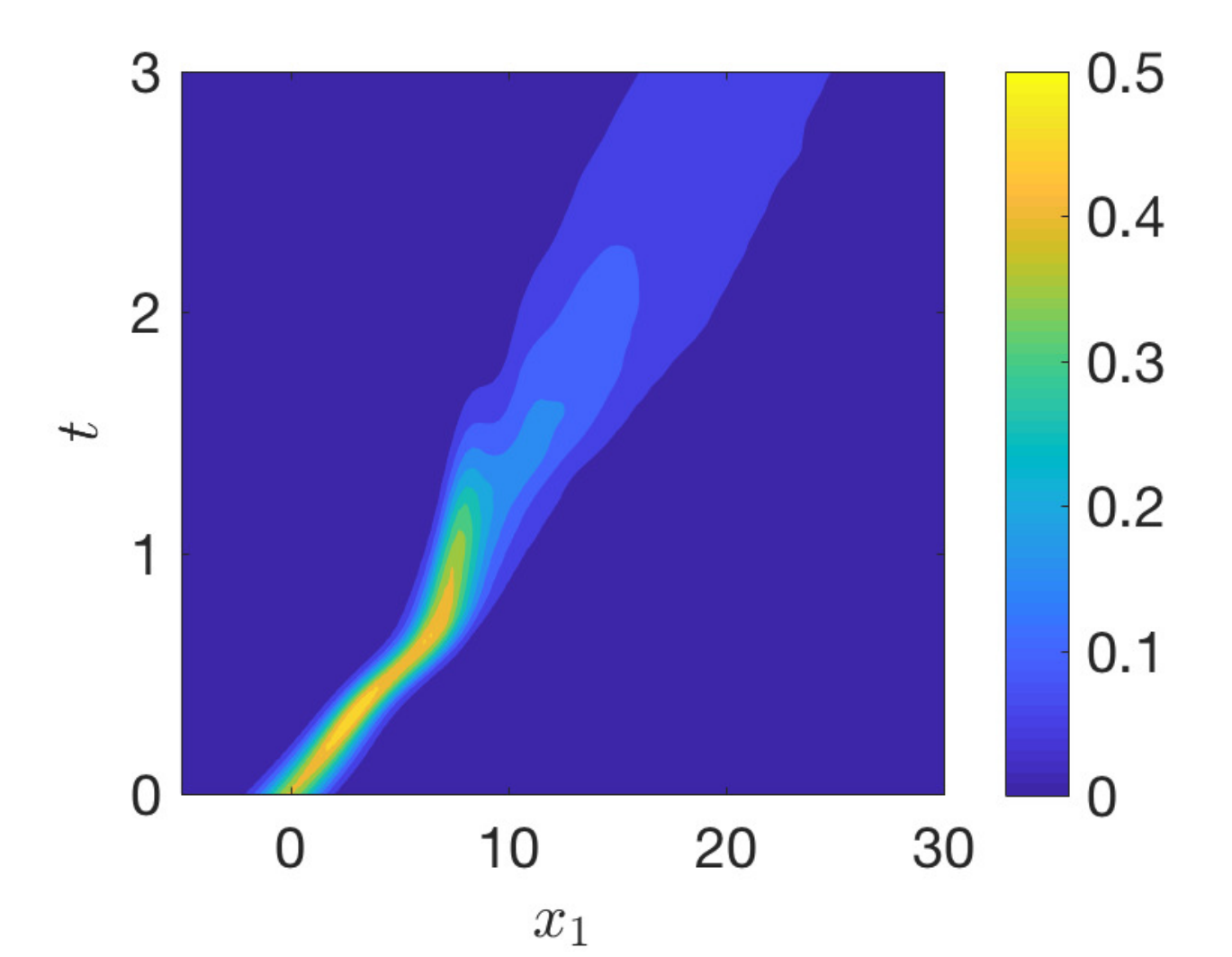}
		\includegraphics[width=0.42\textwidth]{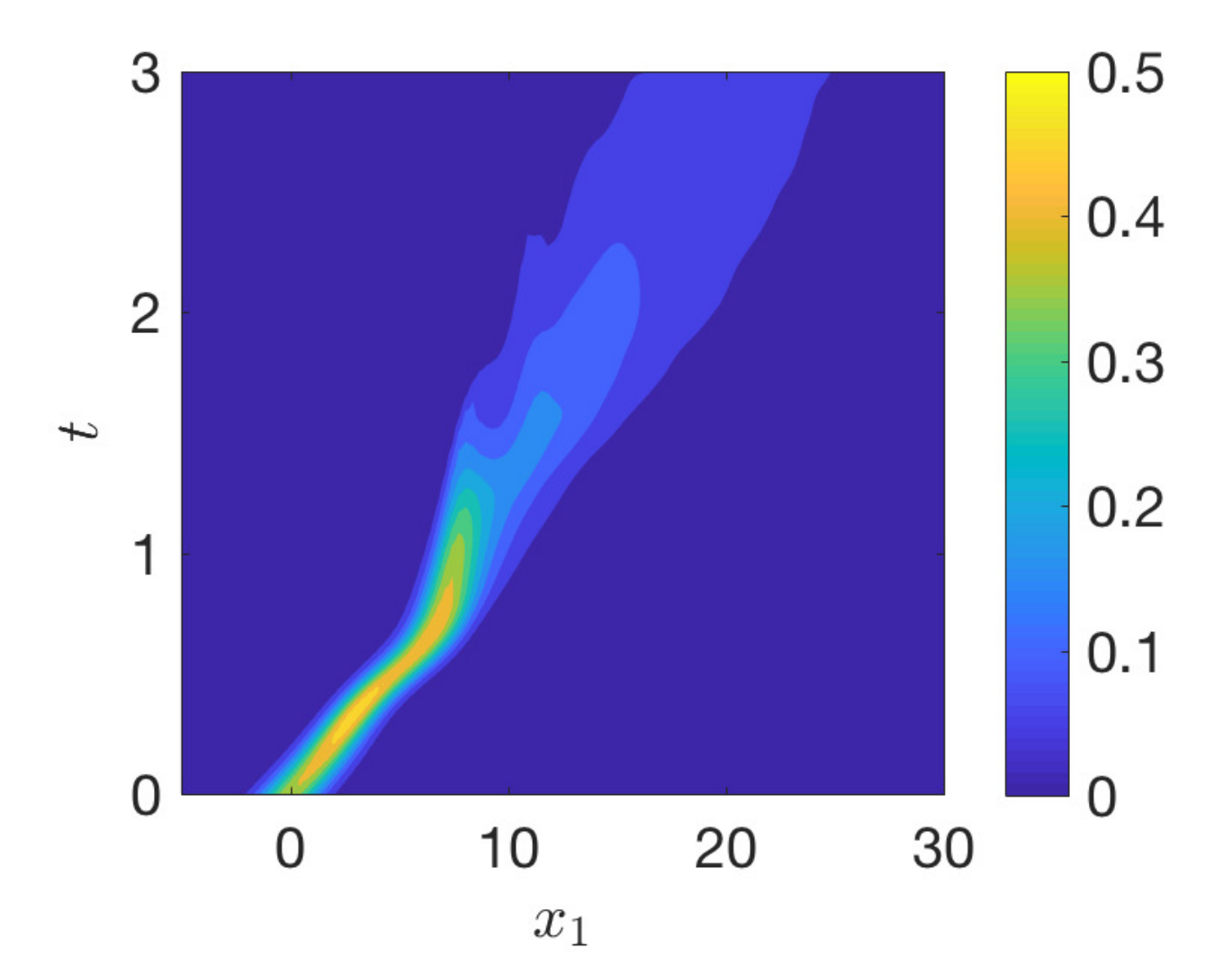}
}
\centering{
\includegraphics[width=0.28\textwidth]{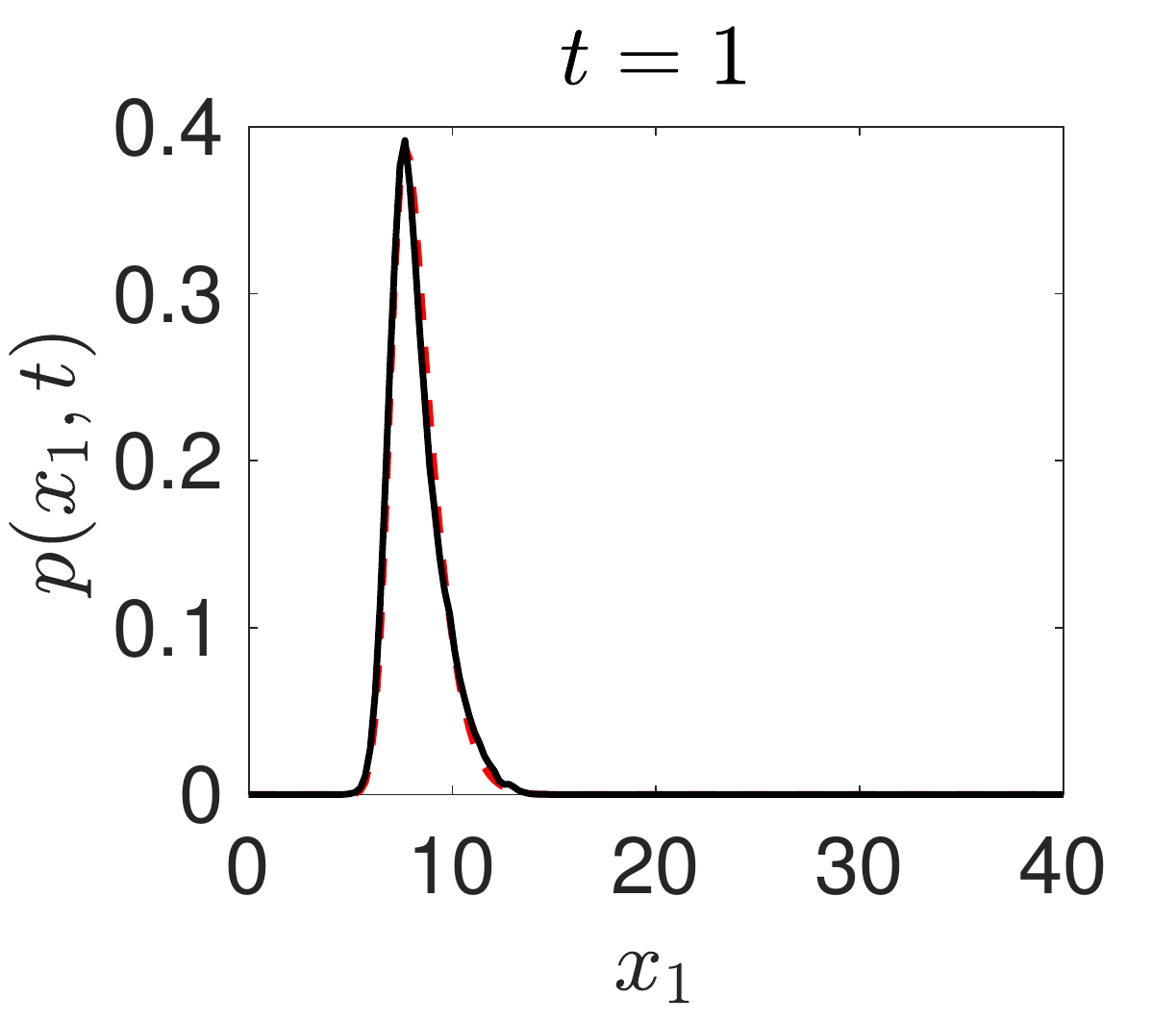}
\includegraphics[width=0.28\textwidth]{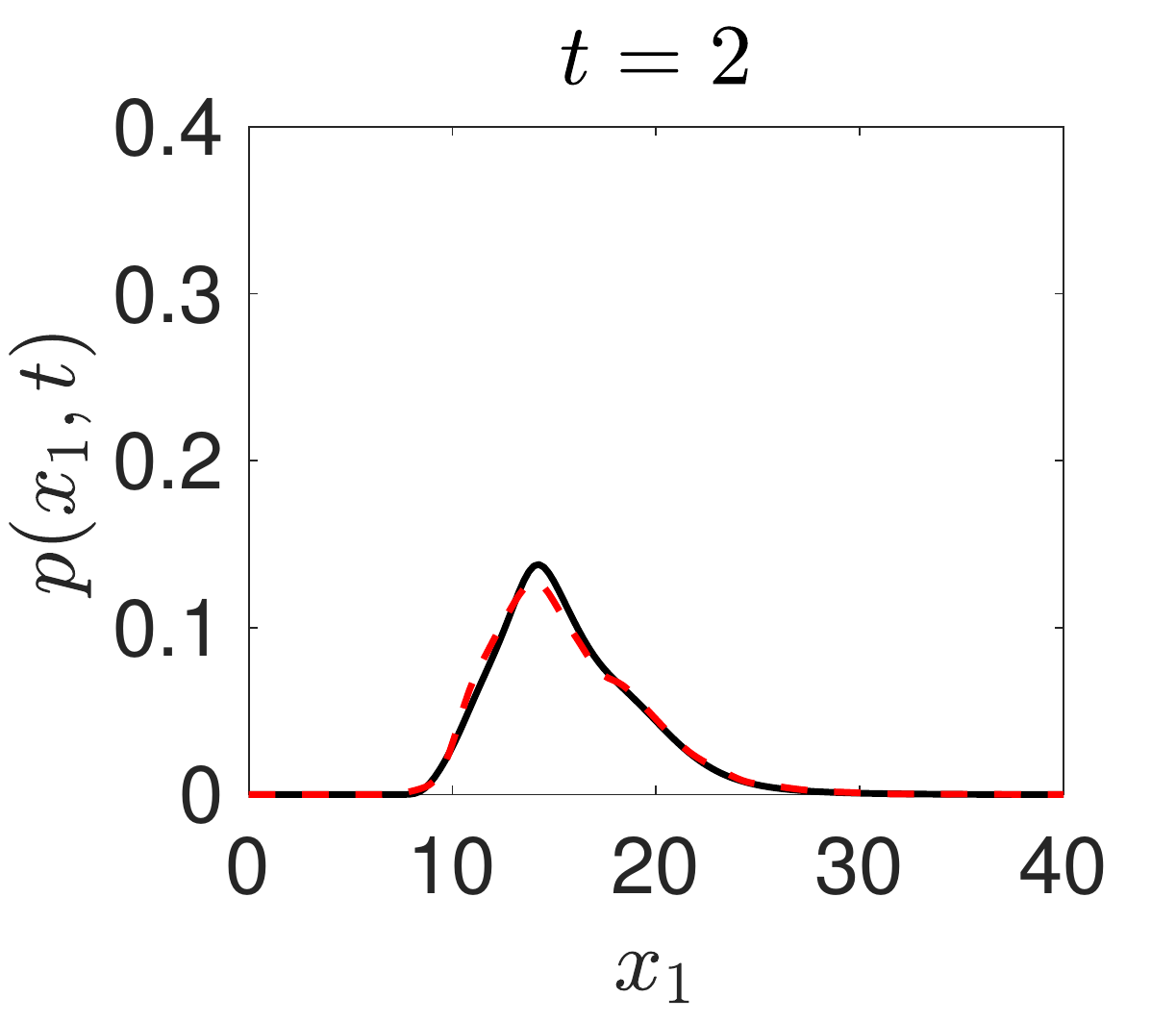}
\includegraphics[width=0.28\textwidth]{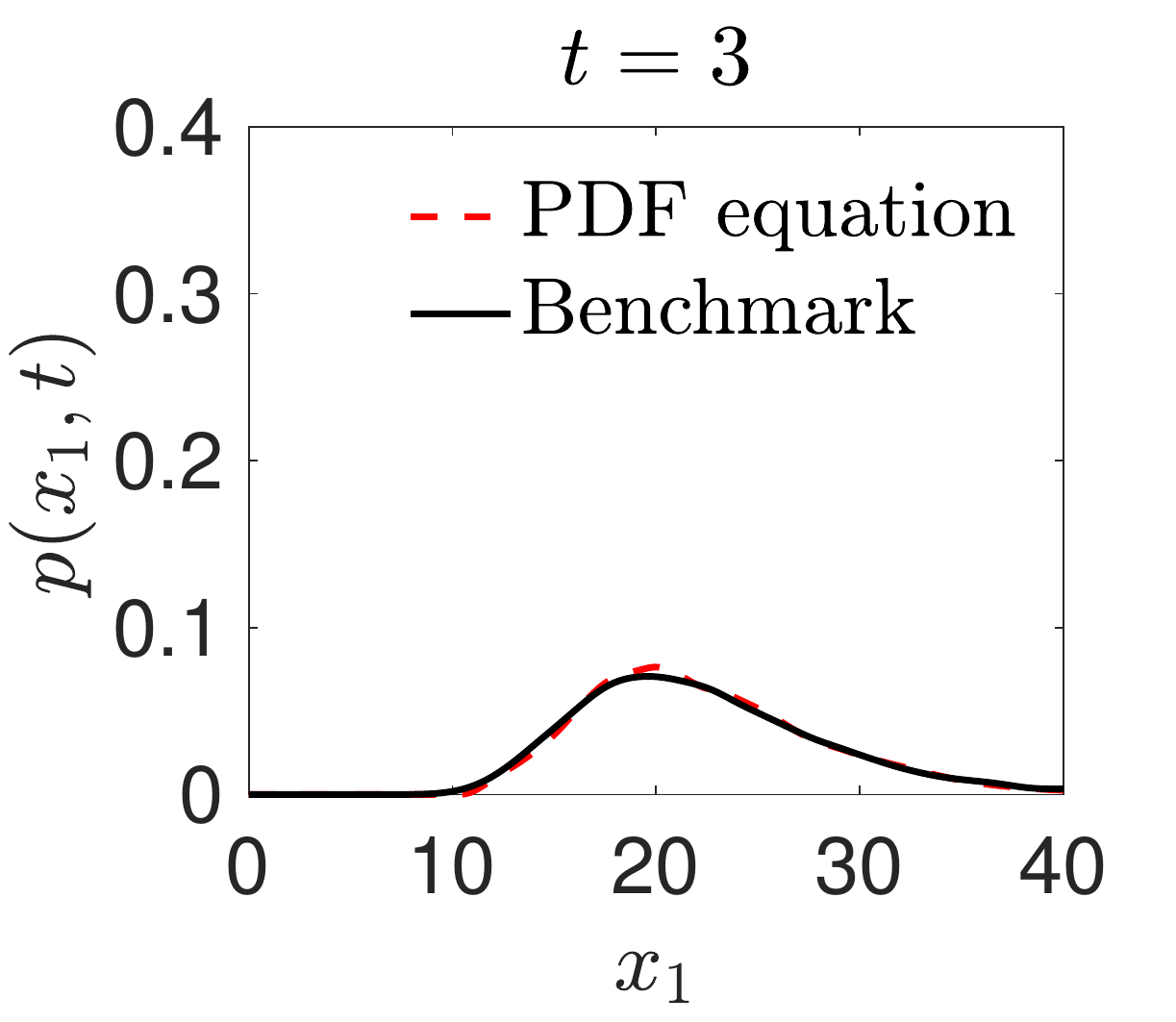}
}

\caption{Nonlinear dynamical system \eqref{D16}. (a) Accurate 
kernel density estimate \cite{Botev} 
of $p(x_1, t)$ based on $20000$ sample trajectories. 
(b) Data-driven solution of the transport equation \eqref{BBGKYHD3}. 
We estimated the conditional expectation 
$\mathbb{E}\left[\sin(x_2(t)) | x_1(t)\right]$ based 
on $5000$ sample trajectories of \eqref{D16} (see Figure \ref{fig:expectations}).}
\label{fig:PDF_D16}
\end{figure}

Consider the following $N$-dimensional nonlinear dynamical system
\begin{equation}
\frac{dx_i}{dt} =  -\sin(x_{i+1})x_i- Ax_i + F,  
\qquad i = 1,...,N,
\label{D16}
\end{equation} 
where $x_{n+1}(t)=x_1(t)$ (periodic boundary conditions). 
Depending on the value of $F$, $A$ and on the number 
of phase variables $N$, this system can exhibit 
different behaviors. Here we set $F=10$, $A=0.2$ 
and $N=1000$. The Liouville transport equation associated 
with \eqref{D16} is 
\begin{equation}
\frac{\partial p(\bm x,t)}{\partial t} = -\sum_{i=1}^{N}
\frac{\partial}{\partial x_i} 
\left[ \left(F-\sin(x_{i+1})x_i- Ax_i\right)  
p(\bm x,t) \right].
\label{BBGD16}
\end{equation}
This equation cannot be solved in a 
tensor product representation because of the high number 
of phase variables and possible lack of regularity of 
the solution. 
The evolution equation for the PDF of each phase variable 
$x_i(t)$ can be obtained by integrating \eqref{BBGD16} with 
respect to all other variables. This yields the unclosed equation
\begin{equation}
\frac{\partial p(x_i,t)}{\partial t}= -\frac{\partial}{\partial x_i} 
\int_{-\infty}^\infty\left[ \left(F-\sin(x_{i+1})x_i- Ax_i \right)  
p(x_i,x_{i+1},t)\right]dx_{i+1}.
\label{BBGKYHD1}
\end{equation}
We can write \eqref{BBGKYHD1} equivalently as
\begin{equation}
\frac{\partial p(x_i,t)}{\partial t}=-F\frac{\partial p(x_i,t)}{\partial x_i} 
+ A\frac{\partial (x_i p(x_i,t))}{\partial x_i}
-\frac{\partial}{\partial x_i} 
x_i \int_{-\infty}^\infty   \sin(x_{i+1}) p(x_i,x_{i+1},t)dx_{i+1}.
\label{BBGKYHD2}
\end{equation}
Note that all equations for $p(x_i,t)$ 
have the same structure, independently of $i$. This means 
that if the random initial state $\bm x_0$ has i.i.d. components, 
then the evolution of each $p(x_i,t)$ does not depend 
on $i$, i.e., it is the same for all $i=1,...,N$. A similar 
conclusion holds for the joint distributions $p(x_i,x_{i+1},t)$, 
which satisfy the equations
\begin{align}
	\frac{\partial p(x_i,x_{i+1},t)}{\partial t}=&
	-\frac{\partial}{\partial x_i} 
	\left[ \left(F-\sin(x_{i+1})x_i- Ax_i\right)  
	p(x_i,x_{i+1},t)\right]-\cdots \nonumber\\
	& \frac{\partial}{\partial x_{i+1}}\int_{-\infty}^\infty  
	\left[ \left(F-\sin(x_{i+2})x_{i+1}- Ax_{i+1}\right)  
	p(x_i,x_{i+1},x_{i+2},t)\right]dx_{i+2}.
	\label{HD-twopdf}
\end{align}
Without loss of generality, let us set $i=1$ in 
equation \eqref{BBGKYHD2} and express the 
integral in terms of the conditional 
expectation of $\sin(x_2(t))$ given $x_1(t)$. 
This yields 
\begin{equation}
\frac{\partial p(x_1,t)}{\partial t}= \frac{\partial}{\partial x_1}\left( 
x_1p(x_1,t) \mathbb{E}\left[\sin(x_2(t)) | x_1(t)\right]\right) + 
\frac{\partial}{\partial x_1} \left[(Ax_1 - F)p(x_1,t) \right],
\label{BBGKYHD3}
\end{equation}
where 
\begin{equation}
\mathbb{E}\left[\sin(x_2(t)) | x_1(t)\right]=\int_{-\infty}^\infty  \sin(x_{2}) p(x_{2}|x_1,t)dx_{2}.
\label{condpdfhd}
\end{equation}
The conditional expectation \eqref{condpdfhd}  
can be estimated from sample trajectories of \eqref{D16} by using 
the mathematical techniques discussed in Section
\ref{sec:estimating_conditional_expectations}. The results are 
summarized in Figure \ref{fig:expectations}. With the conditional 
expectation available, we can compute the numerical solution 
of \eqref{BBGKYHD3} with a Fourier spectral method 
and compare it with an accurate kernel density 
benchmark estimate. This is done in Figure \ref{fig:PDF_D16}. 

\begin{figure}[t]
\centerline{\hspace{0.5cm}(a)\hspace{7cm}(b)}
\centerline{
\includegraphics[width=0.45\textwidth]{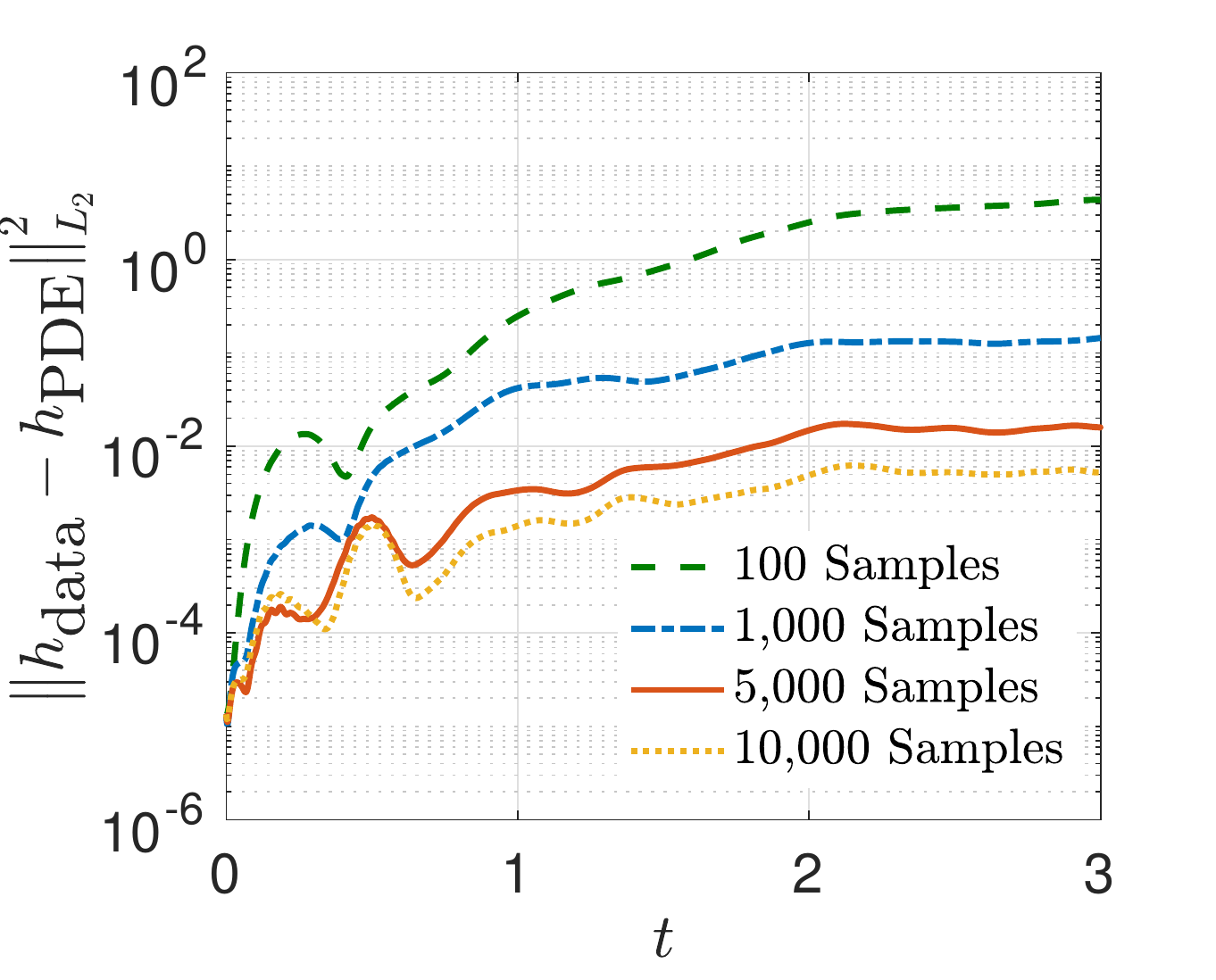}
\includegraphics[width=0.45\textwidth]{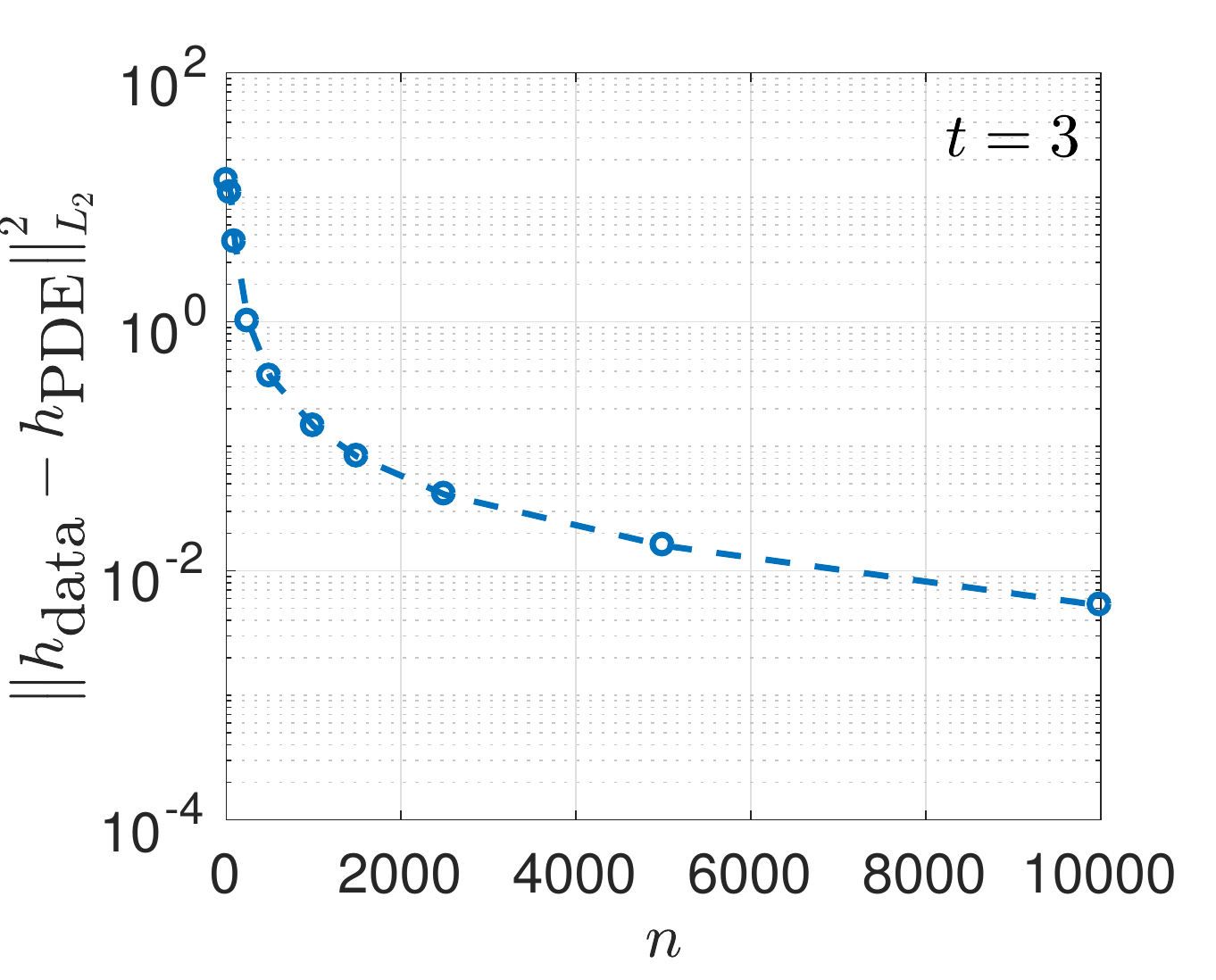}
}
\caption{Nonlinear dynamical system \eqref{D16}. (a) time dependent errors 
in \eqref{h1D16} for a variety of sample sizes. (b) 
Error decay at $t=3$ versus the number of sample trajectories.}
\label{fig:information_D16}
\end{figure}

Alongside the data-driven closure of equation \eqref{BBGKYHD3},
we also studied closures based on a system of equation similar to 
the one we derived in Section \ref{sec:information_content}(a). 
In this case we obtain the hyperbolic system 
\begin{align}
\frac{\partial p(x_1,t)}{\partial t} =&\frac{\partial}{\partial x_1}\left( 
x_1h(x_1,t)\right) + 
\frac{\partial}{\partial x_1} \left[(Ax_1 - F)p(x_1,t) \right],
\label{p5}\\
\frac{\partial h(x_1,t)}{\partial t}  = &  -\frac{\partial}{\partial x_1}  
\left(p(x_1, t) \mathbb{E}\left[\sin(x_2(t))\left(F - \sin(x_2(t))x_1(t) - Ax_1(t)\right)| x_1(t) \right]\right) +\cdots\nonumber \\ 	
&  p(x_1,t) \mathbb{E}\left[ \cos(x_2(t))\left(F-\sin(x_{3}(t))x_{2}(t)- 
Ax_{2}(t)\right) | x_1(t)\right],
\label{u_D16}
\end{align}
where 
\begin{equation}
h(x_1,t)=\int_{-\infty}^{\infty}\sin(x_2)p(x_1,x_2,t)dx_2.
\label{h1D16}
\end{equation}
The numerical solution to  \eqref{p5}-\eqref{u_D16} 
allows us to to measure the information content of the sample 
trajectories we employed to compute the closure approximation 
(see Section \ref{sec:information_content}(b)). 
In Figure \ref{fig:information_D16}, we plot the 
time-dependent error between the function 
\eqref{h1D16} we obtain from data and the 
solution to the system \eqref{p5}-\eqref{u_D16}. Observe that the error decreases as we increase the number of sample trajectories, suggesting that the information content increases with sample size.  
\begin{figure}[t]
\centering{
	\includegraphics[width=0.43\textwidth]{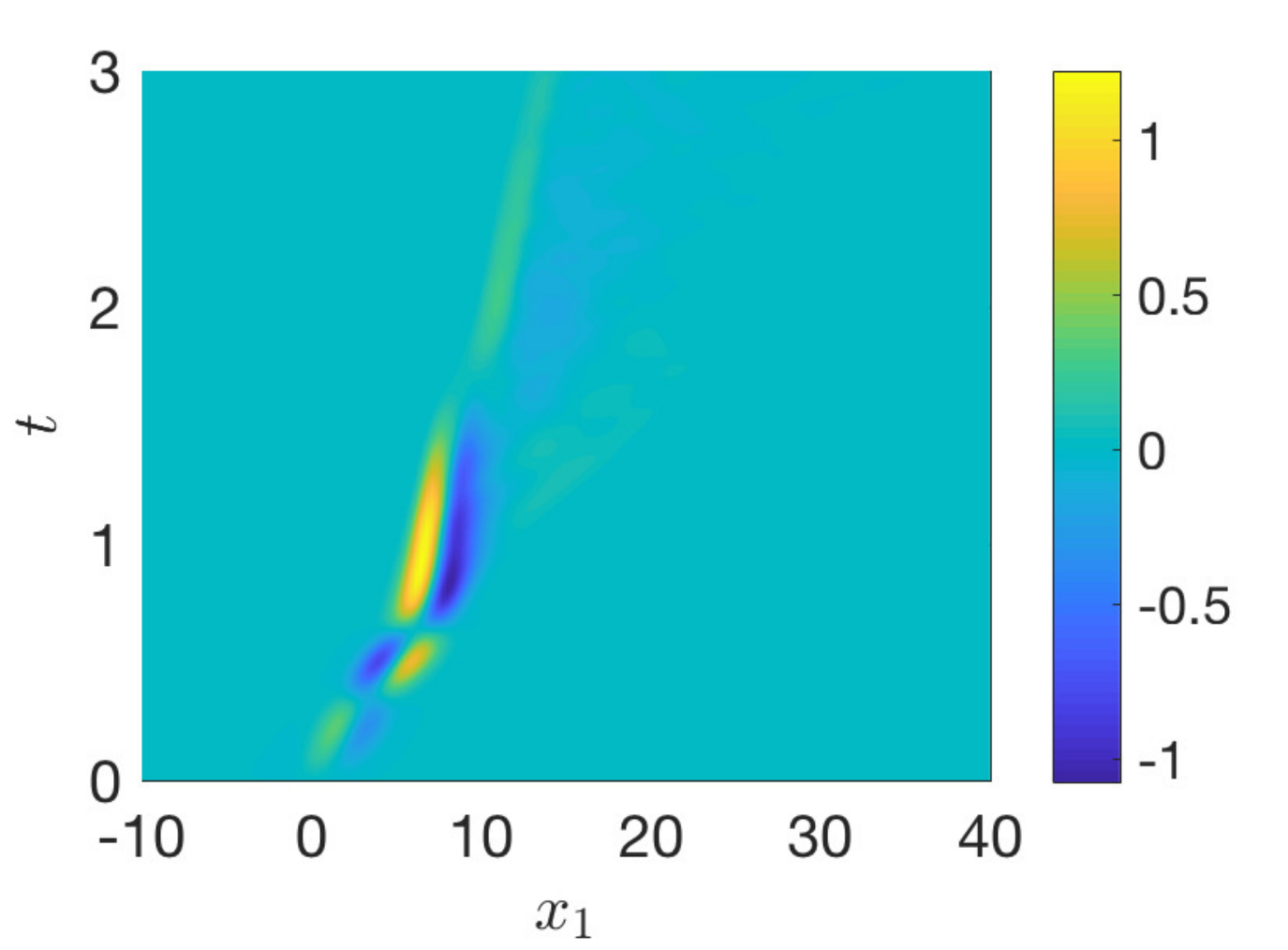}
	\includegraphics[width=0.34\textwidth]{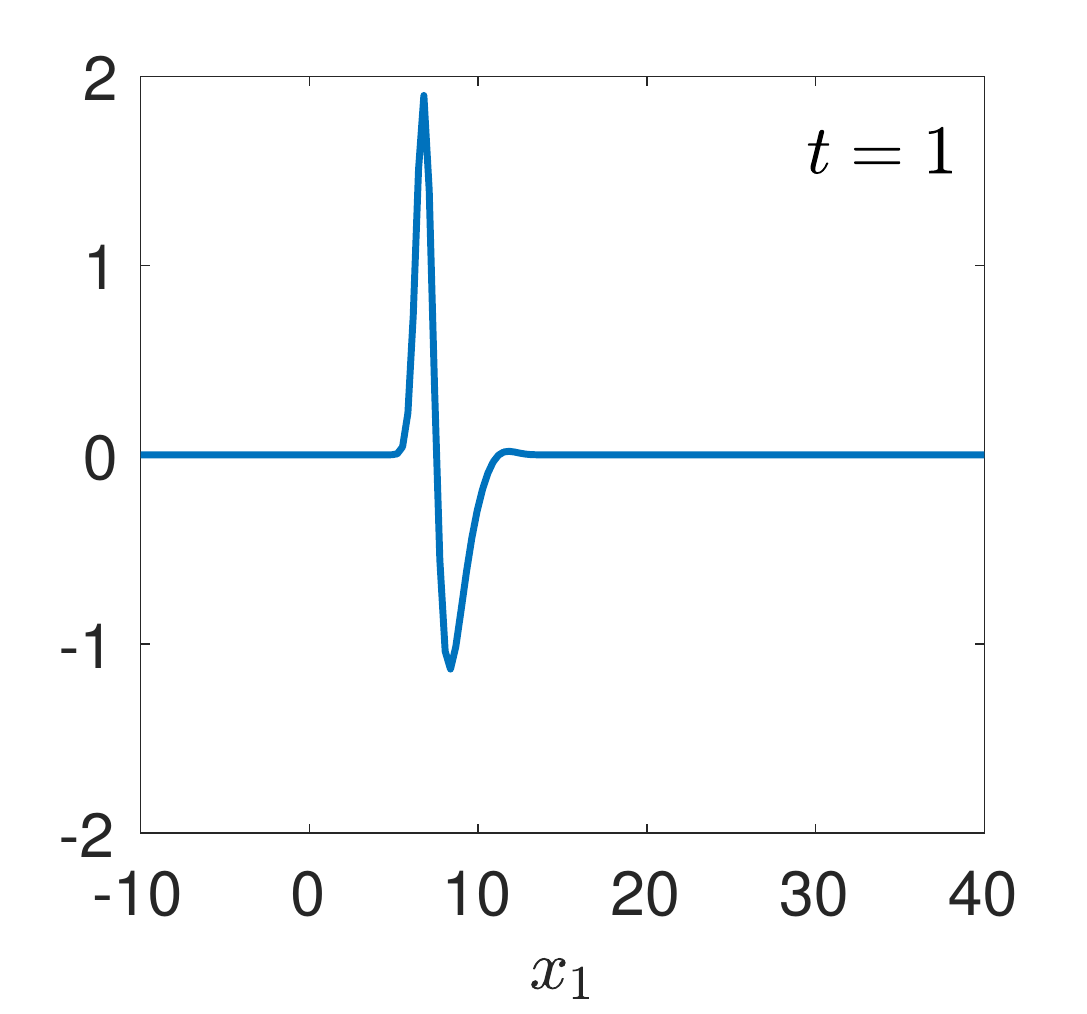}}
	
\caption{Mori-Zwanzig memory integral associated with the exact PDF 
equation for the first component of dynamical system \eqref{D16} 
in $1000$ dimensions. The initial state here is set as i.i.d. Gaussian.}
\label{fig:D16MZ}
\end{figure} 
With accurate estimates for $p(x_1,t)$ and 
$\mathbb{E}[G_1(\bm x(t))|x_1(t)]$ 
available (here $G_1(\bm x)=\sin(x_2)x_1-Ax_1+F$),  
we can immediately 
compute the MZ memory integral by using equations 
\eqref{MZapprox}-\eqref{M}. The results we obtain
are shown in Figure \ref{fig:D16MZ}.

\subsection{Drug resistant malaria propagation model}
The following dynamical system was proposed 
in \cite{Ewungkem,Meara} to model efficacy of intermittent 
preventative treatment (IPT) for battling drug resistant malaria. 
\begin{align} 
\begin{cases}
\displaystyle\frac{dS}{dt} = \mu_h(1 - S) - \beta_h S (M_s + kM_r) - qcS + \sigma(1 - \xi)(I_a + J_a) + rT_a(1-b)  + rT + wR   \vse\\
\displaystyle	\frac{dI_s}{dt} = \lambda \beta_h M_s S + \nu I_a - I_s(pa + \sigma + \mu_h)  \vse\\
\displaystyle	\frac{dI_a}{dt}  =  \beta_h M_s S(1 - \lambda) - I_a ( qc +  \nu + \sigma + \mu_h)  \vse\\
\displaystyle	\frac{dJ_s}{dt}  = \lambda k \beta_h M_r (S + \tau T_s + \tau T +  \tau T_a) +  v J_a - J_s(\sigma + \mu_h)   \vse\\
\displaystyle	\frac{dJ_a}{dt}  = k \beta_h M_r (1-\lambda)(S + \tau T_s + \tau T +  \tau T_a) - J_a(\sigma + \nu + \mu_h) \vse \\
\displaystyle	\frac{dT_s}{dt} = pa I_s - T_s(r + \tau k \beta_h M_r + \mu_h) \vse \\
\displaystyle	\frac{dT}{dt}  = qcS - T(r + \tau k \beta_h M_r + \mu_h ) \vse \\
\displaystyle	\frac{dT_a}{dt}  = qcI_a - T_a(r + \tau k \beta_h M_r + \mu_h)  \vse\\
\displaystyle	\frac{dR}{dt}  = rT_s + brT_a + \xi \sigma(I_a + J_a) + \sigma I_s +  \sigma J_s - R(w + \mu_h) \vse
	\end{cases}
\label{eq:malaria_dyn}
\end{align}	
Each phase variable here represents a proportion of the human population 
that is exposed to the virus (see Table \ref{tab:phase_variables}). 
The remaining equations in the system govern the dynamics of the 
mosquito populations that are responsible for spreading the virus  
\begin{align}
\begin{cases}
\displaystyle\frac{dM_s}{dt }  = \beta_m(1 - M_s - M_r)(I_a + I_s) - \mu_m M_s, \vse\\ 
\displaystyle\frac{dM_r}{dt }  = k \beta_m( 1-M_s - M_r)(J_a + J_s) - \mu_m M_r.
\end{cases}
\label{eq:malaria_dyn1}
\end{align}
$M_s$ here represents the proportion of the mosquito 
population infected by the sensitive strain of malaria 
while $M_r$ represents the proportion infected by the 
resistant strain. We will use a clearance time of $6$ days 
in correspondence with the relatively short half life of popular 
anti-malarial drug chlorproguanil-dapsone (CPG-DDS) \cite{Meara}. 

\begin{table}
\centering
\renewcommand{\arraystretch}{1.5}%
	\begin{tabular}{c | c |  l  } \hline
		Species& 	Variable & Description \\ \hline
		\multirow{2}{*}{\includegraphics[height=1cm]{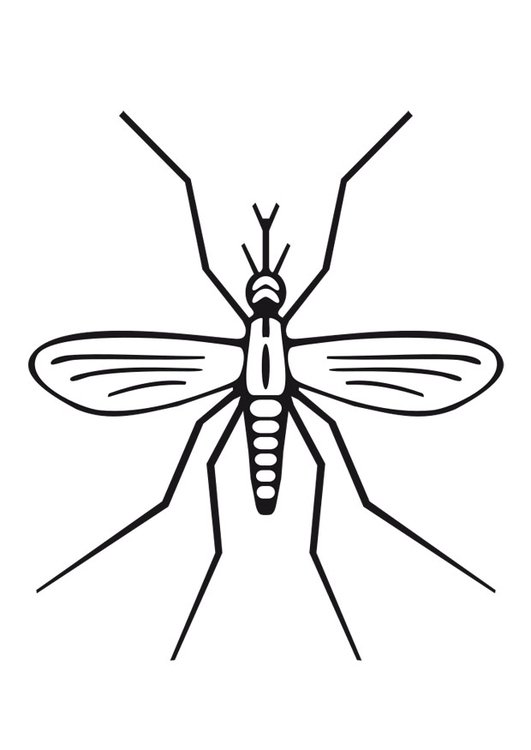}} &$M_s$ & Proportion of mosquitoes infected with the weak strain\\
		& $M_r$ & Proportion of mosquitoes infected with the strong strain \\ \hline
		\multirow{9}{*}{\includegraphics[height=0.8cm]{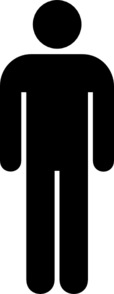}} &	$S$   & Proportion of susceptible humans \\ 
		&$I_s$ & Proportion of symptomatic humans infected with the weak strain \\ 
		&$I_a$ & Proportion of asymptomatic humans infected with the weak strain \\   
		&$J_s$ &  Proportion of symptomatic humans infected with the strong strain \\ 
		&$J_a$ & Proportion of asymptomatic humans infected with the strong strain \\  
		&$T_s$ & Treated symptomatic infectious humans\\ 
		&$T$   & IPT treated susceptible humans \\ 
		&$T_a$   & IPT treated asymptomatic infected humans \\  
		&$R$ & Proportion humans with temporary immunity\\  \hline 
	\end{tabular}
\caption{Definition of all phase variables appearing in the nonlinear system \eqref{eq:malaria_dyn}-\eqref{eq:malaria_dyn1}.}
\label{tab:phase_variables}
\end{table}	
The dynamics defined by \eqref{eq:malaria_dyn}-\eqref{eq:malaria_dyn1} 
can be described as follows: Once a human has been infected, they move from the 
susceptible class, $S$, to one of the infected classes. The infection can be cleared 
naturally or via IPT. The parameter $\xi$ represents the proportion of infections that 
will clear naturally and $\sigma$ is the rate at which the infection will clear. 
The efficacy of IPT depends on the length of time the treatment remains effective 
after administration, $r$, and the factor by which treatment impacts 
transmissability, $\tau$. The number of treatments per day, $c$, is scaled by a 
factor, $q$, to account for the fact that typically only children recieve IPT. 
All of the humans in groups $T$, $T_s$, and $T_a$ have been treated with a drug that is effective against weak strain of malaria, but not against the strong strain.	
Once an infection is cleared, an individual can move into the temporarily immune 
class, $R$, or return to the susceptible class, $S$. The rate at which immunity is lost 
is modeled in the parameter $w$. When members of the temporarily immune 
class lose their immunity, they return to the susceptible class. The parameters 
$\beta_m$ and $\beta_h$ are transmission rates from humans to parasites and vice 
versa,  while $\mu_m$ and $\mu_h$ represent death rates for mosquitoes and 
humans respectively. The transmission rates are multiplied by a 
reduction factor, $k$, to account for the presence of the resistant strain. The 
quantity $p\in[0,1]$ measures treatment efficacy, $\nu$ is the rate at 
which an asymptomatic infection progresses to a symptomatic one and $\lambda$ 
represents the proportion of infections that will be asymptomatic. Finally, 
$b$ represents the proportion of the class $T_a$ that will gain temporary 
immunity, while $(1-b)$ represents the proportion that will return to the 
susceptible population. A thorough parametric study of 
\eqref{eq:malaria_dyn}-\eqref{eq:malaria_dyn1} is 
presented in \cite{Ewungkem}. 

\begin{figure}[t]
\centerline{\hspace{-0.6cm}(a)\hspace{6.4cm}(b)}
\centering{
		\includegraphics[width=0.40\textwidth]{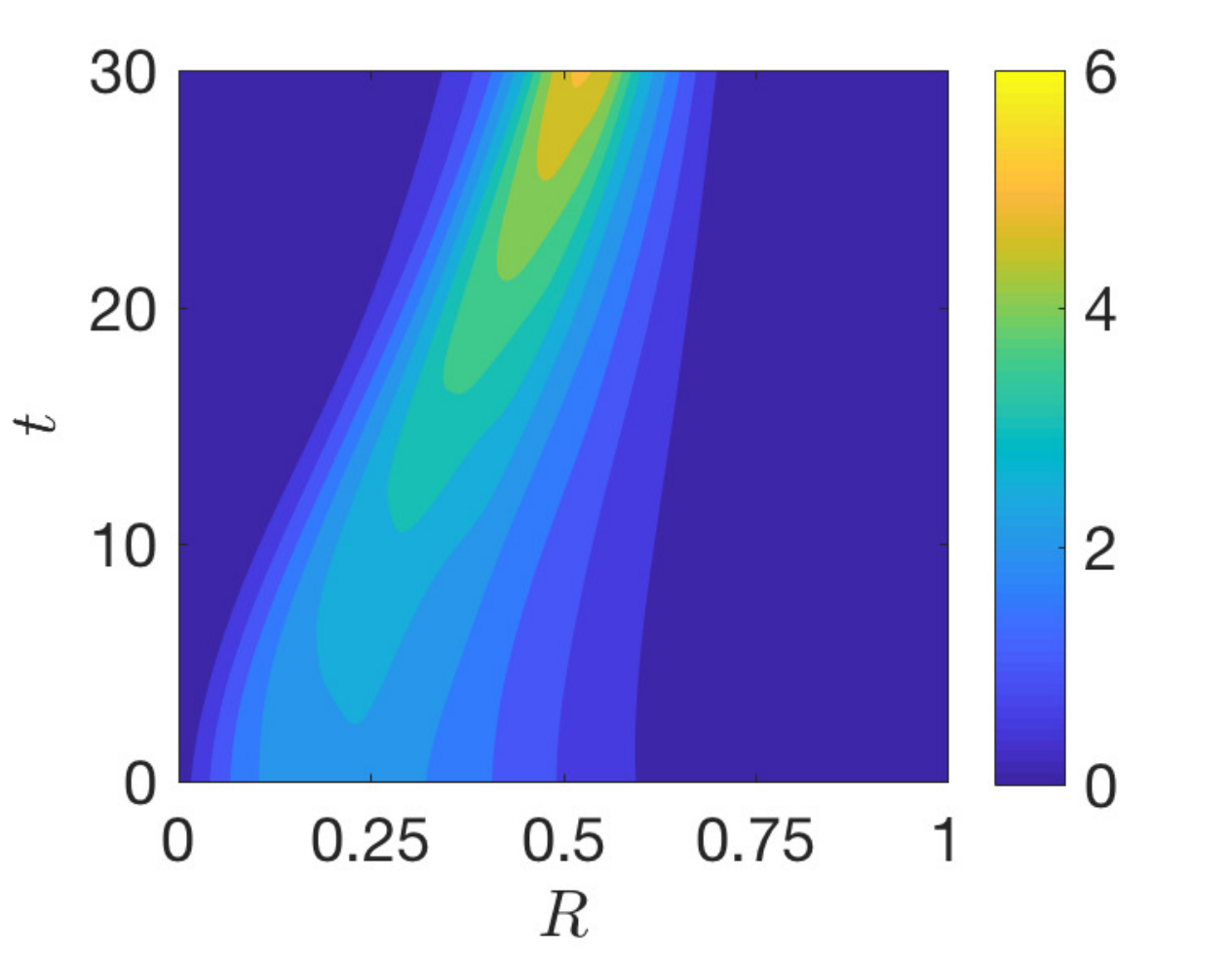}
		\includegraphics[width=0.40\textwidth]{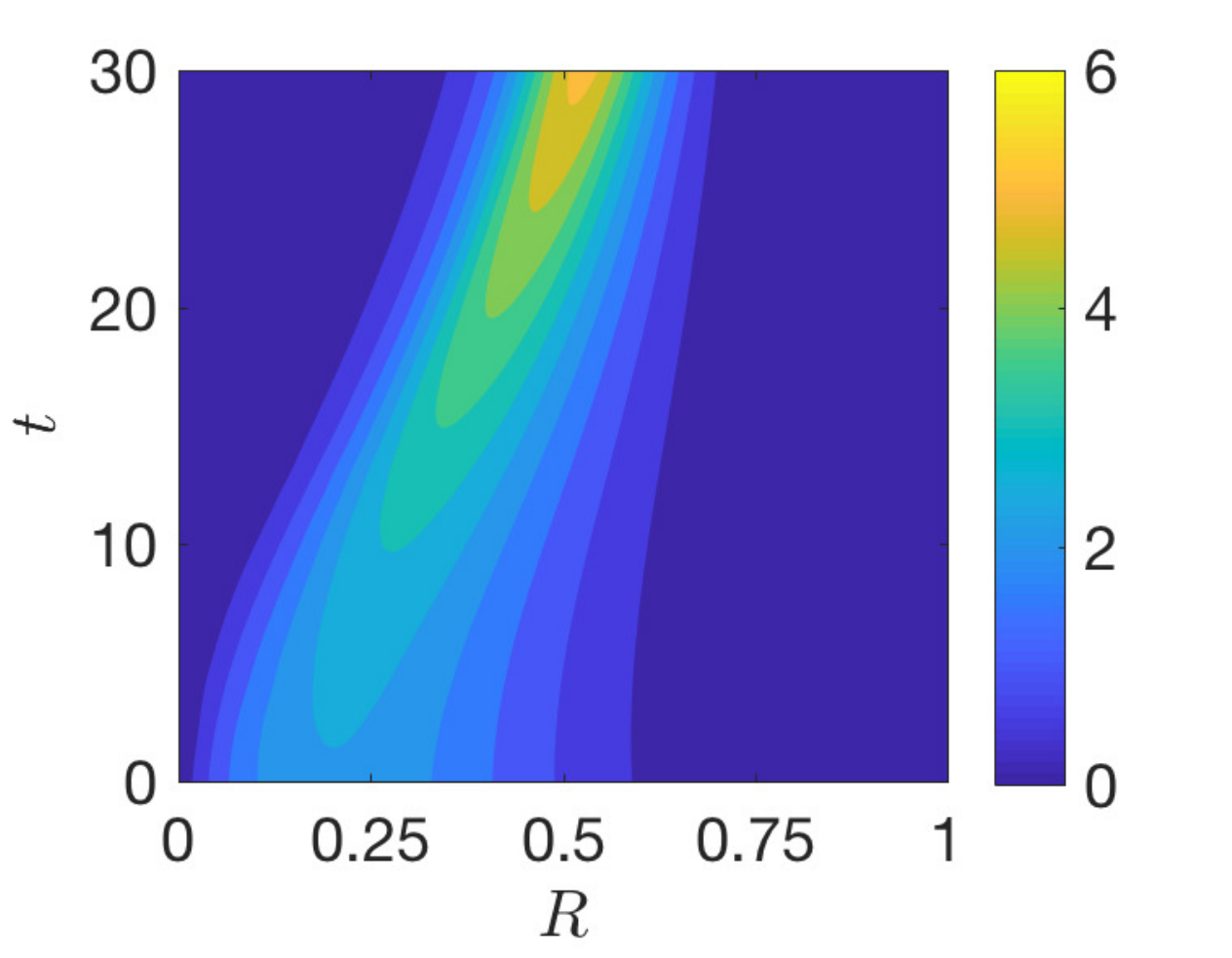}}
	
\centerline{
\includegraphics[width=0.28\textwidth]{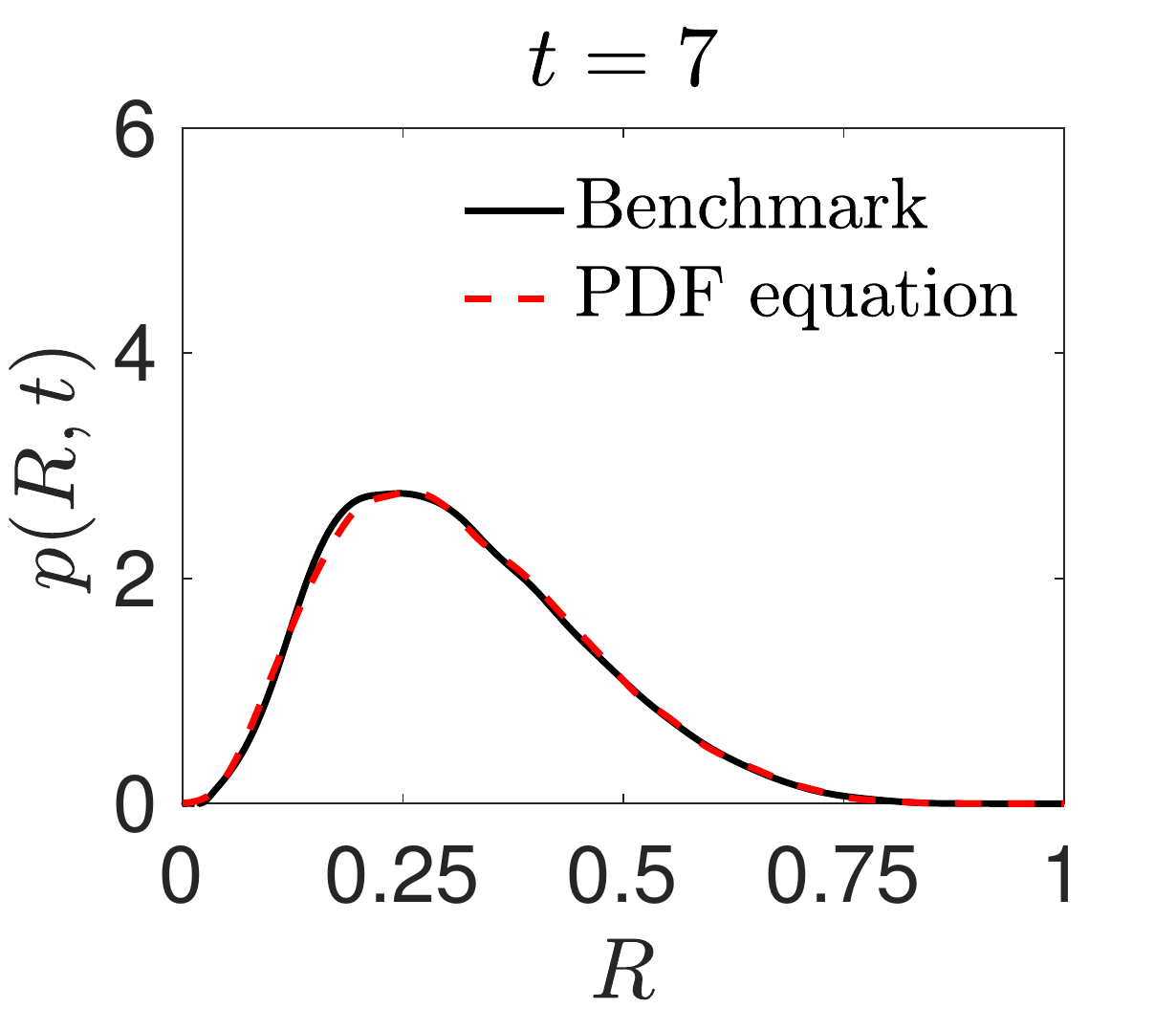}
\includegraphics[width=0.28\textwidth]{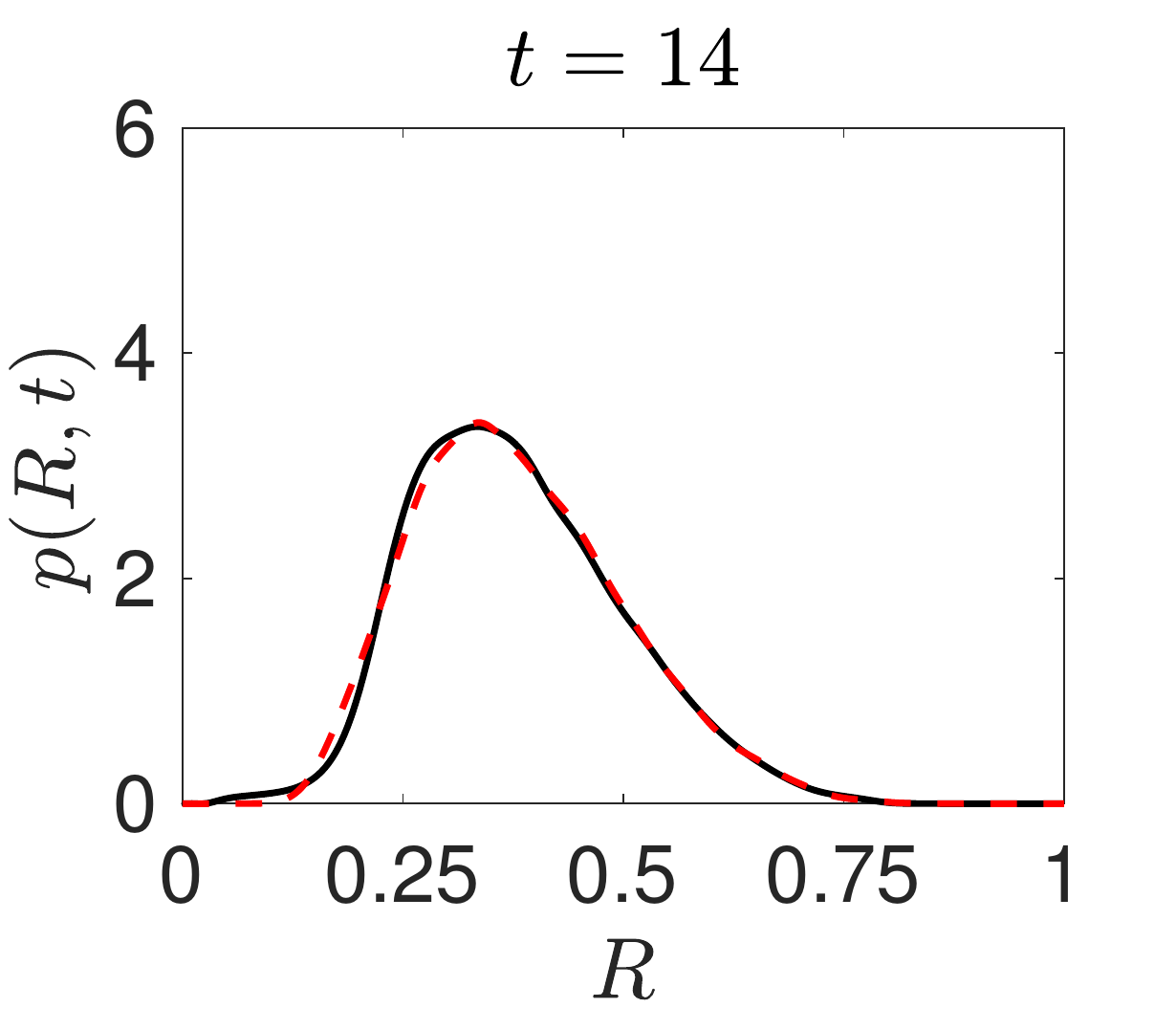}
\includegraphics[width=0.28\textwidth]{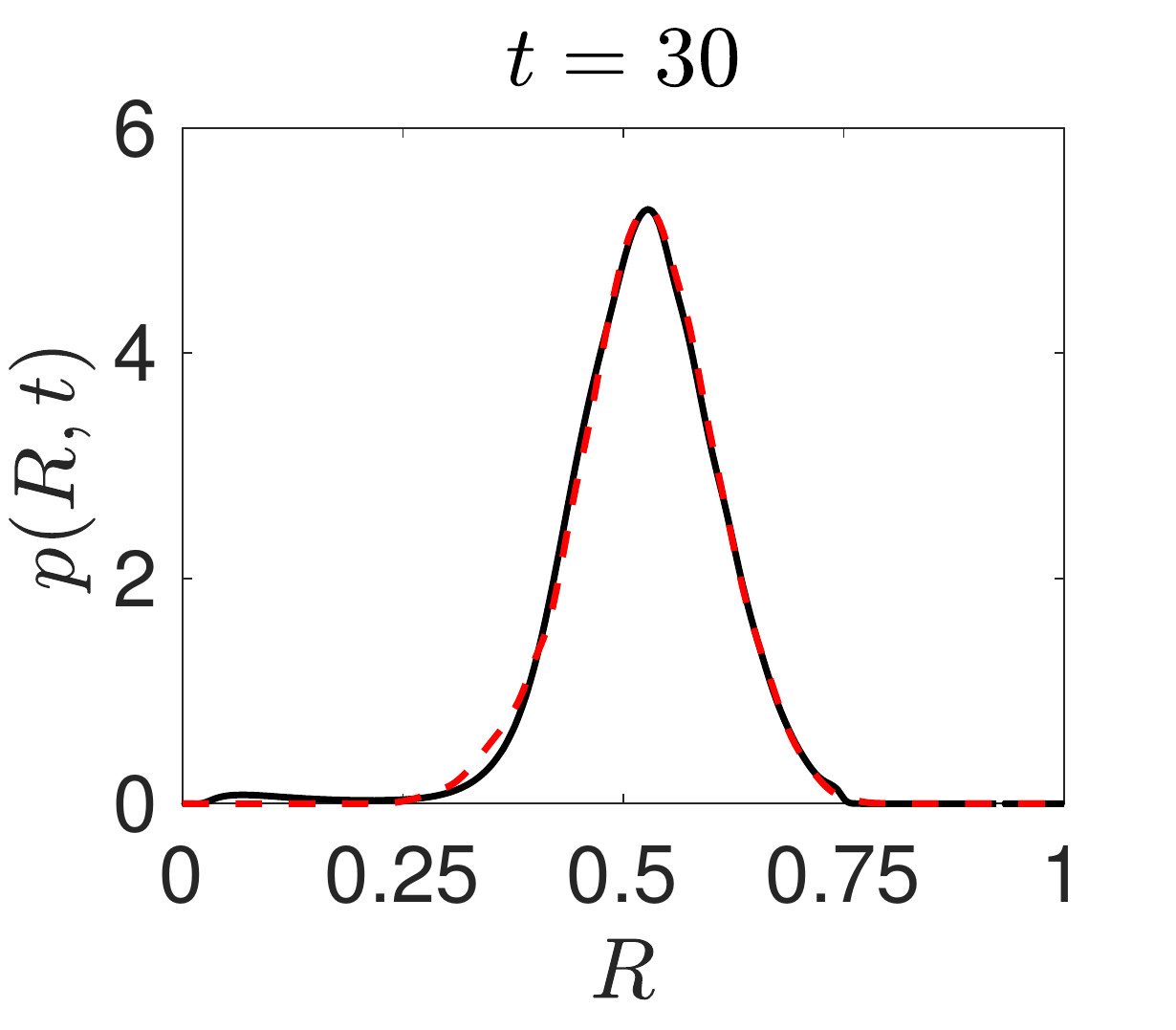}
}
\caption{Drug resistant malaria propagation model. 
(a) Accurate kernel density estimate of $p(R, t)$ 
(humans with temporary immunity) based on $30000$ 
sample trajectories. (b)  Numerical solution of 
\eqref{malaria_pdf}-\eqref{hmalaria} 
obtained estimating the conditional expectation in \eqref{hmalaria} 
with $5000$ sample trajectories. We also plot the time 
snapshots of the PDF $p(R,t)$ at $7$ days, $14$ days, and $30$ days.}.
\label{fig:malariaPDF}
\end{figure}

\begin{figure}[t]
\centerline{\hspace{0.4cm}(a)\hspace{7cm}(b)}
\centerline{
\includegraphics[width=0.45\textwidth]{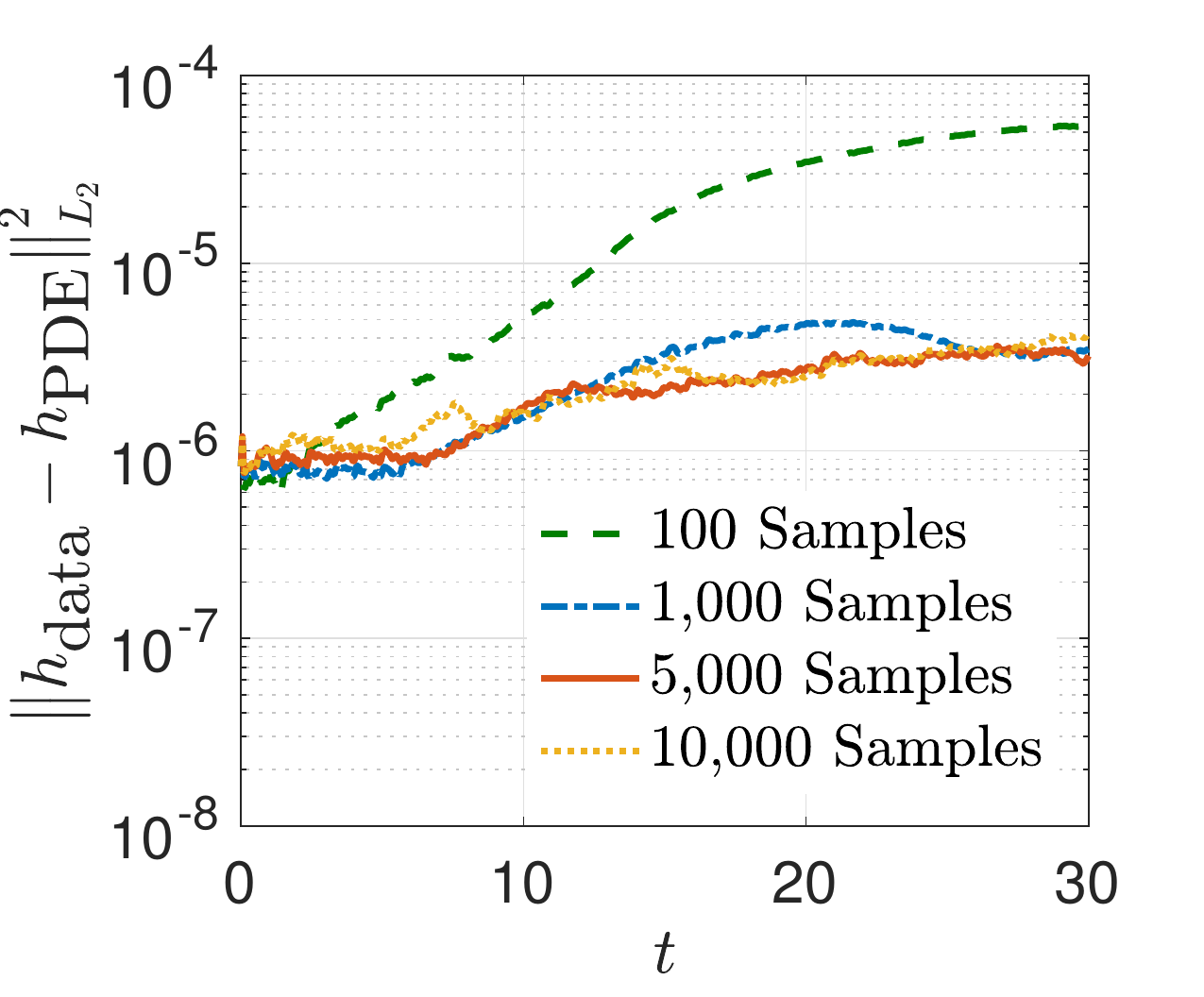}
\includegraphics[width=0.45\textwidth]{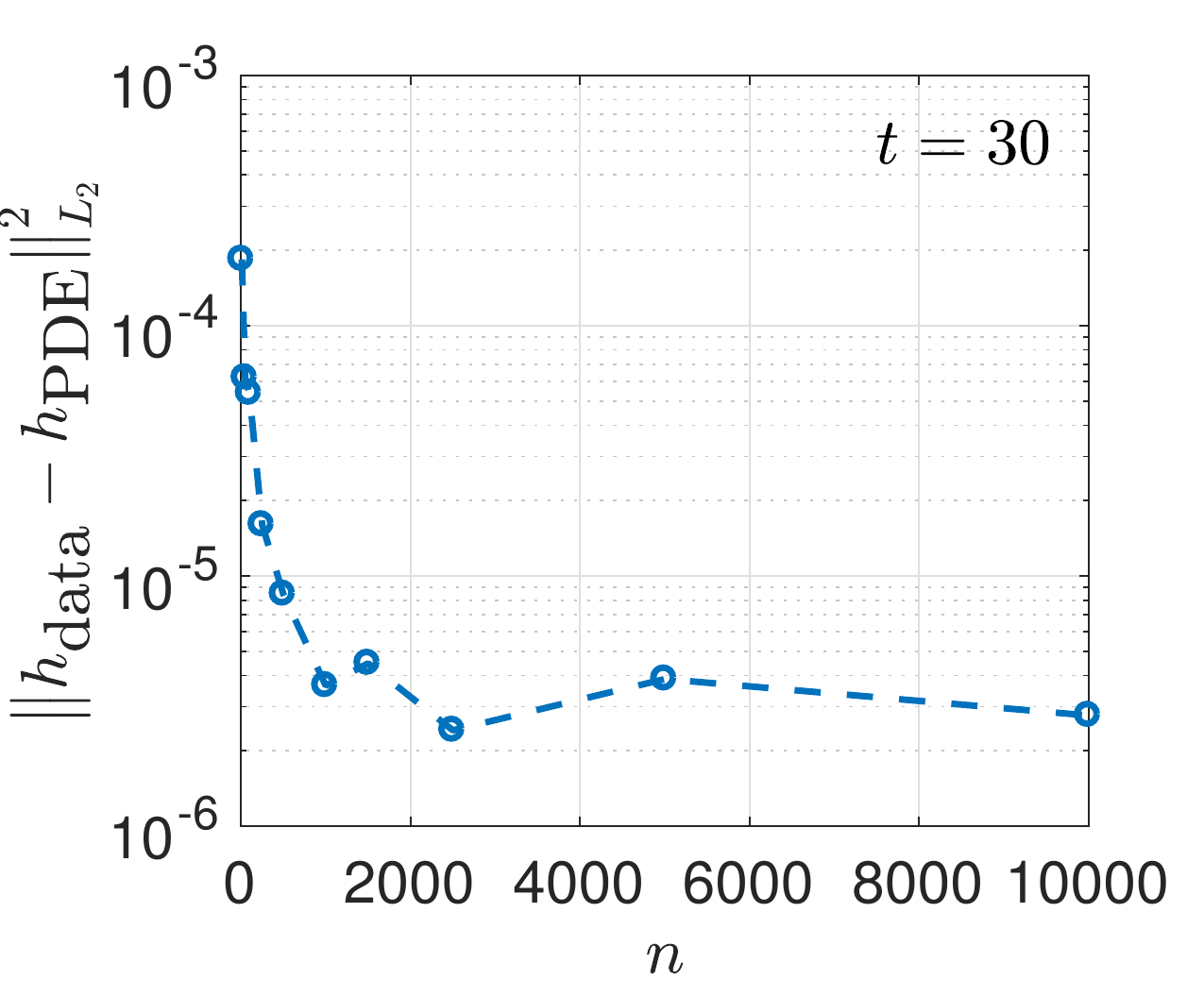}
}
\caption{Drug resistant malaria propagation model. (a) Time dependent errors 
in \eqref{hmalaria1} for a variety of sample sizes. (b) 
Error decay at $t=30$ versus the number of sample trajectories.}
\label{fig:malaria_information}
\end{figure}

\begin{figure}[t]
\centering{
	\includegraphics[width=0.45\textwidth]{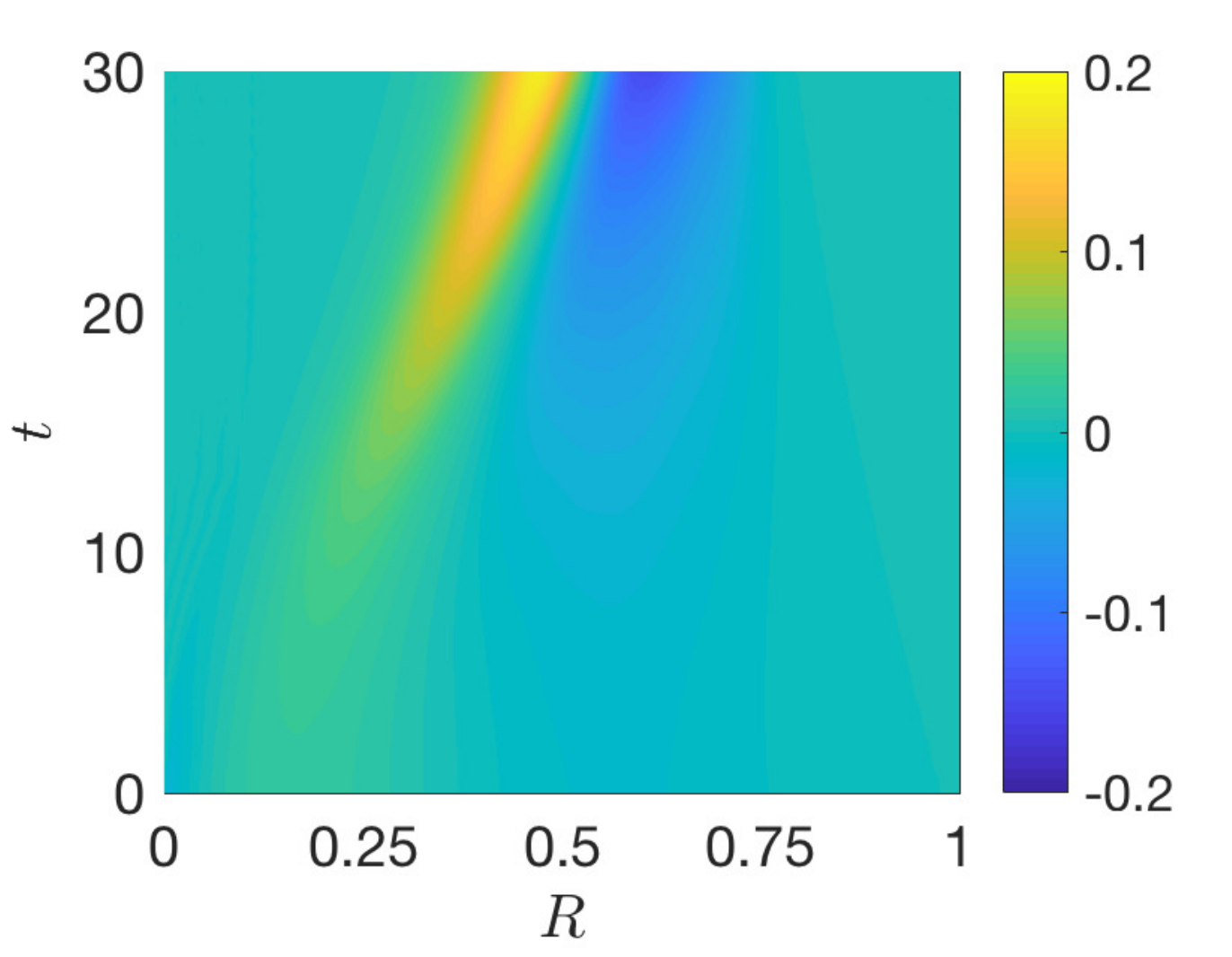}\hspace{0.2cm}
	\includegraphics[width=0.4\textwidth]{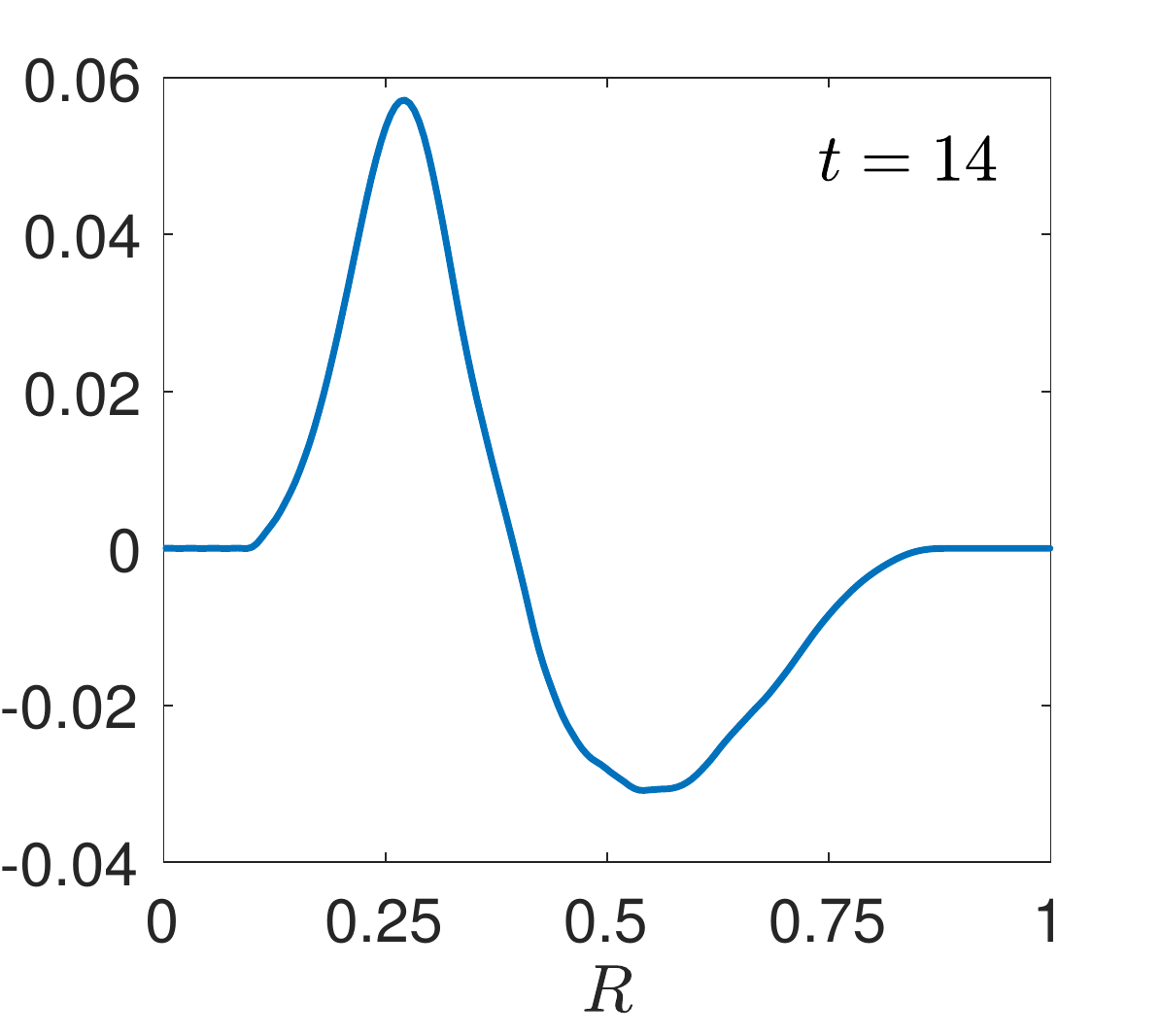}}
	
\caption{Drug resistant Malaria model 
\eqref{eq:malaria_dyn}-\eqref{eq:malaria_dyn1}. 
Mori-Zwanzig memory integral associated with the PDF of the 
phase variable representing the population of humans with temporary immunity.}
\label{fig:MalariaMZ}
\end{figure}

Suppose we are interested in a statistical description of the population 
of humans that have temporary immunity from the virus. Such population 
is represented by the phase variable $R(t)$. The evolution equation 
for the PDF of $R(t)$ is 
\begin{equation} 
\frac{\partial p(R,t)}{\partial t} = \frac{\partial }{\partial R} \left[ \left(w + \mu_h\right)R p(R,t) - h(R,t)\right], 
\label{malaria_pdf}
\end{equation}
where  
\begin{align} 
h(R, t) =\int_{-\infty}^\infty\cdots\int_{-\infty}^\infty \left[(T_s + bT_a)r +  \sigma (\xi I_a + \xi J_a + I_s + J_s) \right] p(R, I_a, \ldots, T_s) dI_a \cdots dT_s.
\label{hmalaria}
\end{align}
Clearly, $h(R,t)$ is an unclosed term which can be written as a product of 
$p(R,t)$ and the conditional expectation  of $(T_s(t)+ bT_a(t))r +  
\sigma \left(\xi I_a(t) + \xi J_a(t) + I_s(t) + J_s(t)\right)$ given $R(t)$, i.e., 
\begin{equation}
h(R,t) = p(R,t)\mathbb{E}\left[(T_s(t)+ bT_a(t))r +  
\sigma \left(\xi I_a(t) + \xi J_a(t) + I_s(t) + J_s(t)\right) | R(t)\right].
\label{hmalaria1}
\end{equation}
The evolution equation for $h(R,t)$ can be obtained by 
differentiating \eqref{hmalaria} with respect to time and using 
the Lioville equation. This yields,
\begin{align}
\frac{\partial h(R,t)}{\partial t}  =   - & \frac{\partial }{\partial R} \Bigg( p(R, t)\mathbb{E}\left[\left.\Big(r(T_s + bT_a) +  \sigma (\xi I_a \!+\! \xi J_a \!+\! I_s \!+\! J_s) \Big)^2 \right|  R(t)\right] - h(R,t)R(\omega  + \mu_h ) \Bigg) \nonumber \\
		& + p(R,t) \mathbb{E}\left[\xi \sigma \Big((1-\lambda)\beta_h M_s(t) S(t) - I_a(t) (qc + \nu + \sigma + \mu_h)\Big) + \dots \right.\nonumber\\
		&   \sigma   \Big(\lambda  \beta_h M_s(t) S(t) + \nu I_a(t) - I_s(t)(ap +\sigma + \mu_h)\Big) + \dots  \nonumber\\
		&   \xi \sigma   \Big( (1-\lambda) k \beta_h M_r(t)(S(t) + \tau T_s(t) + \tau T(t) + \tau T_a(t)) - J_a(t)(\sigma +  \nu + \mu_h) \Big) + \dots\nonumber\\
		&   \sigma   \Big(\lambda k \beta_h M_r(t)(S(t) + \tau T_s(t) + \tau T(t) + \tau T_a(t)) + \nu J_a(t) - J_s(t)(\sigma + \mu_h) \Big) + \dots \nonumber\\
		&rb     \Big(qcI_a(t) - rT_a(t) - \tau k \beta_h T_a(t) M_r(t) - \mu_h T_a(t) \Big) + \dots \nonumber\\
		& \left.  r \Big(apI_s(t) - rT_s(t) - \tau k \beta_h T_s(t) M_r(t) - \mu_h T_s(t) \Big) \Big| R(t) \right].
	\label{evhmalaria}
	\end{align}
In previous applications, we have assumed that each state variable had an 
independent initial condition. Here, we must impose the additional constraints 
$S+I_s+I_a+J_s+J_a+T_s + T + T_a + R = 1$ and 
$M_r + M_s \leq 1$, representing conservation of humans 
and mosquitoes. In doing so, we reduce the degrees of freedom 
by one, resulting in statistically dependent initial conditions. 
We set $I_a(0)$, $I_s(0)$, $J_a(0)$, $J_s(0)$, 
and $S(0)$ to be deterministic, while $R(0)$, $T_a(0)$, 
$M_s(0)$ and $M_r(0)$ are random. In particular, $M_s$ and $M_r$ 
evolve from an initial Gamma distribution with shape parameter $1/2$ 
and scale parameter $1/4$.  
In Figure \ref{fig:malariaPDF} we plot the PDF dynamics of the humans 
with temporary immunity we obtain from data-driven closure approximation, 
and compare it with an accurate benchmark PDF. 
A phase space analysis suggests that the PDF $p(R,t)$ is attracted 
by a stable node. This implies that $p(R,t)$  approaches a Dirac-delta 
distribution asymptotically in time. The information content of the 
sample trajectories of \eqref{eq:malaria_dyn}-\eqref{eq:malaria_dyn1} 
can be measured, as before, by checking the error between 
the function \eqref{hmalaria1} computed from data or from the solution 
of the hyperbolic system \eqref{malaria_pdf}-\eqref{evhmalaria}.
The results we obtain are summarized in Figure \ref{fig:malaria_information}.
Finally, the Mori-Zwanzig memory integral associated with the reduced order 
equation for the PDF $p(R,t)$ (humans with temporary immunity) 
can be computed using \eqref{MZapprox}-\eqref{M}. 
In Figure \ref{fig:MalariaMZ}, we plot the results we obtain 
with $5000$ sample trajectories.

\section{Summary}
In this paper, we developed a new data-driven method 
to compute the probability density function of 
quantities of interest in high-dimensional random 
systems. The method is based on estimating suitable 
system-dependent conditional expectations from 
data, e.g., sample trajectories or experimental data. 
We also addressed the very important question of 
whether enough useful data is being injected into the 
reduced-order PDF equation governing the quantity 
of interest for the purpose of computing an accurate 
numerical solution. To this end, we developed a new 
paradigm which allowed us to measure the information 
content of data a posteriori by solving  
systems of hyperbolic PDEs. 
We applied the proposed mathematical framework to 
the Kraichnan-Orszag three mode problem, to 
a high-dimensional nonlinear dynamical system, 
and to a drug resistant malaria propagation 
model. In all cases we found that the numerical results 
are in agreement with the theory we developed and 
they allow us to compute effectively 
the PDF of the quantity of interest.
A question we did not address in this paper 
is whether the proposed data-driven method approach 
has advantages over probability density function estimators 
purely based on data, e.g., \cite{Botev}. 
Such estimators are computationally efficient in low dimensions,
but they are agnostic about the dynamics in the phase 
space, i.e., they do not take into account the law by which the 
sample trajectories evolve in time. In principle, this opens the 
possibility to develop new classes of estimators that leverage 
on the additional {\em information source} provided by the law 
that governs the dynamics of the system. In this setting, PDF 
estimation can be formulated as a PDE-constrained optimization 
problem, with the constraint being the exact reduced-order PDF 
equation for the quantity of interest. Preliminary numerical 
results we obtained suggest that adding the evolution 
equation (PDE) in the PDF estimation process can reduce 
significantly the number of sample trajectories that are 
necessary to obtain an accurate estimation.

\vs\vs
\noindent 
{\bf Acknowledgements} 
This work was supported by DARPA grant N66001-15-2-4055 
and NSF-TRIPODS grant 81389-444168.



\end{document}